\documentclass[11pt,leqno]{amsart}
\usepackage{epsfig}
\usepackage{amssymb}
\usepackage{amscd}
\usepackage[matrix,arrow]{xy}

\newtheorem{theo}{Theorem}[section]
\newtheorem{lemma}[theo]{Lemma}
\newtheorem{defi}[theo]{Definition}
\newtheorem{prop}[theo]{Proposition}

\newtheorem{cor}[theo]{Corollary}
\newtheorem{remark}[theo]{Remark}

\newtheorem{example}[theo]{Example}
\numberwithin{equation}{section}

\def\Z{\mathbb{Z}}

\def\A{{\mathcal A}}
\def\B{{\mathcal B}}

\def\G{{\mathcal G}}

\def\D{{\mathcal{D}}}

\def\bR{{\mathbf R}}
\def\bL{{\mathbf L}}

\def\PP{{\mathbb P}}

\def\pre-tr{\operatorname{pre-tr}}
\def\h{\operatorname{h}}
\def\Hom{\operatorname{Hom}}
\def\End{\operatorname{End}}
\def\gr{\operatorname{gr}}

\newcommand{\cQ}{{\mathcal Q}}
\newcommand{\cF}{{\mathcal F}}
\newcommand{\cG}{{\mathcal G}}
\newcommand{\cO}{{\mathcal O}}
\newcommand{\cP}{{\mathcal P}}
\newcommand{\cL}{{\mathcal L}}
\newcommand{\cM}{{\mathcal M}}

\newcommand{\cD}{{\mathcal D}}

\newcommand{\cA}{{\mathcal A}}
\newcommand{\cB}{{\mathcal B}}
\newcommand{\cI}{{\mathcal I}}
\newcommand{\cC}{{\mathcal C}}
\newcommand{\cE}{{\mathcal E}}

\newcommand{\cR}{{\mathcal R}}

\newcommand{\cH}{{\mathcal H}}

\newcommand{\cl}{\operatorname{cl}}

\newcommand{\Def}{\operatorname{Def}}
\newcommand{\Perf}{\operatorname{Perf}}

\newcommand{\Ext}{\operatorname{Ext}}

\newcommand{\Sym}{\operatorname{Sym}}

\newcommand{\dgart}{\operatorname{dgart}}
\newcommand{\art}{\operatorname{art}}

\newcommand{\coDef}{\operatorname{coDef}}
\newcommand{\cart}{\operatorname{cart}}

\newcommand{\Ho}{\operatorname{Ho}}

\newcommand{\id}{\operatorname{id}}

\newcommand{\DEF}{\operatorname{DEF}}
\newcommand{\coDEF}{\operatorname{coDEF}}
\newcommand{\ev}{\operatorname{ev}}
\newcommand{\adgalg}{\operatorname{adgalg}}
\newcommand{\Set}{\operatorname{Set}}
\newcommand{\alg}{\operatorname{alg}}
\newcommand{\calg}{\operatorname{calg}}
\newcommand{\ad}{\operatorname{ad}}
\newcommand{\NGr}{\operatorname{NGr}}
\newcommand{\Gr}{\operatorname{Gr}}

\usepackage{epsf}
\usepackage{amscd}

\newcommand{\T}{\mathcal{T}}
\newcommand{\mR}{\mathcal{R}}
\newcommand{\mQ}{\mathcal{Q}}

\newcommand{\m}{\mathfrak{m}}
\newcommand{\mIm}{\mathrm{Im}}
\newcommand{\bfL}{\mathbf{L}}
\newcommand{\bfR}{\mathbf{R}}

\newcommand{\one}{\mathbf{1}}
\newcommand{\rarr}{\rightarrow}

\title[Deformation theory of objects  in homotopy and derived categories II]{Deformation theory of objects  in homotopy and derived categories II: pro-representability of the deformation functor}

\author{Alexander I.~Efimov}
\address{Department of Mechanics and Mathematics, Moscow State University, Moscow,
Russia} \email{efimov@mccme.ru}
\author{Valery A.~Lunts}
\address{Department of Mathematics, Indiana University,
Bloomington, IN 47405, USA} \email{vlunts@indiana.edu}
\author{Dmitri O.~Orlov}
\address{Steklov Mathematical Institute, 8 Gubkina St. Moscow, Russia}
\email{orlov@mi.ras.ru}


\thanks{The first named author was partially supported by grant
NSh-1983.2009.1 and by the Moebius Contest Foundation for Young
Scientists. The second named author was partially supported by the
NSA grant H98230-05-1-0050 and CRDF grant RUM1-2661-MO-05. The
third named author was partially supported by CRDF grant
RUM1-2661-MO-05, grants RFFI 08-01-00297 and NSh-9969.2006.1.}

\begin{document}

\begin{abstract} This is the second paper in a series. In part I we
developed deformation theory of objects in homotopy and derived
categories of DG categories. Here we extend these (derived)
deformation functors to an appropriate bicategory of artinian DG
algebras and prove that these extended functors are
pro-representable in a strong sense.
\end{abstract}

\maketitle

\tableofcontents

\section{Introduction}

In our paper \cite{ELOI} we developed a general deformation theory
of objects in homotopy and derived categories of DG categories. The
corresponding deformation pseudo-functors are defined on the
category of artinian DG algebras $\dgart $ and take values in the
2-category ${\bf Gpd}$ of groupoids. More precisely if $\cA$ is a DG
category and $E$ is a right DG module over $\cA$ we defined four
pseudo-functors $$\Def ^{\h}(E), \coDef ^{\h}(E), \Def (E), \coDef
(E):\dgart \to {\bf Gpd}.$$ The first two are the {\it homotopy}
deformation and co-deformation pseudo-functors, i.e. they describe
deformations (and co-deformations) of $E$ in the homotopy category
of DG $\cA ^{op}$-modules; and the last two are their {\it derived}
analogues. The pseudo-functors $\Def ^{\h}(E)$, $\coDef ^{\h}(E)$
are equivalent and depend only on the quasi-isomorphism class of the
DG algebra $\End (E)$. The derived pseudo-functors $\Def (E)$,
$\coDef (E)$ need some boundedness conditions to give the "right"
answer and in that case they are equivalent to $\Def ^{\h}(F)$ and $
\coDef ^{\h}(F)$ respectively for an appropriately chosen
h-projective or h-injective DG module $F$ which is quasi-isomorphic
to $E$ (one also needs to restrict the pseudo-functors to the
category $\dgart _-$ of negative artinian DG algebras).

In this second paper we would like to discuss the
pro-representability of these pseudo-functors. Recall that
"classically" one defines representability only for functors with
values in the category of sets (since the collection of morphisms
between two objects in a category is a set). For example, given a
moduli problem in the form of a pseudo-functor with values in the
2-category of goupoids one then composes it with the functor $\pi
_0$ to get a set valued functor, which one then tries to (pro-)
represent. This is certainly a loss of information. But in order to
represent the original pseudo-functor one needs the source category
to be a bicategory.

It turns out that there is a natural bicategory $2\text{-}\adgalg$
of augmented DG algebras. (Actually we consider two versions of this
bicategory, $2\text{-}\adgalg$ and $2^\prime\text{-}\adgalg$, but
then show that they are equivalent). We consider its full
subcategory $2\text{-}\dgart _-$ whose objects are negative artinian
DG algebras, and show that the derived deformation functors can be
naturally extended to pseudo-functors
$$\coDEF _-(E):2\text{-}\dgart _- \to {\bf Gpd},\quad \DEF _-(E):2^\prime\text{-}\dgart _- \to {\bf
Gpd}.$$ Then (under some finiteness conditions on the graded
algebra $\Ext(E,E)=H(\cC)$, where $\cC=\bR \Hom (E,E)$), we prove
pro-representability of these pseudo-functors by the DG algebra
$\hat{S}=(B\bar{A}) ^*$ which is the linear dual of the bar
construction $B\bar{A}$ of the minimal $A_{\infty}$-model of $\cC$
(Theorems \ref{pro-repr1}, \ref{pro-repr2}, \ref{pro-repr3},
\ref{pro-repr4}).

This pro-representability appears to be more "natural" for the
pseudo-functor $\coDEF _-$, because the bar complex $B\bar{A}
\otimes _{\tau _{A}}A$ is the "universal co-deformation" of $A$
considered as an $A_{\infty}$-module over $A^{op}$. The
pro-representability of the pseudo-functor $\DEF _-$ may then be
formally deduced from that of $\coDEF _-$, but we can find the
corresponding "universal deformation" (of $A$) only under an
additional assumption on $A$ (Theorem \ref{pro-repr_special}). We
also make the equivalence $\DEF_-(E)\cong
1\text{-}\Hom(\hat{S},-)$ explicit in this case (Corollary
\ref{explicit}).

These theorems describe formal deformation theory of objects in
derived categories. Our formal moduli spaces are in general
"non-commutative DG schemes". In contrast, in the paper \cite{TV}
global {\it commutative} moduli $D^{-}\text{-}$stacks of objects
in DG categories are studied. In \cite{ELOIII} we treat in detail
an example where we can construct a global moduli space of
objects.

Namely, take some vector space $V$ of dimension $n,$ and consider
the object $\cO_{\PP(W)}\in D^{b}_{coh}(\PP(V)),$ where $W\in
\Gr(m,V)(k),$ $1\leq m\leq n-1.$ The corresponding DG algebra
$\hat{S}$ satisfies the following property: $H^i(\hat{S})=0$ for
$i\ne 0,$ and for $m\ne 1$ the algebra $H^0(\hat{S})$ is
non-commutative. This suggests the existence of a non-commutative
space $\NGr(m,V)$ such that there is a $k\text{-}$point $x$
associated with each subspace $W\subset V$ of dimension $m.$ In
\cite{ELOIII} we construct these non-commutative spaces and call
them "non-commutative Grassmanians". These non-commutative
Grassmanians should be treated as true moduli spaces of objects
$\cO_{\PP(W)}\subset D^b_{coh}(\PP(V)).$ One of their properties
is the following: if $x\in \NGr(m,V)(k)$ is the point
corresponding to $W\subset V,$ then we have $\widehat{\cO_x}\cong
H^0(\hat{S}).$

We also note that the space $\NGr(\dim V-1,V),$ which can be
considered as a (dual) non-commutative projective space, is
closely related to the non-commutative projective space of
Kontsevich-Rosenberg \cite{KR}. The example of non-commutative
Grassmanians should admit a generalization to a large class of
families of objects in derived categories, for instance,
"non-commutative Jacobians".

The first part of the paper is devoted to preliminaries on
$A_{\infty}$-algebras, $A_{\infty}$-modules and
$A_{\infty}$-categories. The only non-standard point here is the
DG category of $A_{\infty}$ $A_{\cC}$-modules for an
$A_{\infty}$-algebra $A$ and a DG algebra $\cC$, and the
corresponding derived category $D_{\infty}(A_{\cC})$. We also
discuss certain functors defined by the bar complex of an
augmented $A_{\infty}$-algebra.

In the second part we introduce the Maurer-Cartan pseudo-functor
$\cM\cC(A):\dgart\to {\bf Gpd}$ for a strictly unital
$A_{\infty}$-algebra $A$. The Maurer-Cartan groupoid
$\cM\cC_{\cR}(A)$ can be described by means of some
$A_{\infty}$-category with the same objects, which are solutions
of the generalized Maurer-Cartan equation (Section
\ref{the_definition}). We develop the obstruction theory for the
Maurer-Cartan pseudo-functor (Proposition \ref{obstr}). Finally,
we show the invariance of (quasi-) equivalence classes of the
constructed $A_{\infty}$-categories and Maurer-Cartan
pseudo-functors under the quasi-isomorphisms of
$A_{\infty}$-algebras (Theorems \ref{invar_ainf}, \ref{invar_MC}).

In the third part we define the bicategories $2\text{-}\adgalg$
and $2^\prime\text{-}\adgalg$ and the pseudo-functors $\coDEF _-$
and $\DEF _-$ and discuss their relations. We also obtain here
some results on the equivalences between the homotopy and derived
(co-)deformation functors (Lemma \ref{codef=codef^h}, Theorem
\ref{def=def^h}).

In the fourth part we prove the pro-representability theorems.

We freely use the notation and results of \cite{ELOI}. The reference
to \cite{ELOI} appears in the form I, Theorem ... . As in
\cite{ELOI} our basic reference for bicategories is \cite{Be}.

\part{$A_{\infty}$-structures and the bar complex}

\section{Coalgebras}

\subsection{Coalgebras and comodules}
\label{coalg_comod}

We will consider DG coalgebras. For a DG coalgebra $\cG$ we denote
by $\cG ^{\gr}$ the corresponding graded coalgebra obtained from
$\cG$ by forgetting the differential. Recall that if $\cG$ is a DG
coalgebra, then its graded dual $\cG ^*$ is naturally a DG
algebra. Also given a finite dimensional DG algebra $\cB$ its dual
$\cB ^*$ is a DG coalgebra.

A morphism of DG coalgebras $k\to \cG$ (resp. $\cG \to k$) is called
a co-augmentation (resp. a co-unit) of $\cG$ if it satisfies some
obvious compatibility condition. We denote by $\overline{\cG}$ the
cokernel of the co-augmentation map.

Denote by $\overline{\cG} _{[n]}$ the kernel of the $n$-th iterate
of the co-multiplication map $\Delta ^n:\overline{\cG} \to
\overline{\cG} ^{\otimes n}$. The DG coalgebra $\cG$ is called {\it
co-complete} if
$$\overline{\cG}=\bigcup_{n\geq 2}\overline{\cG} _{[n]}.$$

A $\cG$-{\it comodule} means a left DG comodule over $\cG$.

 A  $\cG^{\gr}$-comodule is {\it cofree} if it is isomorphic to
$\cG\otimes V$ with the obvious comodule structure for some graded
vector space $V$.

Denote by $\cG ^{op}$ the DG coalgebra with the opposite
co-multiplication.

Let $g:\cH \to \cG$ be a homomophism of DG coalgebras. Then $\cH$
is a DG $\cG$-comodule with the co-action $g\otimes 1\cdot \Delta
_{\cH}:\cH\to \cG \otimes \cH$ and a DG $\cG ^{op}$-comodule with
the co-action  $1\otimes g\cdot \Delta _{\cH}:\cH\to \cH \otimes
\cG$.

Let $M$ and $N$ be a right and left DG $\cG$-comodules respectively.
Their cotensor product $M\square _{\cG}N$  is defined as the kernel
of the map
$$\Delta _M\otimes 1-1\otimes \Delta _N:M\otimes N\to M\otimes \cG
\otimes N,$$ where $\Delta _M:M\to M\otimes \cG$ and $\Delta _N:N\to
\cG \otimes N$ are the co-action maps.

A DG coalgebra $\cG$ is a left and right DG comodule over itself.
Given a DG $\cG$-comodule $M$ the co-action morphism $M\to \cG
\otimes M$ induces an isomorphism $M=\cG\square _{\cG} M$.
Similarly for DG $\cG ^{op}$-modules.

\begin{defi}\label{art_coalg} The dual $\cR ^*$ of an artinian DG algebra $\cR$ is
called an {\it artinian} DG coalgebra.
\end{defi}

Given an artinian DG algebra $\cR$, its augmentation $\cR \to k$
induces the co-augmentation $k\to \cR ^*$ and  its unit $k\to \cR$
induces the co-unit $\cR ^*\to k$.

\subsection{From comodules to modules}
\label{comod_to_mod}

If   $P$ is a DG comodule over a DG coalgebra $\cG$, then $P$ is
naturally a DG module over the DG algebra $(\cG ^*)^{op}$. Namely,
the $(\cG ^*)^{op}$-module structure is defined as the composition
$$P\otimes \cG^*\stackrel{\Delta_P\otimes 1}{\longrightarrow}\cG \otimes P\otimes \cG ^*
\stackrel{T\otimes 1}{\longrightarrow}P\otimes \cG \otimes \cG
^*\stackrel{1 \otimes \ev}{\longrightarrow}P,$$ where $T:\cG\otimes
P\to P\otimes \cG$ is the transposition map.

 Similarly, if $Q$ is
a DG  $\cG ^{op}$-comodule, then $Q$ is a DG module over $\cG ^*$.

Let $P$ and $Q$ be a left and right DG $\cG$-comodules respectively.
Then $P\otimes Q$ is a  DG $\cG ^*$-bimodule, i.e. a DG $\cG^*
\otimes \cG ^{*0}$-module by the above construction. Note that its
center
$$Z(P\otimes Q):=\{ x\in P\otimes Q\ \vert \ ax=(-1)^{\bar{a}\bar{x}}xa
\ \  \text{for all} \ \  a\in \cG ^*\}$$ is isomorphic to the
cotensor product $Q\square _{\cG}P$.

\section{Preliminaries on $A_{\infty}$-algebras, $A_{\infty}$-categories and $A_{\infty}$-modules}
\label{prelim}

\subsection{$A_{\infty}$-algebras and $A_{\infty}$-modules}
\label{ainf}

The basic reference for $A_{\infty}$-structures is \cite{Le-Ha}.

Let $A=\bigoplus_{n\in \Z} A^n$ be a $\Z$-graded $k$-vector space.
Put $BA=T(A[1])=\bigoplus_{n\geq 0} A[1]^{\otimes n}$. Then the
graded vector space $BA$ has natural structure of a graded
coalgebra with counit: $$\Delta(a_1\otimes \dots\otimes
a_n)=\sum\limits_{m=0}^n (a_1\otimes \dots\otimes a_m)\otimes
(a_{m+1}\otimes \dots\otimes a_n),$$
$$\varepsilon(a_1\otimes \dots\otimes a_n)=\begin{cases}0 &\text{for }n\geq 1;\\
1 &\text{for }n=0.\end{cases}$$

Here we put $a_1\otimes \dots\otimes a_n=1$ for $n=0$. Put also
$\overline{BA}=BA/k$. Then $\overline{BA}$ is also a graded
coalgebra, but it is non-counital. The most effective way to
define the notion of a $\Z$-graded (non-unital)
$A_{\infty}$-algebra is the following:

\begin{defi}\label{ainfalg} A structure of a (non-unital) $A_{\infty}$-algebra on $\Z$-graded vector space $A$
is a coderivation $b:\overline{BA}\rarr \overline{BA}$ of degree
$1$ such that $b^2=0$, i.e. a structure of a DG coalgebra on the
graded coalgebra $\overline{BA}$.\end{defi}

Such a coderivation is equivalent to a sequence of maps
$b_n=b_n^{A}:A[1]^{\otimes n}\rarr A[1]$, $n\geq 1$, of degree $1$
satisfying for each $n\geq 1$ the following identity:

\begin{equation}\label{axioms_b} \sum\limits_{r+s+t=n}b_{r+1+t}(\one^{\otimes r}\otimes
b_s\otimes \one^{\otimes t})=0.\end{equation}

Note that the coderiveation $b:\overline{BA}\to \overline{BA}$
naturally extends to a coderivation $b:BA\rarr BA$ (which we
denote by the same letter), thus $BA$ also becomes a DG coalgebra,
and $\varepsilon\cdot b=0$. Thus, its dual $\hat{S}=(BA)^*$ is
naturally a DG algebra.

Let $s:A\rarr A[1]$ be the translation map. Identify $A^{\otimes
n}$ with $A[1]^{\otimes n}$ via the map $s^{\otimes n}$, and $A$
with $A[1]$ via the map $s$. Let $m_n=m_n^A:A^{\otimes n}\rarr A$
be the maps corresponding to $b_n$. Then $m_n$ has degree $(2-n)$
and this sequence of maps satisfies for each $n\geq 1$ the
following identity:

\begin{equation}\label{axioms_m}\sum\limits_{r+s+t=n}(-1)^{r+st} m_{r+1+t}(\one^{\otimes r}\otimes
m_s\otimes \one^{\otimes t})=0.\end{equation} In particular, $m_1$
is a differential on $A$, hence $A$ is a complex. Further, $m_2$
is a morphisms of complexes and it is associative up to homotopy
given by $m_3$. Thus, the cohomology $H(A)$ is naturally a
(possibly non-unital) graded algebra. Further, if $m_n=0$ for
$n\geq 3$, then $A$ is a (possibly non-unital) DG algebra.

Let $A_1$, $A_2$ be (non-unital) $A_{\infty}$-algebras. The most
effective way to define the notion of an $A_{\infty}$-morphism
between them is the following:

\begin{defi}\label{ainfmorph} An $A_{\infty}$-morphism $f:A_1\rarr A_2$ is a (counital) homomorphism of DG coalgebras
$f:BA_1\rarr BA_2$ (which we denote by the same letter).\end{defi}

Thus, the assignment $A\mapsto BA$ is the full embedding of the
category of (non-unital) $A_{\infty}$-algebras and
$A_{\infty}$-morphisms to the category of counital DG coalgebras.

An $A_{\infty}$-morphism $f:A_1\rarr A_2$ is equivalent to a
sequence of maps $\widetilde{f_n}:A[1]^{\otimes n}\rarr A[1]$,
$n\geq 1$, of degree zero satisfying for each $n\geq 1$ the
following identity:

\begin{equation}\label{morphisms_1}\sum\limits_{i_1+\dots+i_s=n}
b_s^{A_2}(\widetilde{f_{i_1}}\otimes\dots\otimes
\widetilde{f_{i_s}})=\sum\limits_{r+s+t}
\widetilde{f_{r+1+t}}(\one^{\otimes r}\otimes b_s^{A_1}\otimes
\one^{\otimes t}).\end{equation}

Let $f_n:A_1^{\otimes n}\rarr A_2$ be the maps corresponding to
$\widetilde{f_n}$ with respect to our identifications. Then $f_n$
has degree $(1-n)$ and this sequence of maps satisfies for each
$n\geq 1$ the following identity:

\begin{equation}\label{morphisms_2}\sum\limits_{i_1+\dots+i_s=n}
(-1)^{\epsilon(i_1,\dots,i_{s-1},s)}
m_s^{A_2}(f_{i_1}\otimes\dots\otimes
f_{i_s})=\sum\limits_{r+s+t}(-1)^{r+st+s} f_{r+1+t}(\one^{\otimes
r}\otimes m_s^{A_1}\otimes \one^{\otimes t}),\end{equation} where
$\epsilon(i_1,\dots,i_{s-1},s)=(s-1)i_1+\dots+i_{s-1}+\frac{s(s+1)}2$.

In particular, $f_1$ is a morphism of complexes, and
$H(f):H(A_1)\rarr H(A_2)$ is a morphism of (non-unital) graded
associative algebras. An $A_{\infty}$-morphism $f$ is called
quasi-isomorphism if $f_1$ is a quasi-isomorphism of complexes.

Further, we are going to define the DG category
$A\text{-mod}_{\infty}$ of $A_{\infty}$ $A$-modules for an
$A_{\infty}$-algebra $A$.

\begin{defi}\label{ainfmod} A structure of an $A_{\infty}$-module over $A$ on
the graded vector space $M$ is a differential $b^{M}:BA\otimes
M[1]\rarr BA\otimes M[1]$ of degree $1$, which defines a structure
of a DG $BA$-comodule on the graded cofree $(BA)^{gr}$-comodule
$BA\otimes M[1]$.
\end{defi}

Such a structure is equivalent to a sequence of maps
$b_n=b_n^{M}:A[1]^{\otimes (n-1)}\otimes M[1]\rarr M[1]$, $n\geq
1$, of degree $1$, satisfying for each $n\geq 1$ the identity
(\ref{axioms_b}), where $b_i$ is interpreted as $b_i^A$ or
$b_i^M$, according to the type of its arguments. It is also
equivalent to the sequence of maps $m_n=m_n^M:A^{\otimes
(n-1)}\otimes M\rarr M$, $n\geq 1$, of degree $(2-n)$ satisfying
for each $n\geq 1$ the identity (\ref{axioms_m}), where $m_i$ is
interpreted as $m_i^A$ or $m_i^M$, according to the type of its
arguments. In particular, $(m_1^M)^2=0$, hence $M$ is a complex.
Again, $m_2^M$ is a morphism of complexes and $m_2$ is associative
up to a homotopy given by $m_3^M$. Thus, $H(M)$ is naturally a
graded $H(A)$-module.

If $M$ and $N$ are $A_{\infty}$ $A$-modules, then we put
$$\Hom_{A\text{-mod}_{\infty}}(M,N):=\Hom_{BA\text{-comod}}(BA\otimes M[1], BA\otimes
N[1]).$$ More explicitly,
$$\Hom_{A\text{-mod}_{\infty}}^n(M,N)=\prod\limits_{m\geq 1}\Hom_k^n(A[1]^{\otimes (m-1)}\otimes M[1],N[1]),$$
and for $\phi=(\phi_m)\in \Hom_{A\text{-mod}_{\infty}}^n(M,N)$ one
has
\begin{equation}\label{differential}(d\phi)_m=\sum\limits_{1\leq i\leq m} b_{m-i+1}^N(\one^{\otimes
(l-i)}\otimes \phi_i)-(-1)^n\sum\limits_{r+s+t=m}
\phi_{r+1+t}(\one^{\otimes r}\otimes b_s\otimes \one^{\otimes
t}),\end{equation} where $b_s$ in the RHS is interpreted as
$b_s^A$ or $b_s^M$, according to the type of its arguments. If
$\phi=(\phi_m)\in \Hom_{A\text{-mod}_{\infty}}(M,N)$ and
$\psi=(\psi_m)\in \Hom_{A\text{-mod}_{\infty}}(N,L)$, then

\begin{equation}\label{composition}(\psi\cdot\phi)_m=\sum\limits_{1\leq i\leq
m}\psi_{m-i+1}(\one^{\otimes (m-i)}\otimes \phi_i).\end{equation}

We will write $\Hom_{A}(M,N)$ instead of
$\Hom_{A\text{-mod}_{\infty}}(M,N)$.

The closed morphism $\phi\in \Hom^0(A\text{-mod}_{\infty})$ is
called quasi-isomorphism if $\phi_1$ is a quasi-isomorphism of
complexes.

The homotopy category $K_{\infty}(A)$ is defined as
$\Ho(A\text{-mod}_{\infty})$. It is always triangulated. It turns
out that all acyclic $A_{\infty}$ $A$-modules in $K_{\infty}(A)$
are already null-homotopic. Hence the corresponding derived
category $D_{\infty}(A)$ is the same as $K_{\infty}(A)$. However,
we will write $D_{\infty}(A)$ instead of $K_{\infty}(A)$.

Let $f:A_1\rarr A_2$ be an $A_{\infty}$-morphism. Then we have the
DG functor $f_*:A_2\text{-mod}_{\infty}\rarr
A_1\text{-mod}_{\infty}$, which we call the "restriction of
scalars". Namely, if $M\in A_2\text{-mod}_{\infty}$, then $f_*(M)$
coincides with $M$ as a graded vector space, and the differential
on $BA_1\otimes f_*(M)[1]$ coincides with the differential on
$BA_1\square_{BA_2}(BA_2\otimes M[1])$ after the natural
identification
$$BA_1\otimes
f_*(M)[1]\cong BA_1\square_{BA_2}(BA_2\otimes M[1]).$$

We also have the resulting exact functor $f_*:D_{\infty}(A_2)\rarr
D_{\infty}(A_1)$. If $f$ is a quasi-isomorphism, then the DG
functor $f_*:A_2\text{-mod}_{\infty}\rarr A_1\text{-mod}_{\infty}$
is quasi-equivalence, and hence the functor
$f_*:D_{\infty}(A_2)\rarr D_{\infty}(A_1)$ is an equivalence.

We would like also to define the $A_{\infty}$-bimodules.

\begin{defi} Let $A_1$ and $A_2$ be $A_{\infty}$-algebras. A structure of an $A_{\infty}$ $A_1\text{-}A_2$-bimodule on the
graded vector space $M$ is a differential $b^M:BA_1\otimes
M[1]\otimes BA_2\to BA_1\otimes M[1]\otimes BA_2$ which defines
the structure of a DG comodule over $BA_1\otimes (BA_2)^{op}$ on
the $(BA_1\otimes (BA_2)^{op})^{gr}$-bicomodule $BA_1\otimes
M[1]\otimes BA_2$.
\end{defi}

Such a differential is given by a sequence of maps
$$b_{i,j}:A_1[1]^{\otimes i}\otimes M[1]\otimes A_2[1]^{\otimes
j}\to M[1]$$ satisfying analogous equations. In particlular, we
have a regular $A_1\text{-}A_2$-bimodule $A_1\otimes A_2$. In the
case when $A_1=A_2$, we have a diagonal bimodule $A$. The DG
category $A_1\text{-mod-}A_2$ of $A_{\infty}$
$A_1\text{-}A_2$-bimodules is defined analogously (see also
\cite{KoSo}). Again, we define $K_{\infty}(A_1\text{-}A_2)$ as the
homotopy category $\Ho(A_1\text{-mod-}A_2)$. All acyclic
$A_{\infty}$-bimodules in $K_{\infty}(A_1\text{-}A_2)$ are
null-homotopic and hence the corresponding derived category
$D_{\infty}(A_1\text{-}A_2)$ coincides with
$K_{\infty}(A_1\text{-}A_2)$.

\subsection{Strictly unital $A_{\infty}$-algebras}
\label{unitality}

\begin{defi}\label{unit} An $A_{\infty}$-algebra is called strictly unital if there exists an element $1_A\in A$ of degree
zero satisfying the following properties:\\
(U1) $m_1(1_A)=0$;\\
(U2) $m_2(a,1_A)=m_2(1_A,a)=a$ for each $a\in A$;\\
(U3) for $n\geq 3$, $m_n(a_1,\dots,a_n)$ vanishes if at least one
of $a_i$ equals to $1_A$.

Such an element $1_A$ is called a strict unit.\end{defi}

Clearly, if a strict unit exists then it is unique. An
$A_{\infty}$-morphism $f:A_1\rarr A_2$ of strictly unital
$A_{\infty}$-algebras is called strictly unital if
$f_1(1_{A_1})=1_{A_2}$, and for $n\geq 2$ $f_n(a_1,\dots,a_n)$
vanishes if at least one of $a_i$ equals to $1_{A_1}$. Further, an
$A_{\infty}$-module $M\in A\text{-mod}_{\infty}$ is called
strictly unital if $m_2^M(1_A,m)=m$ for each $m\in M$ and for
$n\geq 3$ $m_n^M(a_1,\dots,a_{n-1},m)=0$ if at least one of $a_i$
equals to $1_A$. If $A$ is strictly unital then we denote by
$D_{\infty}^{su}(A)\subset D_{\infty}(A)$ the full subcategory
which consists of strictly unital $A_{\infty}$ $A$-modules.

Analogously, if $A_1$ and $A_2$ are strictly unital
$A_{\infty}$-algebras, then we have a notion of strictly unital
$A_{\infty}$ $A_1\text{-}A_2$-bimodules, and we define
$D_{\infty}^{su}(A_1\text{-}A_2)\subset
D_{\infty}(A_1\text{-}A_2)$ as the full subcategory which consists
of strictly unital $A_{\infty}$-bimodules.

If $\cC$ is a DG algebra then it is also a strictly unital
$A_{\infty}$-algebra with $m_n=0$ for $n\geq 3$. We have an
obvious DG functor $\cC\text{-mod}\to \cC\text{-mod}_{\infty}.$ It
induces an equivalence $$D(\cC)\stackrel{\sim}{\to}
D_{\infty}^{su}(\cC).$$

Let $A$ be an arbitrary $A_{\infty}$-algebra. Then its unitization
$A_{+}:=k\cdot 1_{+}\oplus A$, which is a strictly unital
$A_{\infty}$-algebra, is defined as follows:

$$m_n^{A_+}(a_1,\dots,a_n)=m_n^{A}(a_1,\dots,a_n)\text{ for any }a_1,\dots,a_n\in A,$$
$$m_1(1_+)=0,$$
$$m_2^{A_+}(1_+,a)=m_2^{A_+}(a,1_+)=a\text{ for each
}a\in A_{+},$$ $$m_n^{A_+}(a_1,\dots,a_n)=0\text{ if at least one
of }a_i\text{ equals to }1_+.$$

Clearly, the assignment $A\mapsto A_+$ defines faithful functor
from the category of $A_{\infty}$-algebras and
$A_{\infty}$-morphisms to the category of strictly unital
$A_{\infty}$-algebras and strictly unital $A_{\infty}$-morphisms.
Further, we have an obvious faithful DG functor
$A\text{-mod}_{\infty}\rarr A_+\text{-mod}_{\infty}$. Its image
consists of strictly unital $A_{\infty}$-modules. The induced
functor $D_{\infty}(A)\to D_{\infty}^{su}(A_+)$ is an equivalence.

We call $A_{\infty}$-algebras of the form $A_{+}$ augmented
$A_{\infty}$-algebras. We also use the notation
$A=\overline{A_{+}}$.

\begin{defi} Let $A$ be an augmented $A_{\infty}$-algebra. Its
bar-cobar construction $U(A)$, which is a DG algebra, together
with a strictly unital $A_{\infty}$ quasi-isomorphism $f_A:A\rarr
U(A)$ are defined by the following universal property. If $\B$ is
a DG algebra, and $f:A\rarr \B$ is a strictly unital
$A_{\infty}$-morphism then there exists a unique morphism of DG
algebras $\varphi:U(A)\rarr \B$ such that $f=\varphi\cdot f_A$.
\end{defi}

More explicitly, $U(A)$ equals to $T(\overline{B\bar{A}}[-1])$ as
a graded algebra, and the differential comes from the differential
and comultiplication on $\overline{B\bar{A}}$. The
$A_{\infty}$-morphism $f_A$ is the obvious one.

\subsection{Minimal models of $A_{\infty}$-algebras.}

An $A_{\infty}$-algebra $A$ is called minimal if $m_1^{A}=0$. Each
(strictly unital) $A_{\infty}$-algebra is quasi-isomorphic to the
minimal (strictly unital) $A_{\infty}$-algebra.

\begin{prop}(\cite{Le-Ha}, Corollaire 1.4.1.4, Proposition
3.2.4.1) Let $A$ be an $A_{\infty}$-algebra. There exists an
$A_{\infty}$-algebra structure on $H(A)$ such that

a) $m_1=0$ and $m_2$ is induced by $m_2^A$;

b) there exists an $A_{\infty}$-quasi-isomorphism of
$A_{\infty}$-algebras $f:H(A)\to A$ such that $f_1$ induces the
identity in cohomology.

Moreover, if $A$ is strictly unital then this
$A_{\infty}$-structure on $H(A)$ and the quasi-isomorphism can be
chosen to be strictly unital.
\end{prop}

\subsection{Perfect $A_{\infty}$-modules and $A_{\infty}$-bimodules.}

Let $A$ be a strictly unital $A_{\infty}$-algebra. The category
$\Perf(A)$ of perfect $A_{\infty}$ $A$-modules is the minimal full
thick triangulated subcategory of $D_{\infty}^{su}(A)$ which
contains $A$.

Further, if $A_1$ and $A_2$ are strictly unital
$A_{\infty}$-algebras then the category $\Perf(A_1\text{-}A_2)$ of
perfect $A_{\infty}$ $A_1\text{-}A_2$-bimodules is the minimal
full thick triangulated subcategory of
$D_{\infty}^{su}(A_1\text{-}A_2)$ which contains $A_1\otimes A_2$.

\subsection{$A_{\infty}$-categories}

The notion of an $A_{\infty}$-category is a straightforward
generalization of the notion of an $A_{\infty}$-algebra. Namely, a
non-unital $A_{\infty}$-category $\A$ is the following data:

- the class of of objects of $\A$;

- for each two objects $X_1,X_2$ the graded vector space
$\Hom(X_1,X_2)$;

- for each finite sequence of objects $X_0,X_1,\dots,X_n\in \A$,
$n\geq 1$, the map
$$m_n^{\A(X_0,\dots,X_n)}:\Hom(X_{n-1},X_n)\otimes\dots\otimes
\Hom(X_0,X_1)\rarr \Hom(X_0,X_n)$$ of degree $(2-n)$, such that
for any $Y_1,\dots,Y_m\in \A$ the graded vector space
$\bigoplus\limits_{1\leq i,j\leq m}\Hom(Y_i,Y_j)$ becomes an
$A_{\infty}$-algebra.

If $\A$ is an $A_{\infty}$-category then $\Ho(\A)$ is a
pre-category, i.e. a "category" which may not have identity
morphisms.

An element $\one_X\in \Hom(X,X)$ of degree zero is called a strict
identity morphism if it satisfies the conditions $(U1), (U2),
(U3)$ from Definition \ref{unit}, where $a$ and $a_i$ are
arbitrary morphisms such that the equalities make sense. An
$A_{\infty}$-category is called strictly unital if each object has
a strict identity morphism. If $\A$ is a strictly unital
$A_{\infty}$-category then $\Ho(\A)$ is a true category.

A (strictly unital) $A_{\infty}$-algebra can be thought of as a
(strictly unital) $A_{\infty}$-category with one object.

Let $\A_1,\A_2$ be $A_{\infty}$-categories. An
$A_{\infty}$-functor $F:\A_1\rarr \A_2$ is the following data:

- an object $F(X)\in \A_2$ for each object $X\in \A_1$;

- for each finite sequence of objects $X_0,X_1,\dots,X_n\in \A_1$,
$n\geq 1$, the map
$$F(X_0,\dots,X_n):\Hom(X_{n-1},X_n)\otimes\dots\otimes
\Hom(X_0,X_1)\rarr \Hom(F(X_0),F(X_n))$$ of degree $(1-n)$, such
that for any $Y_1,\dots,Y_m\in \A_1$ we obtain an
$A_{\infty}$-morphism between $A_{\infty}$-algebras
$$\bigoplus\limits_{1\leq i,j\leq m}\Hom(Y_i,Y_j)\rarr \bigoplus\limits_{1\leq i,j\leq m}\Hom(F(Y_i),F(Y_j)).$$

The definition of a strictly unital $A_{\infty}$-functor between
strictly unital $A_{\infty}$-categories is analogous to the
definition of a strictly unital $A_{\infty}$-morphism between
strictly unital $A_{\infty}$-algebras.

A strictly unital $A_{\infty}$-functor $F:\A_1\rarr \A_2$ between
strictly unital $A_{\infty}$-categories is called
quasi-equivalence if the following conditions hold:

- the map $F(X,Y):\Hom(X,Y)\rarr \Hom(F(X),F(Y))$ is a
quasi-isomorphism of complexes for any $X,Y\in \A_1$;

- the induced functor $\Ho(F):\Ho(\A_1)\rarr \Ho(\A_2)$ is an
equivalence.

\subsection{The tensor product of an $A_{\infty}$-algebra and a DG algebra}
\label{product}

Let $A$ be an $A_{\infty}$-algebra and $\cC$ be a DG algebra. Then
their tensor $A\otimes \cC$ is naturally an $A_{\infty}$-algebra
with the following multiplications:

$$m_1^{A\otimes \cC}=m_1^{A}\otimes \one_{\cC}+\one_{\A}\otimes d_{\cC};$$
$$m_n^{A\otimes \cC}(a_1\otimes c_1,\dots,a_n\otimes c_n)=(-1)^{\epsilon} m_n^A(a_1,\dots,a_n)\otimes (c_1\dots c_n)\text{ for }n\geq 2,$$
where $\epsilon=\sum\limits_{i<j} \bar{a_j}\bar{c_i}$ (all $a_i$
and $c_i$ are homogeneous). If $A$ is strictly unital, then
$A\otimes \cC$ is also strictly unital and $\one_{A\otimes
\cC}=\one_{A}\otimes \one_{\cC}$.

\begin{remark} The constructed tensor product is a specialization
of the complicated construction of the tensor product of
$A_{\infty}$-algebras which was first proposed in \cite{SU}. We
also remark that in the case when $A_1$ and $A_2$ are strictly
unital $A_{\infty}$-algebras, there is a canonical DG model for
$A_1\otimes A_2$:
$$A_1 "\otimes" A_2=\End_{A_1\text{-mod-}A_2^{op}}(A_1\otimes A_2),$$
see \cite{KoSo}.
\end{remark}

\subsection{The category of $A_\cC$-modules for an $A_{\infty}$-algebra $A$ and
a DG algebra $\cC$}

Let $A$ be an $A_{\infty}$-algebra and let $\cC$ be a DG algebra.
We want to define the DG category of $A_{\infty}$
$A_{\cC}$-modules which is analogue of the category of $(A\otimes
\cC)$-modules in the case when $A$ is a DG algebra.

\begin{defi}\label{A_C-modules} A structure of an $A_{\infty}$ $A_{\cC}$-module on the graded
vector space $M$ is the following data:

1) A structure of a $\cC^{gr}$-module on $M$;

2) A differential $b^{M}:BA\otimes M[1]\rarr BA\otimes M[1]$ of
degree $1$ which makes $BA\otimes M[1]$ into a DG comodule over
$BA$ and into a DG module over $\cC$.\end{defi}

If we are already given with the structure of a $\cC^{gr}$-module
on $M$ then such a differential $b^{M}$ is equivalent to the
sequence of maps $b_n=b_n^{M}:A[1]^{\otimes (n-1)}\otimes
M[1]\rarr M[1]$, $n\geq 1$, satisfying the following properties:\\
1) The maps $b_n^M$ satisfy the identities (\ref{axioms_b}) (in the same sense as for $A_{\infty}$
$A$-modules);\\
2) The differential $b_1^M$ makes $M[1]$ into a DG module over $\cC$;\\
3) The maps $b_n^M$ are $\cC^{gr}$-linear for $n\geq 2$.

Further, the corresponding maps $m_n=m_n^{M}:A^{\otimes
(n-1)}\otimes
M\rarr M$ have to satisfy the following properties:\\
1) The maps $m_n^M$ satisfy the identities (\ref{axioms_m}) (in
the same sense as for $A_{\infty}$
$A$-modules);\\
2) The differential $m_1^M$ makes $M$ into a DG module over $\cC$;\\
3) The maps $m_n^M$ are $\cC^{gr}$-linear for $n\geq 2$.

If $M,N$ are $A_{\infty}$ $A_{\cC}$-modules then we put
$$\Hom_{A_{\cC}\text{-mod}_{\infty}}(M,N):=\Hom_{BA\text{-comod}}\cap
\Hom_{\cC\text{-mod}}(BA\otimes M[1], BA\otimes N[1]).$$

More explicitly,
$$\Hom_{A_{\cC}\text{-mod}_{\infty}}^n(M,N)=\prod\limits_{m\geq
1}\Hom_{\cC^{gr}}^n(A[1]^{\otimes (m-1)}\otimes M[1],N[1]),$$ the
differential and the compositions are defined by the formulas
(\ref{differential}) and (\ref{composition}) respectively.

We will write $\Hom_{A_{\cC}}(M,N)$ instead of
$\Hom_{A_{\cC}\text{-mod}_{\infty}}(M,N)$.

Again, the homotopy category $K_{\infty}(A_{\cC})$ is defined as
$\Ho(A_{\cC}\text{-mod}_{\infty})$. The acyclic $A_{\infty}$
$A_{\cC}$-module in $K_{\infty}(A_{\cC})$ are not null-homotopic
in general, hence we define the derived category
$D_{\infty}(A_{\cC})$ as the Verdier quotient of
$K_{\infty}(A_{\cC})$ by the subcategory of acyclic $A_{\infty}$
$A_{\cC}$-modules.

\begin{remark} Notice that the structure of an $A_{\infty}$
$A_{\cC}$-module is not equivalent to the structure of an
$A_{\infty}$ $A\otimes \cC$-module. Moreover, there is a natural
DG functor $A_{\cC}\text{-mod}_{\infty}\to A_+\otimes
\cC\text{-mod}_{\infty}$ which induces an equivalence
$D_{\infty}(A_{\cC})\stackrel{\sim}{\to}
D_{\infty}^{su}(A_+\otimes \cC)$. Also, in the case when $A$ is
strictly unital, the DG functor $A_{\cC}\text{-mod}_{\infty}\to
A\otimes \cC\text{-mod}_{\infty}$ induces an equivalence
$D_{\infty}^{su}(A_{\cC})\stackrel{\sim}{\to}
D_{\infty}^{su}(A\otimes \cC)$.
\end{remark}

\begin{defi} An $A_{\infty}$ $A_{\cC}$-module $M$ is called $h$-projective (resp. $h$-injective) if for each acyclic
$N\in A_{\cC}\text{-mod}_{\infty}$ the complex
$\Hom_{A_{\cC}}(M,N)$ (resp. $\Hom_{A_{\cC}}(N,M)$) is
acyclic.\end{defi}

It turns out that an $A_{\infty}$ $A_{\cC}$-module is h-projective
(resp. h-injective) iff it is such as a DG $\cC$-module.

\begin{prop} Let $M$ be an $A_{\infty}$ $A_{\cC}$-module. Suppose
that $M$ is h-projective (resp. h-injective) as a DG $\cC$-module.
Then $M$ is also h-projective (resp. h-injective) as an
$A_{\infty}$ $A_{\cC}$-module.
\end{prop}

\begin{proof} We will prove Proposition for h-projectives. The
proof for h-injectives is analogous.

So let $M\in A_{\cC}\text{-mod}_{\infty}$ and suppose that $M$ is
h-projective as a DG $\cC$-module. Let $N$ be an acyclic
$A_{\infty}$ $A_{\cC}$-module. The complex
$K^{\cdot}=\Hom_{A_{\cC}}(M,N)$ admits a decreasing filtration by
subcomplexes $$F^pK^{\cdot}=\prod_{n\geq
p}\Hom_{\cC^{gr}}(A^{\otimes n}\otimes M,N).$$ The subquotients
$$F^pK^{\cdot}/F^{p+1}K^{\cdot}=\Hom_{\cC}(A^{\otimes
p}\otimes M,N)$$ are acyclic since the DG modules $A^{\otimes
p}\otimes M$ are h-projective. Since
$$K^{\cdot}=\lim\limits_{\leftarrow}K^{\cdot}/F^pK^{\cdot},$$
the complex $K^{\cdot}$ is also acyclic. Therefore, $M$ is
h-projective as an $A_{\infty}$ $A_{\cC}$-module.
\end{proof}

We denote by $K_{\infty}^P(A_{\cC})\subset K_{\infty}(A_{\cC})$
(resp. by $K_{\infty}^I(A_{\cC})\subset K_{\infty}(A_{\cC})$) the
full subcategory which consists of h-projective (resp.
h-injective) $A_{\infty}$ $A_{\cC}$-modules.

\begin{theo} For each $M\in A_{\cC}\text{-mod}_{\infty}$, there exist
quasi-isomorphisms $M\rarr I$, $P\rarr M$,where $I\in
A_{\cC}\text{-mod}_{\infty}$ is $h$-injective and $P\in
A_{\cC}\text{-mod}_{\infty}$ is $h$-projective. The natural
functor $K_{\infty}^P(A_{\cC})\to D_{\infty}(A_{\cC})$ (resp.
$K_{\infty}^I(A_{\cC})\to D_{\infty}(A_{\cC})$) is an
equivalence.\end{theo}

\begin{proof}  First we construct a quasi-isomorphism $pM\to M$
with h-projective $P$. Namely, let $pM$ be the total complex of
the bicomplex $$\dots\to \cC^{\otimes n}\otimes
M\stackrel{d^n}{\to} \cC^{\otimes n-1}\otimes M\to\dots\to
\cC\otimes M,$$ where $d^n$ is the bar differential. Then $pM$ is
naturally an $A_{\infty}$ $A_{\cC}$-module. A quasi-isomorphism of
complexes $pM\to M$ is a quasi-isomorphism in
$A_{\cC}\text{-mod}_{\infty}$ (with zero components
$f_n:A^{\otimes n-1}\otimes M\to M$ for $n\geq 2$). Further, $pM$
satisfies property (P) as a DG $\cC$-module (I, Definition 3.2).
Hence, $pM$ is an h-projective $A_{\cC}$-module.

The construction $M\to pM$ extends to the functor
$p:K_{\infty}(A_{\cC})\to K_{\infty}^P(A_{\cC})$ which is right
adjoint to the inclusion $K_{\infty}^P(A_{\cC})\to
K_{\infty}(A_{\cC})$. The kernel of $p$ consists of acyclic
$A_{\cC}$-modules. Thus, the functor $K_{\infty}^P(A_{\cC})\to
D_{\infty}(A_{\cC})$ is an equivalence.

Analogously, one can construct a functor $i:K_{\infty}(A_{\cC})\to
K_{\infty}^I(A_{\cC})$ which is left adjoint to the inclusion
$K_{\infty}^I(A_{\cC})\to K_{\infty}(A_{\cC})$. Thus, the functor
$K_{\infty}^I(A_{\cC})\to D_{\infty}(A_{\cC})$ is an equivalence.
Theorem is proved.
\end{proof}

Notice that if $G:K_{\infty}(A_{\cC})\to \T$ is an exact functor
between triangulated categories then we can define its left and
right derived functors $$\bL G:D_{\infty}(A_{\cC})\to \T,\quad \bR
G:D_{\infty}(A_{\cC})\to \T.$$ Namely, for each $M\in
A_{\cC}\text{-mod}_{\infty}$ choose quasi-isomorphisms $P\to M,$
$M\to I$ with h-projective $P$ and h-injective $I,$ and put $$\bL
G(M)=G(P),\quad \bR G(M)=G(I).$$

\begin{prop}\label{AequivU(A)} The derived categories $D_{\infty}(A_{\cC})$ and $D(U(A_+)\otimes
\cC)$ are naturally equivalent.
\end{prop}

\begin{proof} Indeed, the "restriction of scalars" DG functor
$$f_{A*}:(U(A_+)\otimes \cC)\text{-mod}\to
A_{\cC}\text{-mod}_{\infty}$$ admits a right adjoint DG functor
$$f_A^!:A_{\cC}\text{-mod}_{\infty}\to (U(A_+)\otimes
\cC)\text{-mod},$$ given by the formula
$$f_A^!(M)=\Hom_A(U(A_+),M).$$ For any $M\in (U(A_+)\otimes \cC)\text{-mod}$, $N\in
A_{\cC}\text{-mod}_{\infty}$, the adjunction morphisms $M\to
f_A^!f_{A*}M$, $f_{A*}f_A^!N\to N$ are quasi-isomorphisms.
Moreover, both $f_{A*}$ and $f_A^!$ preserve acyclic modules.
Thus, the induced functors $$f_{A*}:D(U(A_+)\otimes \cC)\to
D_{\infty}(A_{\cC}),\quad f_A^!:D_{\infty}(A_{\cC})\to
D(U(A)\otimes \cC)$$ are mutually inverse equivalences.
\end{proof}

\subsection{The bar complex}
\label{bar_complex}

Let $A$ be an $A_{\infty}$-algebra. The graded vector space
$BA\otimes A[1]\otimes BA$ carries a natural differential which
makes it into a DG bicomodule over $BA$. Namely, such a
differential is determined by its components
$$b_{i,j}:A[1]^{\otimes i}\otimes A[1]\otimes A[1]^{\otimes j}\to A[1],$$
and we put $b_{i,j}=b_{i+j+1}^A$.

In particular, $BA\otimes A$ is an $A_{\infty}$-module over
$A^{op}$. It is called the bar complex and is denoted by
$BA\otimes_{\tau_A} A$.

Now let $A$ be an augmented $A_{\infty}$-algebra, and put
$\hat{S}=(B\bar{A})^*$. The graded vector space $B\bar{A}\otimes
A[1]\otimes B\bar{A}$ also carries a natural differential which
makes it into a DG bicomodule over $B\bar{A}$. In particular,
$B\bar{A}\otimes A$ is an $A_{\infty}$-module over $\bar{A}^{op}$.
It is also called the bar complex and is denoted by
$B\bar{A}\otimes_{\tau_A} A$. Note that $B\bar{A}\otimes_{\tau_A}
A$ is a $B\bar{A}$-comodule, and hence is a $\hat{S}^{op}$-module.
This makes it into an object of
$\bar{A}^{op}_{\hat{S}^{op}}\text{-mod}_{\infty}$.

Analgously, we have an $A_{\infty}$ $\bar{A}_{\hat{S}}$-module
$A\otimes_{\tau_A}B\bar{A}$.

\section{Some functors defined by the bar complex}
\label{functors_by_bar}

\begin{defi}\label{admissible} An augmented $A_{\infty}$-algebra $\cC$ is called

a)
 {\it nonnegative} if $\cC ^i=0$
for $i<0$;

b) {\it connected} if $\cC ^0=k$;

c) {\it locally finite} if $\dim _k\cC ^i  <\infty$ for all $i$.

We say that $\cC$ is admissible if it satisfies  a), b), c).
\end{defi}

\subsection{The functor $\Delta$}
\label{Fun_Delta}

Fix an augmented $A_{\infty}$-algebra $\cC$. Consider the bar
construction $B\bar{\cC}$, the corresponding DG algebra
$\hat{S}=(B\bar{\cC})^*$ and the $A_{\infty}$
$\bar{\cC}^{op}_{\hat{S}^{op}}$-module $B\bar{\cC} \otimes _{\tau
_{\cC}} \cC$ (the bar complex). If $\cC$ is connected and
nonnegative, then $B\bar{\cC}$ is concentrated in nonnegative
degrees and consequently $\hat{S}$ is concentrated in nonpositive
degrees.

Let $\cB$ be a DG algebra. Denote by $D_f(\cB^{op})\subset
D(\cB^{op})$ the full triangulated subcategory consisting of DG
modules with finite dimensional cohomology.

\begin{lemma}\label{D_f=<k>} Assume that  DG algebra $\cB$ is augmented and local and
complete.  Also assume that $\cB ^i=0$ for $i>0$. Then the
category $D_f(\cB^{op})$ is the triangulated envelope of the DG
$\cB^{op}$-module $k$.
\end{lemma}

\begin{proof} Denote by $\langle k\rangle\subset D(\cB^{op})$ the
triangulated envelope of $k$.

 Let $M$ be a DG $\cB^{op}$-module with finite
dimensional cohomology. First assume that $M$ is concentrated in one
degree. Then $\dim M<\infty$. Since $\cB^{\gr}$ is a complete local
algebra the module $M$ has a filtration with subquotients isomorphic
to $k$. Thus $M\in \langle k\rangle$.

In the general case by I, Lemma 3.19 we may and will assume that
$M^i=0$ for $\vert i\vert
>>0$. Let $s$ be the least integer such that $M^s\neq 0$. The kernel
$K$ of the differential $d:M^s\to M^{s+1}$ is a DG
$\cB^{op}$-submodule. By the above argument $K\in \langle
k\rangle$. If $K\neq 0$ then by induction on the dimension of the
cohomology we obtain that $M/K \in \langle k\rangle$. Hence also
$M \in \langle k\rangle$. If $K=0$, then the DG $\cB
^{op}$-submodule $\tau _{< s+1}M$ (I, Lemma 3.19) is acyclic, and
hence $M$ is quasi-isomorphic to $\tau _{\geq s+1}M$. But we may
assume that $\tau _{\geq s+1}M \in \langle k\rangle$ by descending
induction on $s$.
\end{proof}

Choose a quasi-isomorphism of $A_{\infty}$
$\bar{\cC}^{op}_{\hat{S}^{op}}$-modules $B\bar{\cC} \otimes _{\tau
_{\cC}} \cC \to J$, where $J$ satisfies the property (I) as
$\hat{S}^{op}$-module (hence is h-injective).

Consider the contravariant DG functor $\Delta
:\hat{S}^{op}\text{-mod}\to \bar{\cC} ^{op}\text{-mod}_{\infty}$
defined by
$$\Delta (M):=\Hom _{\hat{S}^{op}}(M,J)$$
This functor extends trivially to derived categories $\Delta
:D(\hat{S}^{op})\to D_{\infty}(\bar{\cC} ^{op})$.

\begin{theo}\label{Delta} Assume that the DG algebra $\cC$ is admissible. Then

 a) The contravariant
functor $\Delta $ is full and faithful on the category
$D_f(\hat{S}^{op})$.

b) $\Delta(k)$ is isomorphic to $\cC$.
\end{theo}

\begin{proof} By Lemma \ref{D_f=<k>} the category $D_f(\hat{S}^{op})$ is the triangulated
envelope of the DG $\hat{S}^{op}$-module $k$. So for the first
statement of the theorem it suffices to prove that the map $\Delta
:\Ext _{\hat{S}^{op}}(k,k)\to \Ext _{\cC ^{op}}(\Delta (k),\Delta
(k))$ is an isomorphism. The following proposition implies the
theorem.

\begin{prop}\label{admissible_prelim} Under the assumptions of the above theorem the following holds.

a) The complex $\bR \Hom _{\hat{S}^{op}}(k,k)$ is quasi-isomorphic
to $\cC$.

b) The natural morphism of complexes $\Hom
_{\hat{S}^{op}}(k,B\bar{\cC} \otimes _{\tau _{\cC}}\cC)\to \Hom
_{\hat{S}^{op}}(k,J)$ is a quasi-isomorphism.

c) $\Delta (k)$ is quasi-isomorphic to $\cC$.

d) $\Delta :\Ext _{\hat{S}^{op}}(k,k)\to \Ext _{\cC ^{op}}(\Delta
(k),\Delta (k))$ is an anti-isomorphism.
\end{prop}

\begin{proof} a) Recall the $A_{\infty}$ $\bar{\cC}_{\hat{S}}$-module
$\cC \otimes _{\tau _{\cC}}B\bar{\cC}$ (subsection
\ref{bar_complex}). Consider the corresponding $A_{\infty}$
$\bar{\cC}^{op}_{\hat{S}^{op}}$-module $P:=\Hom _k(\cC \otimes
_{\tau _{\cC}}B\bar{\cC}, k)$. Since $\cC$ is locally finite and
bounded below and $B\bar{\cC}$ is bounded below the graded
$\hat{S}^{op}$-module $P^{\gr}$ is isomorphic to $(\hat{S}\otimes
\Hom _k(\cC ,k))^{\gr}$. Since the complex $\Hom _k(\cC ,k)$ is
bounded above and the DG algebra $\hat{S}$ is concentrated in
nonpositive degrees the DG $\hat{S}^{op}$-module $P$ has the
property (P) (and hence is h-projective). Thus $\bR \Hom
_{\hat{S}^{op}}(k,k)=\Hom _{\hat{S}^{op}}(P,k)=\Hom _k(\Hom _k(\cC
,k),k)=\cC $. This proves a).

b) Since $\Hom _{\hat{S}^{op}}(k,B\bar{\cC} \otimes _{\tau
_{\cC}}\cC)=\cC$ the assertion follows from a).

c) follows from b).

d) follows from a) and c).
\end{proof}

This proves the theorem.
\end{proof}

\begin{remark}\label{univ_codef}
Notice that for any augmented $A_{\infty}$-algebra $\cC$ we have
$\Hom _{\hat{S}^{op}}(k, B\bar{\cC} \otimes _{\tau
_{\cC}}\cC)=\cC$. Thus the $A_{\infty}$
$\bar{\cC}^{op}_{\hat{S}^{op}}$-module $B\bar{\cC} \otimes _{\tau
_{\cC}}\cC$ is a "homotopy $\hat{S}$-co-deformation" of $\cC$. The
Proposition \ref{admissible_prelim} implies that for an admissible
$\cC$ this $A_{\infty}$ $\bar{\cC}^{op}_{\hat{S}^{op}}$-module is
a "derived $\hat{S}$-co-deformation" of $\cC$.
 (Of course we have only defined co-deformations along artinian DG
algebras.)
\end{remark}

\subsection{The functor $\nabla$}
\label{Fun_Nabla}

Now we define another functor $\nabla :D(\hat{S}^{op})\to
D_{\infty}(\bar{\cC} ^{op})$, which is closely related to
$\Delta$.

Denote by $m$ the augmentation ideal of $\hat{S}$. For a DG
$\hat{S}^{op}$-module $M$ denote $M_n:=M/m^nM$ and
$$\hat{M}=\lim_{\stackrel{\longleftarrow}{n}}M_n.$$
 Fix a DG
$\hat{S}^{op}$-module $N$. Choose a quasi-isomorphism $P\to N$
with an h-projective $P$. Define
$$\nabla
(N):=\lim_{\rightarrow}\Delta(P_n)=\lim_{\rightarrow}\Hom_{\hat{S}^{op}}(P_n,J).$$

Denote by $\Perf(\hat{S}^{op})\subset D(\hat{S}^{op})$ the minimal
full triangulated subcategory which contains the DG
$\hat{S}^{op}$-module $\hat{S}$ and is closed with respect to
taking of direct summands.

\begin{theo}\label{Nabla}
Assume that the $A_{\infty}$-algebra $\cC$ is admissible and
finite dimensional. Then

a) The contravariant functor $\nabla :D(\hat{S}^{op})\to
D_{\infty}(\bar{\cC} ^{op})$ is full and faithful on the
subcategory $\Perf(\hat{S}^{op})$.

b) $\nabla (\hat{S})$ is isomorphic to $k$.
\end{theo}

\begin{proof}
Denote by $m\subset \hat{S}^{op}$ the maximal ideal and put
$S_n:=\hat{S}^{op}/m^n\hat{S}^{op}$. Since the
$A_{\infty}$-algebra $\cC$ is finite dimensional $S_n$ is also
finite dimensional for all $n$. We need a few lemmas.

\begin{lemma}\label{hom=rhom1} Let $K$ be a DG $\hat{S}^{op}$-module such that $\dim
_kK<\infty$. Then the natural morphism of complexes
$$\Hom _{\hat{S}^{op}}(K,B\cC \otimes _{\tau _{\cC}}\cC)\to \Hom _{\hat{S}^{op}}(K,J)$$
is a quasi-isomorphism.
\end{lemma}

\begin{proof} Notice that since the algebra $\hat{S}$ is local,
every element $x\in m$ acts on $K$ as a nilpotent operator. It
follows that $m^nK=0$ for $n>>0$. For the same reason the DG
$\hat{S}^{op}$-module $K$ has a filtration with subquotients
isomorphic to $k$. Thus we may prove the assertion by induction on
$\dim K$. If $K=k$, then this is part b) of Proposition
\ref{admissible_prelim}. Otherwise we can find a short exact
sequence of DG $\hat{S}^{op}$-modules
$$0\to M\to K\to N\to 0,$$
such that $\dim M,\dim N <\dim K$.

\medskip

\noindent{\bf Sublemma.} {\it The sequence of complexes
$$0\to \Hom _{\hat{S}^{op}}(N, B\cC \otimes _{\tau _{\cC}}\cC)\to \Hom _{\hat{S}^{op}}(K, B\cC
\otimes _{\tau _{\cC}} \cC) \to \Hom _{\hat{S}^{op}}(M, B\cC
\otimes _{\tau _{\cC}}\cC) \to 0$$ is  exact.}

\medskip

\begin{proof} We only need to prove the surjectivity of the map
$$\Hom _{\hat{S}^{op}}(K, B\cC \otimes _{\tau _{\cC}}
\cC) \to \Hom _{\hat{S}^{op}}(M, B\cC \otimes _{\tau
_{\cC}}\cC).$$

  Let
$n>>0$ be such that $m^nK=m^nM=0$. Let ${}_n(B\cC \otimes _{\tau
_{\cC}} \cC)\subset (B\cC \otimes _{\tau _{\cC}}\cC)$ denote the
DG $\hat{S}^{op}$-submodule consisting of elements $x$ such that
$m^nx=0$. Then ${}_n(B\cC \otimes _{\tau _{\cC}}\cC)$
 is a DG $S_n$-module and $\Hom _{\hat{S}^{op}}(K, B\cC \otimes _{\tau _{\cC}}
 \cC)=\Hom_{S_n}(K,{}_n(B\cC \otimes _{\tau _{\cC}}\cC))$ and similarly for $M$.

Note that ${}_n(B\cC \otimes _{\tau _{\cC}}\cC)$
 as a graded $S_n$-module is isomorphic to $S^*_n\otimes \cC$,
 hence is a finite direct sum of shifted copies of the injective graded module
 $S_n ^*$.  Hence the above map of complexes is surjective.
 \end{proof}

 Now we can prove the lemma.

 Consider the commutative diagram of complexes
 $$\begin{array}{ccccccccc} 0 & \to & \Hom _{\hat{S}^{op}}(N, B\cC \otimes _{\tau _{\cC}}\cC) & \to &
 \Hom _{\hat{S}^{op}}(K, B\cC \otimes _{\tau _{\cC}}
\cC) & \to & \Hom _{\hat{S}^{op}}(M, B\cC \otimes _{\tau _{\cC}}\cC) &  \to &  0\\
 & & \downarrow \alpha & & \downarrow \beta & & \downarrow \gamma &
 & \\
 0 & \to &  \Hom _{\hat{S}^{op}}(N, J) & \to & \Hom _{\hat{S}^{op}}(K, J) & \to & \Hom _{\hat{S}^{op}}(M, J) & \to &
 0,
 \end{array}$$
where the bottom row is exact since $J^{gr}$ is an injective
graded $\hat{S}^{op}$-module (because $J$ satisfies property (I)).
By the induction assumption $\alpha $ and $\gamma $ are
quasi-isomorphisms. Hence also $\beta $ is such.
\end{proof}

We are ready to prove the theorem.

It follows from Lemma \ref{hom=rhom1} that $\nabla(\hat{S})$ is
quasi-isomorphic to
$$\lim_{\rightarrow}\Hom
_{\hat{S}^{op}}(S_n,B\bar{\cC}\otimes _{\tau
_{\cC}}\cC)=\lim_{\rightarrow}\Hom
_{S_n}(S_n,{}_n(B\bar{\cC}\otimes _{\tau _{\cC}}
\cC))=\lim_{\rightarrow}({}_n(B\bar{\cC}\otimes _{\tau _{\cC}}
\cC))=B\bar{\cC}\otimes _{\tau _{\cC}} \cC.$$ This proves the
second assertion. The first one follows from the next lemma.

\begin{lemma}\label{end(k)=S} For any augmented $A_{\infty}$-algebra $\cC$
the complex $\Hom _{\bar{\cC} ^{op}}(k,k)$ is quasi-isomorphic to
$\hat{S}^{op}$.
\end{lemma}

\begin{proof} This follows straightforwardly from the definition of the DG
category of $A_{\infty}$ $\bar{\cC}^{op}$-modules.\end{proof}

%Note that the DG $\cC ^{op}$-module $B\cC \otimes _{\tau
%_{\cC}}\cC$ has the property (P). Hence
%$$\bR\Hom _{\cC ^{op}}(k,k)=\Hom _{\cC ^{op}}(B\cC \otimes _{\tau _{\cC}} \cC ,
%B\cC \otimes _{\tau _{\cC}} \cC).$$ We will define a homomorphism of
%DG algebras
%$$\theta :\hat{S}^{op}=\Hom _k(B\cC,k)\to \Hom _{\cC ^{op}}(B\cC \otimes _{\tau _{\cC}} \cC ,
%B\cC \otimes _{\tau _{\cC}} \cC)$$ and will prove that it is a
%quasi-isomorphism. Put $\theta (f)(a\otimes
%b)=f(a_{(1)})a_{(2)}\otimes b,$ where $\Delta (a)=a_{(1)}\otimes
%a_{(2)}\in B\cC \otimes B\cC$ is the comultiplication. One checks
%that $\theta $ is a homomorphism of DG algebras. It remains to prove
%that $\theta $ is a quasi-isomorphism of complexes.

%Notice the isomorphism of complexes
%$$\Hom _{\cC ^{op}}(B\cC \otimes _{\tau _{\cC}} \cC ,
%B\cC \otimes _{\tau _{\cC}} \cC)= \Hom _k(B\cC ,B\cC \otimes _{\tau
%_{\cC}}\cC).$$

%Denote by $\eta :B\cC\to k$ the counit (the unit in
%$\hat{S}^{op}$) and by $\epsilon :\cC \to k$ the augmentation.

%Consider the morphism of complexes
%$$\delta :\Hom _k(B\cC ,B\cC \otimes _{\tau _{\cC}}\cC)\to \Hom
%_k(B\cC ,k)$$ defined by $\delta (g)=\eta \otimes \epsilon \cdot g$.

%Then we have $\delta \cdot \theta =\id _{\hat{S}^{op}}$.

%On the other hand the inclusion $k\hookrightarrow B\cC \otimes
%_{\tau _{\cC}}\cC$ induces an embedding of complexes
%$$\gamma :\Hom _k(B\cC ,k)\to \Hom _k(B\cC ,B\cC \otimes _{\tau
%_{\cC}}\cC)$$ so that the composition $\gamma \cdot \delta$ is a
%quasi-isomorphism. Hence $\theta $ and $\delta$ are
%quasi-isomorphisms.

This proves the theorem.
\end{proof}

\subsection{The functor $\Psi$}
\label{Fun_Psi}

Finally consider the covariant functor $\Psi :D(\hat{S})\to
D_{\infty}(\bar{\cC} ^{op})$ defined by
$$\Psi (M):=(B\bar{\cC} \otimes _{\tau _{\cC}} \cC )\stackrel{\bL}{\otimes
}_{\hat{S}}M.$$

\begin{theo}\label{Psi} For any augmented $A_{\infty}$-algebra $\cC$ the following holds.

a) The functor $\Psi$ is full and faithful on the subcategory
$\Perf(\hat{S})$.

b) $\Psi(\hat{S})=k.$
\end{theo}

\begin{proof} b) is obvious and a) follows from Lemma \ref{end(k)=S} above.
\end{proof}

\part{Maurer-Cartan pseudo-functor for $A_{\infty}$-algebras}

\section{The definition}
\label{the_definition}

Let $A$ be a strictly unital $A_{\infty}$-algebra, and $\mR$ be an
artinian DG algebra with the maximal ideal $\m$. Recall that
$A\otimes \mR$ is naturally a strictly unital $A_{\infty}$-algebra
(see subsection \ref{product}). We define the set $MC(A\otimes
\m)$ as the set of $\alpha\in (A\otimes \m)^1$ such that the
generalized Maurer-Cartan equation holds:

\begin{equation}\label{MC} \sum\limits_{n\geq 1}(-1)^{\frac{n(n+1)}2}
m_n(\alpha,\dots,\alpha)=0.\end{equation}

This equation is well defined since $\m\subset\mR$ is nilpotent
ideal. Below for convenience we will write $\alpha^n$ instead of
$\underbrace{\alpha,\dots,\alpha}_{n}$.

There is a natural $A_{\infty}$-category
$\cM\cC_{\infty}^{\mR}(A)$ with the set of objects $MC(A\otimes
\m)$. Namely, for $\alpha_1,\alpha_2\in MC(A\otimes \m)$ we define
$$\Hom_{\cM\cC_{\infty}^{\mR}(A)}(\alpha_1,\alpha_2):=(A\otimes \mR)^{gr}$$ as a
graded vector space. Further, for
$\alpha_0,\alpha_1,\dots,\alpha_m\in MC(A\otimes \m)$ and for
homogeneous
$x_1\in\Hom(\alpha_0,\alpha_1),\dots,x_n\in\Hom(\alpha_{n-1},\alpha_n)$
we define
$$m_n^{\cM\cC_{\infty}^{\mR}(A)(\alpha_0,\dots,\alpha_n)}(x_n,\dots,x_1)=\sum\limits_{i_0,\dots,i_n\geq
0}(-1)^{\epsilon} m_{n+i_0+\dots+i_n}^{A\otimes
\mR}(\alpha_n^{i_n},x_n,\alpha_{n-1}^{i_{n-1}},\dots,\alpha_1^{i_1},x_1,\alpha_0^{i_0}),$$
where $$\epsilon=\sum\limits_{n\geq k>j\geq
0}(\bar{x_k}+i_k)i_j+\sum\limits_{k=0}^n
\frac{i_k(i_k+1)}2+\sum\limits_{k=1}^n ki_k.$$ One checks without
difficulties that this indeed defines an $A_{\infty}$-category and
that $\one\in (A\otimes \mR)^{gr}=\Hom(\alpha,\alpha)$ is a strict
identity for each $\alpha\in \cM\cC_{\infty}^{\mR}(A)$. Below we
will write $m_n^{\alpha_0,\dots,\alpha_n}$ instead of
$m_n^{\cM\cC_{\infty}^{\mR}(A)(\alpha_0,\dots,\alpha_n)}$

\begin{remark} The Maurer-Cartan equation and the formulas for
higher multiplications are the same as in the definition of the
$A_{\infty}$-category of one-sided twisted complexes, see
\cite{Ko}. Note that in the case of one-sided twisted complexes
all the solutions of Maurer-Cartan equation are automatically
"nilpotent".
\end{remark}

Now we define the Maurer-Cartan pseudo-functor
$\cM\cC(A):\dgart\rarr {\bf Gpd}$ as follows. Let $\mR$ and $\m$
be as above. The objects of the groupoid $\cM\cC_{\mR}(A)$ are the
same as the objects of $\cM\cC_{\infty}^{\mR}(A)$. For
$\alpha,\beta\in \cM\cC_{\infty}^{\mR}(A)$, let $G(\alpha,\beta)$
be the set of elements $g\in \one+(A\otimes \m)^0$ such that
$$m_1^{\alpha,\beta}(g)=\sum\limits_{i_0,i_1\geq 0}
(-1)^{i_0i_1+\frac{i_0(i_0+1)}{2}+\frac{i_1(i_1-1)}{2}}m_{1+i_0+i_1}^{A\otimes
\cR}(\beta^{i_1},g,\alpha^{i_0})=0.$$ Then we have an obvious
action of the group $(A\otimes \m)^{-1}$ on the set
$G(\alpha,\beta)$:
$$h:g\mapsto g+m_1^{\alpha,\beta}(h)=g+\sum\limits_{i_0,i_1\geq 0}
(-1)^{i_0i_1+\frac{i_0(i_0-1)}{2}+\frac{i_1(i_1-1)}{2}}m_{1+i_0+i_1}^{A\otimes
\cR}(\beta^{i_1},g,\alpha^{i_0}).$$

We define $\Hom_{\cM\cC_{\mR}(A)}(\alpha,\beta)$ as the set of
orbits $G(\alpha,\beta)/(A\otimes \m)^{-1}$. The composition of
morphisms in $\cM\cC_{\mR}(A)$ is induced by
$m_2^{\cM\cC_{\infty}^{\mR}(A)}$. It follows from the axioms of
$A_{\infty}$-structures that we obtain a well defined category.

\begin{prop}\label{MC_R(A)-groupoid} The category $\cM\cC_{\mR}(A)$ is a groupoid.\end{prop}

\begin{proof}Let $g\in \Hom_{\cM\cC_{\mR}(A)}(\alpha,\beta)$. Prove
that it has a left inverse $g'\in
\Hom_{\cM\cC_{\mR}(A)}(\beta,\alpha)$.

Let $\tilde{g}\in G(\alpha,\beta)$ be a lift of $g$. First prove
that there exists $\tilde{g'}\in 1+(A\otimes \m)^0$ such that
\begin{equation}\label{left_inverse}m_2^{\alpha,\beta,\alpha}(\tilde{g'},\tilde{g})=1.\end{equation}
Let $n$ be the minimal positive integer such that $\m^n=0$. The
proof is by induction over $n$.

For $n=1$, there is nothing to prove.

Suppose that the induction hypothesis holds for $n=m\geq 1$. Prove
it for $n=m+1$. From the induction hypothesis it follows that
there exists $\tilde{g'}\in \Hom_{\cM\cC_{\mR}(A)}(\beta,\alpha)$
such that $m_2^{\alpha,\beta,\alpha}(\tilde{g'},\tilde{g})=1+x$,
where $x\in (A\otimes \m^{n-1})^0$. Then we obviously have
$$m_2^{\alpha,\beta,\alpha}(\tilde{g'}-x,\tilde{g})=1.$$ Thus, the
induction hypothesis is proved for $n=m+1$.The statement is
proved.

Further, take $\tilde{g'}\in 1+(A\otimes \m)^0$ such that
(\ref{left_inverse}) holds. To prove that $g$ has a left inverse
it suffices to prove that
$$m_1^{\beta,\alpha}(\tilde{g'})=0.$$
From the equality (\ref{left_inverse}), and since
$m_1^{\alpha,\beta}(g)=0$, we obtain that
$$m_2^{\alpha,\beta,\alpha}(m_1^{\beta,\alpha}(\tilde{g'}),\tilde{g})=0.$$ Suppose that
$m_1^{\beta,\alpha}(\tilde{g'})\ne 0$. Take the maximal positive
integer $m$ such that $m_1^{\beta,\alpha}(\tilde{g'})\in (A\otimes
\m^{m})^0$. Then we obviously obtain that
$m_2^{\alpha,\beta,\alpha}(m_1^{\beta,\alpha}(\tilde{g'}),\tilde{g})\in
(A\otimes \m^{m})^0\setminus (A\otimes \m^{m+1})^0$, this leads to
contradiction.

Thus, $g$ has a left inverse. Analogously, it has a right inverse,
hence $g$ is invertible. Therefore, the category $\cM\cC_{\mR}(A)$
is a groupoid.
\end{proof}

Clearly, the assignment $\mR\mapsto \cM\cC_{\mR}(A)$ defines a
pseudo-functor from $\dgart$ to {\bf Gpd}. We denote this
pseudo-functor by $\cM\cC(A)$ and call it Maurer-Cartan
pseudo-functor.

Notice that if $A$ is a DG algebra, i.e. $m_n^A=0$ for $n\geq 3$,
then $\cM\cC_{\infty}^{\mR}(A)$ is a DG category. Further, for
$\phi\in \Hom(\alpha,\beta)$ we have

$$d^{\cM\cC^{\infty}_{\mR}(A)}(x)=d^{A\otimes \cR}(x)+\beta
x-(-1)^{\bar{x}}x\alpha,$$ and the composition in
$\cM\cC^{\infty}_{\mR}(A)$ is just the product in $A\otimes \mR$.
It follows that the constructed Maurer-Cartan pseudo-functor
coincides in this case with that constructed in \cite{ELOI},
Section 5.

\begin{remark} The Maurer-Cartan groupoid $\cM\cC_{\cR}(A)$ can be extended to a
$\infty$-groupoid $\cM\cC_{\cR}^{\infty}(A)$ so that
$\cM\cC_{\cR}(A)=\pi_0(\cM\cC_{\cR}^{\infty}(A))$. Further, the
assignment $\cR\to \cM\cC_{\cR}^{\infty}(A)$ defines a pseudo-
functor $\cM\cC^{\infty}(A):\dgart\to {\bf Gpd} ^{\infty}$, where
${\bf Gpd} ^{\infty}$ is a $\infty$-category of
$\infty$-groupoids.
\end{remark}

\section{Obstruction theory}

Fix a strictly unital $A_{\infty}$-algebra $A$.

Let $\mR$ be an artinian DG algebra with the maximal ideal $\m$.
Further, let $n$ be the minimal positive integer such that
$\m^{n+1}=0$. Put $\cI=\m^n$, $\bar{\mR}=\mR/\cI$, and
$\pi:\mR\rarr \bar{\mR}$ --- the projection morphism. The next
Proposition describes the obstruction theory for lifting of
objects and morphisms along the functor
$$\pi^*:\cM\cC_{\mR}(A)\rarr \cM\cC_{\bar{\mR}}(A).$$

\begin{prop}\label{obstr} 1). There exists a map
$o_2:Ob(\cM\cC_{\bar{\mR}}(A))\rarr H^2(A\otimes \cI)$ such that
$\alpha\in \cM\cC_{\bar{\mR}}(A)$ is in the image of $\pi^*$ if
and only if $o_2(\alpha)=0$. Furthermore, if $\alpha,\beta\in
\cM\cC_{\bar{\mR}}(A)$ are isomorphic then $o_2(\alpha)=0$ iff
$o_2(\beta)=0$.

2). Let $\xi\in Ob(\cM\cC_{\bar{\mR}}(A))$ be such that the fiber
$(\pi^*)^{-1}(\xi)$ is non-empty. Then there exists a simply
transitive action of the group $Z^1(A\otimes \cI)$ on the set
$Ob((\pi^*)^{-1}(\xi))$. Let $\xi_1,\xi_2\in
Ob(\cM\cC_{\bar{\mR}}(A))$ be isomorphic objects such that both
fibers $(\pi^*)^{-1}(\xi_1)$, $(\pi^*)^{-1}(\xi_2)$ are non-empty,
and let $f:\xi_1\rarr \xi_2$ be a morphism. Take the action of
$Z^1(A\otimes \cI)$ on $Ob((\pi^*)^{-1}(\xi_2))$ as above and the
action on $Ob((\pi^*)^{-1}(\xi_1))$ which is inverse to the above
action. Then there is a (non-canonical) $Z^1(A\otimes
\cI)$-equivariant map
$$\widetilde{o_1}:Ob((\pi^*)^{-1}(\xi_1))\times Ob((\pi^*)^{-1}(\xi_2))\rarr
Z^1(A\otimes \cI),$$ such that the composition of it with the
projection $$Z^1(A\otimes \cI)\rarr H^1(A\otimes \cI),$$ which we
denote by $o_1^{f}$, is canonically defined and satisfies the
following property: for $\alpha_1\in Ob((\pi^*)^{-1}(\xi_1))$,
$\alpha_2\in Ob((\pi^*)^{-1}(\xi_2))$ there exists a morphism
$\gamma:\alpha_1\rarr\alpha_2$ such that $\pi^*(\gamma)=f$ iff
$o_1^f(\alpha_1,\alpha_2)=0$.

3) Let $\tilde{\alpha},\tilde{\beta}\in \cM\cC_{\mR}(A)$ be
objects and let $f:\alpha\rarr \beta$ be a morphism from
$\alpha=\pi^*(\tilde{\alpha})$ to $\beta=\pi^*(\tilde{\beta})$.
Suppose that the set $(\pi^*)^{-1}(f)$ of morphisms
$\tilde{f}:\tilde{\alpha}\rarr \tilde{\beta}$ such that
$\pi^*(\tilde{f})=f$ is non-empty. Then there is a simple
transitive action of the group $\mIm(H^0(A\otimes \cI)\rarr
H^0(A\otimes\m,m_1^{\alpha,\beta}))$ on the set $(\pi^*)^{-1}(f)$.
In particular, the difference map
$$o_0:(\pi^*)^{-1}(f)\times (\pi^*)^{-1}(f)\rarr \mIm(H^0(A\otimes
\cI)\rarr H^0(A\otimes\m,m_1^{\alpha,\beta}))$$ satisfies the
following property: if $\tilde{f},\tilde{f'}\in (\pi^*)^{-1}(f)$
then $\tilde{f}=\tilde{f'}$ iff
$o_0(\tilde{f},\tilde{f'})=0$.\end{prop}

\begin{proof} 1). Let $\alpha\in \cM\cC_{\bar{\mR}}(A)$. Take some $\tilde{\alpha}\in (A\otimes \m)^1$ such that
$\pi(\tilde{\alpha})=\alpha$. Then we have
$$\sum\limits_{n\ge 1}(-1)^{\frac{n(n+1)}2+1} m_n^{A\otimes \mR}(\tilde{\alpha},\dots,\tilde{\alpha})\in (A\otimes \cI)^2.$$
A straightforward applying of (\ref{axioms_m}) shows that
$$\sum\limits_{n\ge 1}(-1)^{\frac{n(n+1)}2+1} m_n^{A\otimes
\mR}(\tilde{\alpha},\dots,\tilde{\alpha})\in Z^2(A\otimes \cI).$$
Further, if $\tilde{\alpha'}\in A\otimes \m$ is another lift of
$\alpha$ then \begin{equation}\label{diff_o_1} \sum\limits_{n\ge
1}(-1)^{\frac{n(n+1)}2+1} m_n^{A\otimes
\mR}(\tilde{\alpha'},\dots,\tilde{\alpha'})-\sum\limits_{n\ge
1}(-1)^{\frac{n(n+1)}2+1} m_n^{A\otimes
\mR}(\tilde{\alpha},\dots,\tilde{\alpha})=m_1^{A\otimes
\mR}(\tilde{\alpha'}-\tilde{\alpha}).\end{equation} Hence, we
obtain the well defined element $o_2(\alpha)\in H^2(A\otimes \cI)$
and therefore the map $o_2:Ob(\cM\cC_{\bar{\mR}}(A))\rarr
H^2(A\otimes \cI)$. The first property of $o_2$ is obviously
satisfied.

Further, let $\alpha,\beta\in \cM\cC_{\bar{\mR}}(A)$, and
$f:\alpha\rarr \beta$ be a morphism. Suppose that $o_2(\alpha)=0$.
Take some $\tilde{\alpha}\in (\pi^*)^{-1}(\alpha)$. Further, take
some $\tilde{f}\in 1+(A\otimes \m)^0$ such that $\pi(\tilde{f})$
represents $f$, and $\tilde{\beta}\in (A\otimes \m)^1$ such that
$\pi(\tilde{\beta})=\beta$. We have that
$$\sum\limits_{i_0,i_1\geq 0}(-1)^{i_0i_1+\frac{i_0(i_0+1)}2+\frac{i_1(i_1+3)}2} m_{1+i_0+i_1}^{A\otimes \mR}(\tilde{\beta}^{i_1},\tilde{f},\tilde{\alpha}^{i_0})\in (A\otimes \cI)^1$$ A straightforward applying of
(\ref{axioms_m}) shows that \begin{multline*}m_1^{A\otimes
\mR}(\sum\limits_{i_0,i_1\geq
0}(-1)^{i_0i_1+\frac{i_0(i_0+1)}2+\frac{i_1(i_1+3)}2}
m_{1+i_0+i_1}^{A\otimes
\mR}(\tilde{\beta}^{i_1},\tilde{f},\tilde{\alpha}^{i_0}))=\\=
m_2^{A\otimes \mR}(\sum\limits_{n\ge 1}(-1)^{\frac{n(n+1)}2+1}
m_n^{A\otimes
\mR}(\tilde{\beta},\dots,\tilde{\beta}),\tilde{f})=\sum\limits_{n\ge
1}(-1)^{\frac{n(n+1)}2+1} m_n^{A\otimes
\mR}(\tilde{\beta},\dots,\tilde{\beta}).\end{multline*} Therefore,
$o_2(\beta)=0$. This proves 1).

2). Let $\eta\in Z^1(A\otimes \cI)$. It follows from
(\ref{diff_o_1}) that the formula $$\eta:\alpha\mapsto
\alpha+\eta$$ defines a simply transitive action of the group
$Z^1(A\otimes \cI)$ on the set $Ob((\pi^*)^{-1}(\xi))$. Let
$\xi_1$, $\xi_2$, $f$  be as in Proposition. Take some
$\tilde{f}\in 1+(A\otimes \m)^0$ such that
$\overline{\pi(\tilde{f})}=f$. Define $\widetilde{o_1}$ by the
formula
$$\widetilde{o_1}(\alpha,\beta)=o_1^{\tilde{f}}(\alpha,\beta)=m_1^{\alpha,\beta}(\tilde{f}).$$

It is easy to see that the image of $o_1^{\tilde{f}}$ lies in
$Z^1(A\otimes \cI)$ and that $o_1^{\tilde{f}}$ is $Z^1(A\otimes
\cI)$-equivariant.

If $\tilde{f'}$ is another lift of $f$, then there exists $h\in
(A\otimes \m)^{-1}$ such that
$$v=\tilde{f'}-\tilde{f}-m_1^{\alpha,\beta}(h)\in (A\otimes \cI)^0.$$
Further,
$$o_1^{\tilde{f'}}(\alpha,\beta)-o_1^{\tilde{f}}(\alpha,\beta)=m_1^{A\otimes \mR}(v),$$
hence the map $o_1^f$ is canonically defined.

Suppose that $o_1^f(\alpha,\beta)=0$ for some $\alpha\in
(\pi^*)^{-1}(\xi_1), \beta\in (\pi^*)^{-1}(\xi_2)$. Let
$\tilde{f}$ be as above. Then there exists $x\in (A\otimes \cI)^0$
such that
$$m_1^{\alpha,\beta}(\tilde{f})=m_1^{A\otimes
\mR}(x).$$ We have $\tilde{f}-x\in G(\alpha,\beta)$, and
$\pi^*(\overline{\tilde{f}-x})=f$.

Conversely, suppose that there exists a morphism $\gamma\in
\Hom_{\cM\cC_{\mR}(A)}(\alpha,\beta)$ for some $\alpha\in
(\pi^*)^{-1}(\xi_1)$, $\beta\in (\pi^*)^{-1}(\xi_2)$, such that
$\pi^*(\gamma)=f$. Let $\tilde{\gamma}\in G(\alpha,\beta)$ be a
representative of $\gamma$. Then we have
$o_1^{\tilde{\gamma}}(\alpha,\beta)=0$, hence
$o_1^{f}(\alpha,\beta)=0$. This proves 2).

3). First we define the action of the group $Z^0(A\otimes \cI)$ on
the set $(\pi^*)^{-1}(f)$ by the formula
$$\eta:\bar{\tilde{f}}\rarr \overline{\tilde{f}+\eta},$$ where $\eta\in Z^0(A\otimes \cI)$,
and $\tilde{f}\in G(\alpha,\beta)$ is such that
$\pi^*(\bar{\tilde{f}})=f$. Clearly, this is correct. Further, if
$\eta=m_1^{A\otimes \mR}(\zeta)$ for some $\zeta\in (A\otimes
\m)^{-1}$, then
$$\eta(\bar{\tilde{f}})=\overline{\tilde{f}+m_1^{\alpha,\beta}(\zeta)}=\bar{\tilde{f}}.$$ Hence, we
have an action of $\mIm(H^0(A\otimes \cI)\rarr
H^0(A\otimes\m,m_1^{\alpha,\beta}))$ on the set $(\pi^*)^{-1}(f)$.

Tautologically, this action is simple.

Prove that it is transitive. Let $\tilde{f},\tilde{f'}\in
G(\alpha,\beta)$ be such that
$\pi^*(\bar{\tilde{f}})=\pi^*(\bar{\tilde{f'}})=f$. Then, by
definition, there exists $h\in (A\otimes \m)^{-1}$ such that
$$\tilde{f'}-\tilde{f}-m_1^{\alpha,\beta}(h)\in (A\otimes \cI)^0.$$
Replacing $\tilde{f}$ by $\tilde{f}+m_1^{\alpha,\beta}(h)$, we
obtain $\tilde{f'}=\tilde{f}+\eta$, where $\eta\in (A\otimes
\cI)^0$. Since $\tilde{f},\tilde{f'}\in G(\alpha,\beta)$, we have
that $\eta\in Z^0(A\otimes \cI)$. This shows transitivity and
proves 3).

Proposition is proved.
\end{proof}

\begin{remark} One can also construct the obstruction theory for
lifting of objects and all $k$-morphisms along the
$\infty$-functor
$$\pi^*:\cM\cC_{\cR}^{\infty}(A)\to \cM\cC_{\bar{\cR}}^{\infty}(A).$$
\end{remark}

\section{Invariance Theorems}

Let $A_1,A_2$ be strictly unital $A_{\infty}$-algebras and
$f:A_1\rarr A_2$ be a strictly unital $A_{\infty}$-morphism
between them given by a sequence of maps $$f_n:A_1^{\otimes
n}\rarr A_2.$$ Further, let $\mR$ be an artinian DG algebra with
the maximal ideal $\m$.

Then we have a (strictly unital) $A_{\infty}$-functor
$$f_{\mR}^*:\cM\cC_{\infty}^{\mR}(A_1)\rarr \cM\cC_{\infty}^{\mR}(A_2)$$ defined by the
formulas $$f_{\mR}^*(\alpha)=\sum\limits_{n\geq
1}(-1)^{\frac{n(n-1)}2} f_n(\alpha,\dots,\alpha);$$
$$f_{\mR}^*(\alpha_0,\dots,\alpha_n)(x_1,\dots,x_n)=\sum\limits_{i_0,\dots,i_n\geq
0}(-1)^{\epsilon}
f_{n+i_0+\dots+i_n}(\alpha_n^{i_n},x_n,\alpha_{n-1}^{i_{n-1}},\dots,\alpha_1^{i_1},x_1,\alpha_0^{i_0}),$$
where $$\epsilon=\sum\limits_{n\geq k>j\geq
0}(\bar{x_k}+i_k)i_j+\sum\limits_{k=0}^n
\frac{i_k(i_k-1)}2+\sum\limits_{k=1}^n ki_k.$$ One checks without
difficulties that these formulas indeed define a strictly unital
$A_{\infty}$-functor.

It induces a functor $f_{\mR}^*:\cM\cC_{\mR}(A_1)\rarr
\cM\cC_{\mR}(A_2)$ and we obtain a morphism of pseudo-functors
$$f^*:\cM\cC(A_1)\rarr \cM\cC(A_2).$$

The following theorems show that for quasi-isomorphic strictly
unital $A_{\infty}$-algebras the corresponding Maurer-Cartan
$A_{\infty}$-categories (resp. Maurer-Cartan pseudo-functors) are
quasi-equivalent (resp. equivalent).

\begin{theo}\label{invar_ainf} Let $f:A_1\rarr A_2$ be a strictly unital quasi-isomorphism of strictly unital $A_{\infty}$-algebras and let
$\mR$ be an artinian DG algebra with the maximal ideal $\m$. Then
the $A_{\infty}$-functor
$$f_{\mR}^*:\cM\cC_{\infty}^{\mR}(A_1)\rarr \cM\cC_{\infty}^{\mR}(A_2)$$ is
a quasi-equivalence.
\end{theo}
\begin{proof} 1). Prove that for any $\alpha,\beta\in \cM\cC_{\infty}^{\mR}(A_1)$ the morphism of complexes
$$f_{\mR}^*(\alpha,\beta):\Hom_{\cM\cC_{\infty}^{\mR}(A_1)}(\alpha,\beta)\rarr \Hom_{\cM\cC_{\infty}^{\mR}(A_2)}(f^*(\alpha),f^*(\beta))$$ is
quasi-isomorphism. Note that both complexes have finite
filtrations by subcomplexes $A_1\otimes \m^i$ and $A_2\otimes
\m^i$. The morphism $f_{\mR}^*(\alpha,\beta)$ is compatible with
these filtrations and induces quasi-isomorphisms on the
subquotients. Hence, it is quasi-isomorphism.

2). Now we prove that the functor
$$\Ho(f^*):\Ho(\cM\cC_{\infty}^{\mR}(A_1))\rarr \Ho(\cM\cC_{\infty}^{\mR}(A_2))$$ is
an equivalence. We have already proved that it is full faithful,
hence it remains to prove that it is essentially surjective. We
will prove the stronger statement: the functor
$$f_{\mR}^*:\cM\cC_{\mR}(A_1)\rarr \cM\cC_{\mR}(A_2)$$
is essentially surjective.

Let $n$ be the minimal positive integer such that $\m^n=0$. The
proof is by induction over $n$.

For $n=1$, there is nothing to prove.

Suppose that the induction hypothesis holds for $n=m$. Prove it
for $n=m+1$. Let $\cI$, $\bar{\mR}$, $\bar{\m}$, $\pi$ be as
above. A straightforward checking shows that the following diagram
commutes:

\begin{equation} \label{o_2_comm}
\begin{CD}
Ob(\cM\cC_{\bar{\mR}}(A_1)) @>f_{\bar{\mR}}^*>> Ob(\cM\cC_{\bar{\mR}}(A_2))\\
@Vo_2VV                              @Vo_2VV\\
H^2(A_1\otimes \cI) @>\sim>> H^2(A_2\otimes \cI),
\end{CD}
\end{equation}
where the map $o_2$ is defined in Prpoosition \ref{obstr}.

Let $\alpha\in \cM\cC_{\mR}(A_2)$. By the induction hypothesis,
there exists $\beta\in \cM\cC_{\bar{\mR}}(A_1)$ such that
$f_{\bar{\mR}}^*(\beta)$ is isomorphic to $\pi^*(\alpha)$ in
$\cM\cC_{\bar{\mR}}(A_2)$. Since the diagram (\ref{o_2_comm})
commutes, we have that $o_2(\beta)=0$. Thus, by Proposition
\ref{obstr}, the fiber $(\pi^*)^{-1}(\beta)$ is nonempty. Fix some
$\tilde{\beta}\in (\pi^*)^{-1}(\beta)$. Let
$\gamma:f_{\bar{\mR}}^*(\beta)\rarr \pi^*(\alpha)$ be a morphism.
A straightforward checking shows that the following diagram
commutes:

\begin{equation} \label{o_1_comm}
\begin{CD}
Ob((\pi^*)^{-1}(\beta)) @>f_{\mR}^*>> Ob((\pi^*)^{-1}(f_{\bar{\mR}}^*(\beta)))\\
@VVo_1^{id_{\beta}}(*,\tilde{\beta})V                  @VVo_1^{\gamma}(f_{\cR}^*(*),\alpha)-o_1^{\gamma}(f_{\mR}^*(\tilde{\beta}),\alpha)V\\
H^1(A_1\otimes \cI) @>\sim>> H^1(A_2\otimes \cI),
\end{CD}
\end{equation}
where the vertical arrows are defined in Proposition \ref{obstr}.
Since the map $o_1^{id_{\beta}}(*,\tilde{\beta})$ is surjective
and the diagram (\ref{o_1_comm}) commutes, there exists an object
$\tilde{\beta'}\in Ob((\pi^*)^{-1}(\beta))$ such that
$o_1^{\gamma}(f_{\mR}^*(\tilde{\beta'}),\alpha)=0$. Then, by
Proposition \ref{obstr}, there exists a morphism
$\tilde{\gamma}:f_{\mR}^*(\tilde{\beta'})\rarr \alpha$ (such that
$\pi^*(\tilde{\gamma})=\gamma$). Therefore, the functor
$f_{\mR}^*$ is essentially surjective, and the induction
hypothesis is proved for $n=m+1$. The statement is proved.

Theorem is proved.
\end{proof}

\begin{theo}\label{invar_MC} Let $f:A_1\rarr A_2$ be a strictly unital quasi-isomorphism of strictly unital $A_{\infty}$-algebras. Then
the morphism of pseudo-functors $$f^*:\cM\cC(A_1)\rarr
\cM\cC(A_2)$$ is an equivalence.\end{theo}
\begin{proof} Fix an artinian DG algebra $\mR$ with the maximal ideal $\m$. We must prove that the
functor
$$f_{\mR}^*:\cM\cC_{\mR}(A_1)\rarr \cM\cC_{\mR}(A_2)$$ is an equivalence.

In the proof of the previous Theorem we have already shown that it
is essentially surjective. So it remains to prove that it is full
and faithful.

Let $n$ be the minimal positive integer such that $\m^n=0$. The
proof is by induction over $n$.

For $n=1$, there is nothing to prove.

Suppose that the induction hypothesis holds for $n=m\geq 1$. Prove
it for $n=m+1$.\\
{\bf Full.} Let $\alpha,\beta\in \cM\cC_{\mR}(A_1)$ and let
$\gamma:f_{\mR}^*(\alpha)\rarr f_{\mR}^*(\beta)$ be a morphism. By
induction hypothesis, there exists a morphism
$g:\pi^*(\alpha)\rarr \pi^*(\beta)$ such that
$$f_{\bar{\mR}}^*(g)=\pi^*(\gamma).$$ A straightforward checking
shows that the following diagram commutes

\begin{equation}\label{o_1_comm'}
\begin{CD}
(\pi^*)^{-1}(\pi^*(\alpha))\times (\pi^*)^{-1}(\pi^*(\beta)) @>>>
(\pi^*)^{-1}(\pi^*(f_{\bar{\mR}}^*(\alpha)))\times
(\pi^*)^{-1}(\pi^*(f_{\bar{\mR}}^*(\beta)))\\
@VVo_1^{g}V                  @VVo_1^{\pi^*(\gamma)}V\\
H^1(A_1\otimes \cI) @>\sim>> H^1(A_2\otimes \cI).
\end{CD}\end{equation}

By Proposition \ref{obstr} and since the diagram (\ref{o_1_comm'})
commutes there exists a morphism $\tilde{g}:\alpha\rarr \beta$
such that $\pi^*(\tilde{g})=g$. Further, a straightforward
checking shows that the following diagram commutes:

\begin{equation} \label{o_0_comm}
\begin{CD}
(\pi^*)^{-1}(g) @>o_0(*,\tilde{g}) >> \mIm(H^0(A_1\otimes
\cI)\rarr
H^0(A_1\otimes\m,m_1^{\alpha,\beta})) \\
@VVf_{\mR}^*V @V\sim VV\\(\pi^*)^{-1}(\pi^*(\gamma))
@>o_0(f_{\cR}^*(*),\gamma)-o_0(f_{\mR}^*(\tilde{g}),\gamma) >>
\mIm(H^0(A_2\otimes \cI)\rarr
H^0(A_2\otimes\m,m_1^{f_{\mR}^*(\alpha),f_{\mR}^*(\beta)}))
\end{CD}
\end{equation}

Since the upper arrow is surjective, there exists a morphism
$\tilde{g'}\in (\pi^*)^{-1}(g)$ such that
$$o_0(f_{\mR}^*(\tilde{g'}),\gamma)=0,$$
i.e. $f_{\mR}^*(\tilde{g'})=\gamma$. Hence, the functor
$f_{\mR}^*$ is full.\\
{\bf Faithful.} Let $\gamma_1,\gamma_2:\alpha\rarr \beta$ be two
morphisms in $\cM\cC_{\mR}(A_1)$. Suppose that
$f_{\mR}^*(\gamma_1)=f_{\mR}^*(\gamma_2).$ Then we have also
$f_{\mR}^*(\pi^*(\gamma_1))=f_{\mR}^*(\pi^*(\gamma_2)),$ hence by
induction hypothesis $\pi^*(\gamma_1)=\pi^*(\gamma_2)$. A
straightforward checking shows that the following diagram
commutes:

\begin{equation} \label{o_0_comm'}
\begin{CD}
(\pi^*)^{-1}(\pi^*(\gamma_1))\times (\pi^*)^{-1}(\pi^*(\gamma_1))
@>o_0>> \mIm(H^0(A_1\otimes \cI)\rarr
H^0(A_1\otimes\m,m_1^{\alpha,\beta})) \\
@VVf_{\mR}^*V @V\sim
VV\\(\pi^*)^{-1}(\pi^*(f_{\mR}^*(\gamma_1)))\times
(\pi^*)^{-1}(\pi^*(f_{\mR}^*(\gamma_1))) @>o_0>>
\mIm(H^0(A_2\otimes \cI)\rarr
H^0(A_2\otimes\m,m_1^{f_{\mR}^*(\alpha),f_{\mR}^*(\beta)})).
\end{CD} \end{equation}

By Proposition \ref{obstr} and since the diagram (\ref{o_0_comm'})
commutes we have that $$o_0(\gamma_1,\gamma_2)=0,$$ hence
$\gamma_1=\gamma_2$. Thus, the functor $f_{\mR}^*$ is full.

The induction hypothesis is proved for $n=m+1$. The statement is
proved.

Theorem is proved.
\end{proof}

\begin{remark} It can be proved that an $A_{\infty}$-quasi-isomorphism
$f:A_1\to A_2$ induces an equivalence of $\infty$-groupoids
$f_{\cR}^*:\cM\cC_{\cR}^{\infty}(A_1)\to
\cM\cC_{\cR}^{\infty}(A_2)$.
\end{remark}

\section{Twisting cochains} Let $\G$ be a co-augmented DG
coalgebra. Let $A$ be an arbitrary $A_{\infty}$-algebra. Then the
graded vector space $Hom_{k}(\G,A)$ has natural structure of an
$A_{\infty}$-algebra. If $\dim\G<\infty$ or $\dim A<\infty$, then
$Hom_k (\G,A)$ is canonically identified with $A\otimes \G^*$ as
an $A_{\infty}$-algebra.

Suppose that the DG coalgebra $\G$ is co-complete. The map $\tau:
\G\rarr A$ of degree $1$ is called a twisting cochain if it passes
through $\bar{\G}$ and satisfies the generalized Maurer-Cartan
equation (\ref{MC}) as an element of $A_{\infty}$-algebra
$Hom_{k}(\G,A)$. This is well defined since $\G$ is co-complete.
If $\mR$ is an artinian DG algebra and $A$ is strictly unital,
then we have natural bijection between the set of twisting
cochains $\tau:\mR^*\rarr A$ and the set $MC(A\otimes \m)$. In the
case when $A$ is augmented, the twisting cochain is called
admissible if it passes through $\bar{A}$. Tautologically,
admissible twisting cochains $\G\to A$ are in one-to-one
correspondence with twisting cochains $\G\to \bar{A}$.

\begin{prop}\label{tw_coch} Let $A$ be an $A_{\infty}$-algebra The composition $\tau_{A}:BA\rarr A$ of the natural projection
$BA\rarr A[1]$ with the shift map $A[1]\rarr A$ is the universal
twisting cochain. That is, if $\G$ is a co-augmented co-complete
DG coalgebra and $\tau:\G\rarr A$ is a twisting cochain then there
exists a unique homomorphism $g_{\tau}:\G\rarr BA$ of co-augmented
DG coalgebras, such that $\tau_{A}\cdot g_{\tau}=\tau$.

It follows that if $A$ is augmented then the composition of
$\tau_{\bar{A}}$ with the embedding $\bar{A}\hookrightarrow A$,
which we also denote by $\tau_{A}$, is the universal admissible
twisting cochain in the same sense.
\end{prop}
\begin{proof} A straightforward checking.\end{proof}

Further, if $\G$ is a co-augmented co-complete DG coalgebra, and
$\tau:\G\rarr A$ is a twisting cochain then
$$\G\otimes_{\tau} A:=\G\square_{BA}
(BA\otimes_{\tau_A} A)$$ is an object of
$A^{op}_{(\G^*)^{op}}\text{-mod}_{\infty}$.

\begin{prop}\label{invar_h1} Let $f:A_1\rarr A_2$ be an $A_{\infty}$-quasi-isomorphism of $A_{\infty}$-algebras,
$\G$ be a co-augmented co-complete DG coalgebra, and $\tau:\G\rarr
A_1$ be a twisting cochain. Then there is a natural homotopy
equivalence in $A_{1(\G^*)^{op}}^{op}\text{-mod}_{\infty}$:
$$\G\otimes_{\tau} A_1\rarr f_*(\G\otimes_{f\cdot \tau} A_2).$$
\end{prop}
\begin{proof} We have a natural homotopy equivalence of DG
bicomodules over $BA_1$:
$$BA_1\otimes A_1[1]\otimes BA_1\rarr BA_1\square_{BA_2} (BA_2\otimes A_2[1]\otimes BA_2)\square_{BA_2} BA_1.$$
Co-tensoring it on the left by $\G$, we obtain the required
homotopy equivalence.
\end{proof}

Now let $A$ be an augmented $A_{\infty}$-algebra. If $\G$ is a
co-augmented co-complete DG coalgebra and $\tau:\G\rarr A$ is an
admissible twisting cochain then
$$\G\otimes_{\tau} A:=\G\square_{B\bar{A}}
(B\bar{A}\otimes_{\tau_A} A)$$ is an object of
$\bar{A}^{op}_{(\G^*)^{op}}\text{-mod}_{\infty}$.

\begin{prop}\label{invar_h} Let $f:A_1\rarr A_2$ be an $A_{\infty}$-quasi-isomorphism of augmented
$A_{\infty}$-algebras, $\G$ be a co-augmented co-complete DG
coalgebra and $\tau:\G\rarr A_1$ be an admissible twisting
cochain. Then there is a natural homotopy equivalence in
$\bar{A_1}^{op}_{(\G^*)^{op}}\text{-mod}_{\infty}$:
$$\G\otimes_{\tau} A_1\rarr f_*(\G\otimes_{f\cdot \tau} A_2).$$
\end{prop}
\begin{proof} We have a natural homotopy equivalence of DG
bicomodules over $B\bar{A_1}$:
$$B\bar{A_1}\otimes A_1[1]\otimes B\bar{A_1}\rarr B\bar{A_1}\square_{B\bar{A_2}} (B\bar{A_2}\otimes A_2[1]\otimes B\bar{A_2})\square_{B\bar{A_2}} B\bar{A_1}.$$
Co-tensoring it on the left by $\G$, we obtain the required
homotopy equivalence.
\end{proof}

Let $\mR$ be an artinian DG algebra, and $\tau:\mR^*\rarr A$ be an
admissible twisting cochain. Then by Proposition \ref{tw_coch} we
have a natural morphism of DG coalgebras $g_{\tau}:\mR^*\rarr
B\bar{A}$. Further, we have the dual morphism of DG algebras
$g_{\tau}^*:\hat{S}\rarr \mR$. In particular, $\mR$ becomes a DG
$\hat{S}^{op}$-module.

\begin{lemma}\label{remark} In the above notation  $A_{\infty}$ $\bar{A}^{op}_{\hat{S}^{op}}$-modules
$\Hom_{\hat{S}^{op}}(\mR,B\bar{A}\otimes_{\tau_A} A)$ and
$\mR^*\otimes_{\tau} A$ are isomorphic.
\end{lemma}
\begin{proof} Evident.\end{proof}

If $\tau:\mR^*\rarr A$ is an admissible twisting cochain and
$\alpha\in \cM\cC_{\mR}(A)$ is the corresponding object, then we
will write also $A\otimes_{\alpha} \mR^*$ instead of
$\mR^*\otimes_{\tau} A$.

Further, for $\alpha\in \cM\cC_{\mR}(A)$ corresponding to an
admissible twisting cochain we put
$$A\otimes_{\alpha} \mR:=\Hom_{\mR}(\mR^*,
A\otimes_{\alpha}\mR^*).$$ This is an object of
$\bar{A}^{op}_{\mR^{op}}\text{-mod}_{\infty}$. Its
$(\mR^{op})^{gr}$-module structure is obvious and
$A_{\infty}$-module structure can also be given by the explicit
formulas:

\begin{equation}\label{twisted_product_infty} m_n^{A\otimes_{\alpha} \mR}(m,a_1,\dots,a_{n-1})=
m_n^{0,\dots,0,\alpha}(m,a_1\otimes
\one_{\mR},\dots,a_{n-1}\otimes \one_{\mR}).\end{equation}

\begin{prop}\label{invar_h_def} Let $f:A_1\rarr A_2$ be an $A_{\infty}$-quasi-isomorphism of augmented
$A_{\infty}$-algebras, $\mR$ be an artinian DG algebra and let
$\alpha\in \cM\cC_{\mR}(A_1)$. Then there is a natural homotopy
equivalence in $\bar{A_1}_{\mR^{op}}^{op}\text{-mod}_{\infty}$:
$$A_1\otimes_{\alpha} \mR\rarr f_*(A_2\otimes_{f_{\mR}^*(\alpha)} \mR).$$
\end{prop}
\begin{proof} The required homotopy
equivalence is obtained by applying the functor
$\Hom_{\mR}(\mR^*,-)$ to the homotopy equivalence
$$A_1\otimes_{\alpha} \mR^*\rarr
f_*(A_2\otimes_{f_{\mR}^*(\alpha)}\mR^*)$$ from Proposition
\ref{invar_h}.
\end{proof}

Note that if $A$ is a DG algebra, then  $A\otimes_{\alpha} \mR^*$
and $A\otimes_{\alpha} \mR$ are the DG modules from
$\coDef^h_{\mR}(A)$ and $\Def^h_{\mR}(A)$ respectively, which
correspond to $\alpha$.

Finally, if $A$ is a strictly unital but not necessarily augmented
$A_{\infty}$-algebra, $\mR$ is an artinian DG algebra, $\alpha$ is
an object of $\cM\cC_{\mR}(A)$ and $\tau:\mR^*\to A$ is the
corresponding twisting cochain then we also write
$A\otimes_{\alpha}\mR^*$ instead of $\mR^*\otimes_{\tau} A$.
Further, we put
$$A\otimes_{\alpha}\mR=\Hom_{\mR}(\mR^*,A\otimes_{\alpha}\mR^*).$$
This is the object of $A^{op}_{\mR^{op}}\text{-mod}_{\infty}$.
Again, its $(\mR^{op})^{gr}$-module structure is obvious and the
$A_{\infty}$-module structure is given by the formulas
\ref{twisted_product_infty}. The following Proposition is
absolutely analogous to the previous one and we omit the proof.

\begin{prop}\label{invar_h_def1} Let $f:A_1\rarr A_2$ be a strictly unital $A_{\infty}$-morphism of
strictly unital $A_{\infty}$-algebras, $\mR$ be an artinian DG
algebra and let $\alpha\in \cM\cC_{\mR}(A_1)$. Then there is a
natural homotopy equivalence in
$A^{op}_{1\mR^{op}}\text{-mod}_{\infty}$:
$$A_1\otimes_{\alpha}\mR\rarr f_*(A_2\otimes_{f_{\mR}^*(\alpha)}\mR).$$
\end{prop}

\part{The pseudo-functors $\DEF $ and $\coDEF$}

\section{The bicategory $2\text{-}\adgalg$ and deformation
pseudo-functor $\coDEF$} \label{2-adgalg_codef}

Let $\cE$ be a bicategory and $F,G:\cE \to {\bf Gpd}$ two
pseudo-functors. A morphism $\epsilon :F\to G$ is called an
equivalence if for each $X\in Ob\cE$ the functor $\epsilon _X:
F(X)\to G(X)$ is an equivalence of categories.

\begin{defi}\label{2-adgalg} We define the bicategory $2\text{-}\adgalg$ of augmented DG algebras as follows.
The objects are augmented DG algebras. For DG algebras $\cB, \cC$
the collection of 1-morphisms $1\text{-}\Hom(\cB,\cC)$ consists of
pairs $(M,\theta)$, where
\begin{itemize}
\item $M\in D(\cB ^{op}\otimes \cC)$ is such that there exists an
isomorphism (in $D(\cC)$)  $\cC\to \nu _*M$ (where $\nu _*:D(\cB
^{op} \otimes \cC)\to D(\cC)$ is the functor of restriction of
scalars corresponding to the natural homomorphism $\nu :\cC  \to
\cB ^{op} \otimes \cC$); \item and $\theta
:k{\stackrel{\bL}{\otimes }}_{\cC}M\to k$ is an isomorphism in
$D(\cB ^{op})$.
\end{itemize}
The  composition of 1-morphisms
$$1\text{-}\Hom(\cB,\cC)\times 1\text{-}\Hom(\cC,\cD)\to
1\text{-}\Hom(\cB,\cD)$$
 is defined by the tensor product $\cdot
 \stackrel{\bL}{\otimes}_{\cC}\cdot $.
 Given
1-morphisms $(M_1,\theta _1), (M_2,\theta _2)\in
1\text{-}\Hom(\cB,\cC)$ a 2-morphism $f: (M_1,\theta _1)\to
(M_2,\theta _2)$ is an isomorphism (in $D(\cB ^{op} \otimes \cC)$)
$f:M_2\to M_1$ (not from $M_1$ to $M_2$!) such that $\theta
_1\cdot k{\stackrel{\bL}{\otimes }}_{\cC}(f)=\theta _2$. So in
particular the category $1\text{-}\Hom(\cB ,\cC)$ is a groupoid.
Denote by $2\text{-}\dgart$ the full subbicategory of
$2\text{-}\adgalg$ consisting of artinian DG algebras. Similarly
we define the full subbicategories $2\text{-}\dgart_+$,
$2\text{-}\dgart_-$, $2\text{-}\art$, $2\text{-}\cart$ (I,
Definition 2.3).
\end{defi}

\begin{remark}\label{remark_on_adgalg} Assume that augmented DG algebras $\cB$ and $\cC$ are
such that $\cB ^i=\cC ^i=0$ for $i>0$, $\dim \cB ^i,\dim \cC^i
<\infty$ for all $i$ and $\dim H(\cC)<\infty$. Denote by $\langle
k\rangle \subset D(\cB ^{op} \otimes \cC )$ the triangulated
envelope of the DG $\cB ^{op}\otimes \cC$-module $k$. Let
$(M,\theta )\in 1\text{-}\Hom (\cB ,\cC)$. Then by I, Corollary
3.22 $M\in \langle k\rangle$.
\end{remark}

For any augmented DG algebra $\cB$ we obtain a pseudo-functor
$h_{\cB}$ between the bicategories $2\text{-}\adgalg$ and {\bf Gpd}
defined by $h_{\cB}(\cC)=1\text{-}\Hom (\cB ,\cC)$.

Note that a usual homomorphism of augmented
 DG algebras $\gamma :\cB
\to \cC$ defines the structure of a DG $\cB^{op}$-module on $\cC$
with the canonical isomorphism of DG $\cB^{op}$-modules $\id :
k{\stackrel{\bL}{\otimes }}_{\cC}\cC \to k$. Thus it defines a
1-morphism $(\cC ,\id)\in 1\text{-}\Hom (\cB ,\cC)$. This way we
get a pseudo-functor $\cF:\adgalg \to 2\text{-}\adgalg$, which is
the identity on objects.

\begin{lemma}\label{2-adgalg_vs_alg} Assume that augmented DG algebras $\cB$ and $\cC$ are concentrated in degree zero
(hence have zero differential). Also assume that these algebras
are local (with maximal ideals being the augmentation ideals).
Then

a) the map $\cF:\Hom (\cB ,\cC)\to \pi _0(1\text{-}\Hom (\cB,\cC))$
is surjective, i.e. every 1-morphism from $\cB$ to $\cC$ is
isomorphic to $\cF(\gamma)$ for a homomorphism of algebras $\gamma$;

b) the 1-morphisms $\cF(\gamma _1)$ and $\cF(\gamma _2)$ are
isomorphic if and only if $\gamma _2$ is the composition of $\gamma
_1$ with the conjugation by an invertible element in $\cC$;

c) in particular, if $\cC$ is commutative then the map of sets
$\cF:\Hom (\cB ,\cC)\to \pi _0(1\text{-}\Hom (\cB,\cC))$ is a
bijection.
\end{lemma}

\begin{proof} a) For any  $(M, \theta)\in 1\text{-}\Hom (\cB,\cC)$ the DG
$\cB ^{op} \otimes \cC$-module $M$ is isomorphic (in $D(\cB ^{op}
\otimes \cC)$ to  $H^0(M)$. Thus we may assume that $M$ is
concentrated in degree 0. By assumption there exists an
isomorphism of $\cC$-modules
 $\cC \to M$. Multiplying this isomorphism by a scalar we may assume that
 it is compatible with the isomorphisms $\id :k{\stackrel{\bL}{\otimes }}_{\cC}\cC\to k$
 and $\theta :k{\stackrel{\bL}{\otimes }}_{\cC}M\to k$.
  A choice of such an isomorphism defines a homomorphism of algebras
 $\cB ^{op}\to \End _{\cC}(\cC)=\cC ^{op}$. Since $\cB$ and $\cC$ are local
 this is a homomorphism of augmented algebras.  Thus
$(M,\theta)$ is isomorphic to $\cF (\gamma)$.

b) Let $\gamma _1,\gamma _2:\cB \to \cC$ be homomorphisms of
algebras. A 2-morphism $f:\cF (\gamma _1)\to \cF (\gamma _2)$ is
simply an isomorphism of the corresponding $\cB ^{op}\otimes \cC
$-modules $f:\cC \to \cC$, which commutes with the augmentation.
Being an isomorphism of $\cC$-modules it is the right
multiplication by an invertible element $c\in \cC$. Hence for
every $b \in \cB$ we have $c^{-1}\gamma _1(b)c=\gamma _2(b)$.

c) This follows from a) and b).
\end{proof}

\begin{remark}\label{quasi-functor} If in the definition of 1-morphisms $1\text{-}\Hom(\cB,\cC)$ we do not fix an
isomorphism $\theta$, then we obtain a special case of a
"quasi-functor" between the DG categories $\cB \text{-mod}$
 and $\cC \text{-mod}$. This notion was first introduced by Keller
 in [Ke] for DG modules over general DG categories.
 \end{remark}

The next proposition asserts that the deformation functor $\coDef$
has a natural "lift" to the bicategory $2\text{-}\dgart$.

\begin{prop}\label{coDEF} There exist a pseudo-functor  $\coDEF(E)$ from
$2\text{-}\dgart$ to {\bf Gpd} which is an extension to
$2\text{-}\dgart$ of the pseudo-functor $\coDef$ i.e. there is  an
equivalence of pseudo-functors $\coDef(E)\simeq \coDEF(E)\cdot
\cF$.
\end{prop}

\begin{proof} Given artinian DG
algebras $\cR, \cQ$  and $M=(M,\theta )\in 1\text{-}\Hom (\cR ,\cQ)$
we need to define the corresponding functor
$$M^!:\coDef _{\cR}(E)\to \coDef _{\cQ}(E).$$
Let $S=(S,\sigma)\in \coDef _{\cR}(E)$. Put
$$M^!(S):=\bR \Hom _{\cR^{op}}(M,S)\in D(\cA _{\cQ}^{op}).$$
We claim that $M^!(S)$ defines an object in $\coDef _{\cQ}(E)$,
i.e. $\cR \Hom _{\cQ^{op}}(k,M^!(S))$ is naturally isomorphic to
$E$ (by the isomorphisms $\theta $ and $\sigma$).

 Indeed, choose quasi-isomorphisms $P\to k$ and $S\to I$
for $P\in \cP(\cA _{\cQ}^{op})$ and $I\in \cI(\cA _{\cR}^{op})$.
Then
$$\bR\Hom _{\cQ^{op}}(k,M^!(S))=\Hom _{\cQ^{op}}(P,\Hom _{\cR^{op}}(M,I)).$$
By I, Lemma 3.17 the last term is equal to $\Hom _{\cR
^{op}}(P\otimes _{\cQ}M,I)$. Now the isomorphism $\theta$ defines
an isomorphism between $P\otimes _{\cQ}M=k\stackrel{\bL}{\otimes}
_{\cQ}M$ and $k$, and we compose it with the isomorphism $\sigma
:E \to \bR\Hom _{\cR ^{op}}(k,I)=i^!S$.

So $M^!$ is a functor from $\coDef_{\cR}(E)$ to $\coDef _{\cQ}(E)$.

Given another artinian DG algebra $\cQ ^\prime$ and $M^\prime \in
1\text{-}\Hom (\cQ ,\cQ ^\prime)$ there is a natural isomorphism of
functors
$$(M^\prime
\stackrel{\bL}{\otimes}_{\cQ}M)^!(-) \simeq M^{\prime !}\cdot M
^!(-).$$ (This follows again from I, Lemma 3.17).

Also a 2-morphism $f\in 2\text{-}\Hom(M,M_1)$ between objects
$M,M_1\in 1\text{-}\Hom (\cR ,\cQ)$ induces an isomorphism of the
corresponding functors $M^!\stackrel{\sim}{\to}M_1^!$.

Thus we obtain a pseudo-functor $\coDEF (E):2\text{-}\dgart \to {\bf
Gpd}$, such that $\coDEF (E)\cdot \cF=\coDef(E)$.
\end{proof}

We denote by $\coDEF_+(E)$, $\coDEF_-(E)$, $\coDEF_0(E)$,
$\coDEF_{\cl}(E)$ the restriction of the pseudo-functor $\coDEF (E)$
to subbicategories $2\text{-}\dgart _+$, $2\text{-}\dgart _-$,
$2\text{-}\art$ and $2\text{-}\cart$ respectively.

\begin{prop}\label{coDEF_invar} A quasi-isomorphism $\delta :E_1\to E_2$ of DG $\cA
^{op}$-modules induces an equivalence of pseudo-functors
$$\delta ^*:\coDEF (E_2)\to \coDEF (E_1)$$
defined by $\delta ^*(S,\sigma )=(S,\sigma \cdot \delta)$.
\end{prop}

\begin{proof} This is clear.
\end{proof}

\begin{prop}\label{coDEF_invar1} Let $F:\cA \to \cA ^\prime$ be a DG functor which
induces a quasi-equivalence $F^{\text{pre-tr}}:\cA
^{\text{pre-tr}}\to \cA ^{\prime \text{pre-tr}}$ (this happens for
example if $F$ is a quasi-equivalence). Then for any $E\in D(\cA
^{op})$ the pseudo-functors $\coDEF _-(E)$ and $\coDEF _-(\bR
F^!(E))$ are equivalent (hence also $\coDEF (F_*(E^\prime))$ and
$\coDEF _-(E^\prime)$ are equivalent for any $E^\prime \in D(\cA
^{\prime 0})$).
\end{prop}

\begin{proof} The proof is similar to the proof of I, Proposition
10.11. Namely let $R,Q \in 2\text{-}\dgart _-$ and $M\in
1\text{-}\Hom (R,Q)$. The DG functor $F^!$ induces a commutative
functorial diagram
$$\begin{array}{ccc}
D(\cA ^{op}_{\cR}) & \stackrel{\bR (F\otimes
\id)^!}{\longrightarrow} &
D(\cA ^{\prime 0}_{\cR})\\
\bR i^! \downarrow & & \downarrow \bR i^! \\
D(\cA ^{op}) & \stackrel{\bR F^!}{\longrightarrow} & D(\cA
^{\prime 0})
\end{array}$$
(and a similar diagram for $\cQ$ instead of $\cR$) which is
compatible with the functors
$$M^!:D(\cA _{\cR}^{op}) \to D(\cA ^{op}_{\cQ}), \quad \text{and}\quad
M^!:D(\cA _{\cR}^{\prime 0}) \to D(\cA ^{\prime 0}_{\cQ}).$$ Thus we
obtain a morphism of pseudo-functors
$$F^!:\coDEF _-(E)\to \coDEF _-(\bR F^!(E)).$$

By I, Corollary 3.15 the functors $\bR F^!$ and $\bR (F\otimes \id
)^!$ are equivalences.
\end{proof}

\begin{cor}\label{codef(B)=codef(C)} Assume that DG algebras $\cB$ and $\cC$ are
quasi-isomorphic. Then the pseudo-functors $\coDEF _-(\cB)$ and
$\coDEF _-(\cC)$ are equivalent.
\end{cor}

\begin{proof} We may assume that there exists a homomorphism of DG algebras $\cB \to \cC$ which is a quasi-isomorphism.
Then put $\cA=\cB$ and $\cA ^\prime =\cC$ in the last proposition.
\end{proof}

The following Lemma is stronger then I, Corollary 11.15 for the
pseudo-functors $\coDef_-$ and $\coDef_-^h$.

\begin{lemma}\label{codef=codef^h} Let $\B$ be a DG algebra. Suppose that the
following conditions hold:

a) $H^{-1}(\B)=0$;

b) the graded algebra $H(\B)$ is bounded below.

Then the pseudo-functors $\coDef_-(\B)$ and $\coDef_-^h(\B)$ are
equivalent.
\end{lemma}
\begin{proof} Fix some negative artinian DG algebra $\mR\in \dgart_-$.
Take some $(T,id)\in \coDef_{\mR}^h(\B)$. Due to I, Corollary 11.4
b) it suffices to prove that $i^!T=\bfR i^!T$ Let $A$ be a
strictly unital minimal model of $\B$, and let $f:A\rarr \B$ be a
strictly unital $A_{\infty}$ quasi-isomorphism. By our assumption
on $H(\B)$, $A$ is bounded below.

By Theorem \ref{invar_MC} there exists an object $\alpha\in
\cM\cC_{\mR}(A)$ such that $S\cong
\B\otimes_{f_{\mR}^*(\alpha)}\mR^*$. The DG $\mR^{op}$-modules
$\B\otimes_{f_{\mR}^*(\alpha)}\mR^*$ and
$f_*(\B\otimes_{f_{\mR}^*(\alpha)}\mR^*)$ are naturally
identified. Further,  by Proposition \ref{invar_h_def1} we have
natural homotopy equivalence (in
$A^{op}_{\mR^{op}}\text{-mod}_{\infty}$)
$$\gamma:A\otimes_{\alpha} \mR^*\rarr f_*(\B\otimes_{f_{\mR}^*(\alpha)}\mR^*)$$
Thus, it remains to prove that
$$i^!(A\otimes_{\alpha} \mR^*)=\bfR i^!(A\otimes_{\alpha}
\mR^*).$$ We claim that $A\otimes_{\alpha} \mR^*$ is h-injective.
Indeed, since $A$ is bounded below and $\mR\in \dgart_-$, this DG
$\mR^{op}$-module has a decreasing filtration by DG
$\mR^{op}$-submodules $A^{\geq i}\otimes \mR^*$ with subquotients
being cofree DG $\mR^{op}$-modules $A^i\otimes \mR^*$. Thus
$A\otimes_{\alpha} \mR^*$ satisfies property (I) as DG
$\mR^{op}$-module and hence is h-injective. Lemma is proved.
\end{proof}

The next result implies stronger statement for pseudo-functor
$\coDef_-$ then I, Proposition 11.16.

\begin{prop}\label{codef(B)=codef(E)} Let $E\in \cA ^{op}\text{-mod}$. Assume that

a) $\Ext ^{-1}(E,E)=0$;

b) the graded algebra $\Ext(E,E)$ is bounded below;

b) there exists a bounded below h-projective or h-injective DG
$\cA ^{op}$-module $F$ which is quasi-isomorphic to $E$.

Put $\B=\End(F)$. Then the pseudo-functors $\coDEF _-(\B)$ and
$\coDEF _-(E)$ are equivalent.
\end{prop}

\begin{proof} Consider the DG functor
$$\cL :=\Sigma ^F\cdot \psi ^*:\B ^{op}\text{-mod}\to \cA
^{op}\text{-mod}, \quad \cL(N)=N\otimes _{\cC}F$$ as in I, Remark
11.17.
 It induces the equivalence of pseudo-functors
$$\coDef ^{\h}(\cL):\coDef ^{\h}_-(\B)\stackrel{\sim}{\to}\coDef
^{\h}_-(F),$$ i.e. for every artinian DG algebra $\cR \in \dgart _-$
the corresponding DG functor
$$\cL _{\cR}:(\B\otimes \cR )^{op}\text{-mod}\to \cA
_{\cR}^{op}\text{-mod}$$ induces the equivalence of groupoids
$\coDef^{\h} _{\cR}(\B)\stackrel{\sim}{\to}\coDef^{\h}_{\cR}(F)$
(I, Propositions 9.2, 9.4). By I, Theorem 11.6 b) there is a
natural equivalences of pseudo-functors $$\coDef ^{\h}_-(F)\simeq
\coDef _-(E).$$ By Lemma \ref{codef=codef^h} there is an
equivalence of pseudo-functors
$$\coDef ^{\h}_-(\B)\simeq \coDef _-(\B).$$ Hence the functor
$\bL\cL$ induces the equivalence
$$\bL\cL:\coDef _-(\B)\stackrel{\sim}{\to}\coDef _-(E).$$

Fix $\cR,\cQ \in 2\text{-}\dgart _-$ and $M\in
1\text{-}\Hom(\cR,\cQ)$. We need to show that there exists a natural
isomorphism
$$\bL\cL _{\cQ}\cdot M^!\simeq M^!\cdot \bL\cL _{\cR}$$
between functors from $\coDef _{\cR}(\B)$ to $\coDef _{\cQ}(E)$.

Since the cohomology of $M$ is finite dimensional, and the DG
algebra $\cR \otimes \cQ$ has no components in positive degrees, by
I, Corollary 3.21 we may assume that $M$ is finite dimensional.

\begin{lemma}\label{M^!=hom(M,-)} Let $(S, \id)$ be an object in $\coDef
^{\h}_{\cR}(\B)$ or in $\coDef ^{\h}_{\cR}(F)$. Then $S$ is
acyclic for the functor $\Hom _{\cR ^{op}}(M,S)$, i.e.
$M^!(S)=\Hom _{\cR ^{op}}(M,S)$.
\end{lemma}

\begin{proof} In the proof of Lemma \ref{codef=codef^h} (resp. in I, Lemma 11.8) we showed that $S$ is
h-injective when considered as a DG $\cR^{op}$-module.
\end{proof}

Choose $(S,\id)\in \coDef ^{\h}(\B)$. By the above lemma
$M^!(S)=\Hom _{\cR^{op}}(M,S)$.

We claim that the DG $\B^{op}$-module $\Hom _{\cR^{op}}(M,S)$ is
h-projective. Indeed, first notice that the graded $\cR
^{op}$-module $S$ is injective being isomorphic to a direct sum of
copies of shifted graded $\cR^{op}$-module $\cR ^*$ (the abelian
category of graded $\cR ^{op}$-modules is locally notherian, hence
a direct sum of injectives is injective). Second, the DG $\cR
^{op}$-module $M$ has a (finite) filtration with subquotients
isomorphic to $k$. Thus the DG $\B ^{op}$-module $\Hom _{\cR
^{op}}(M,S)$ has a filtration with subquotients isomorphic to
$\Hom _{\cR ^{op}}(k,S)=i^!S\simeq \B$. So it has property (P).

Hence $\bL\cL\cdot M^!(S)=\Hom _{\cR^{op}}(M,S)\otimes _{\B}F$.
For the same reasons $M^!\cdot \bL\cL _{\cR}(S)=\Hom
_{\cR^{op}}(M, S\otimes _{\B}F)$. The isomorphism
$$\Hom_{\cR^{op}}(M,S)\otimes _{\B}F=\Hom _{\cR^{op}}(M, S\otimes _{\B}F)$$
follows from the fact that $S$ as a graded module is a tensor
product of graded $\cC^{op}$ and $\cR^{op}$ modules and also
because $\dim _kM<\infty$.
\end{proof}

\section{Deformation pseudo-functor $\coDEF$ for an augmented
$A_{\infty}$-algebra}

Let $A$ be an augmented $A_{\infty}$-algebra. We are going to
define the pseudo-functor $\coDEF(A):2\text{-}\dgart\rarr {\bf
Gpd}$.

Let $\mR$ be an artinian DG algebra. An object of the groupoid
$\coDEF_{\mR}(A)$ is a pair $(S,\sigma)$, where $S\in
D_{\infty}(\bar{A}^{op}_{\mR^{op}})$, and $\sigma$ is an
isomorphism (in $D_{\infty}(\bar{A}^{op})$)
$$\sigma:A\rarr \bfR i^!(S).$$

A morphism $f:(S,\sigma)\rarr (T,\tau)$ in $\coDEF_{\mR}(A)$ is an
isomorphism (in $D(\bar{A}^{op}_{\cR^{op}})$) $f:S\rarr T$ such
that
$$\bfR i^!(f)\circ \sigma=\tau.$$

This defines the pseudo-functor $\coDEF(A)$ on objects. Further,
let $(M,\theta)\in 1\text{-}\Hom(\mR,\mQ)$. Define the
corresponding functor
$$M^!:\coDEF_{\mR}(A)\rarr \coDEF_{\mQ}(A)$$ as follows. For an
object $(S,\sigma)\in \coDEF_{\mR}(A)$ put
$$M^!(S)=\bfR \Hom_{\mR^{op}}(M,S)\in D_{\infty}(\bar{A}_{\mQ^{op}}^{op}).$$

Then we have natural isomorphisms in $D_{\infty}(\bar{A})$:
$$\bfR \Hom_{\mQ^{op}}(k,M^!(S))\cong \bfR \Hom_{\mR^{op}}(k\otimes_{\mR^{op}}^{\bfL}M,S)\cong \bfR \Hom_{\mR^{op}}(k,S)=\bfR i^!(S)$$
(the second isomorphism is induced by $\theta$). Thus, $M^!$ is a
functor form $\coDEF_{\mR}(A)$ to $\coDEF_{\mQ}(A)$.

If $\mQ'$ is another artinian DG algebra and $(M',\theta')\in
1\text{-}\Hom(\mQ,\mQ')$ then there is a natural isomorphism of
functors
$$(M'\otimes_{\mQ}^{\bfL} M)^!\cong M'^!\cdot M^!.$$

Further, if $f\in 2\text{-}\Hom((M,\theta),(M,\theta_1))$ is a
$2$-morphism between objects $(M,\theta),(M,\theta_1)\in
1\text{-}\Hom(\mR,\mQ)$ then it induces an isomorphism between the
corresponding functors $M^!\simeq M_1^!$.

Thus we obtain a pseudo-functor $\coDEF(A):2\text{-}\dgart\rarr
{\bf Gpd}$. We denote by $\coDEF_-(A)$ its restriction to the
sub-$2$-category $2\text{-}\dgart_-$.

\begin{prop}\label{codef-invar} Let $A$ be an augmented $A_{\infty}$-algebra and $U(A)$ its bar-cobar construction. Then there is a natural
equivalence of pseudo-functors $\coDEF(U(A))\cong  \coDEF(A)$.
\end{prop}
\begin{proof} Let $f_A:A\rarr U(A)$ be the universal strictly
unital $A_{\infty}$-morphism. Let $\mR$ be an artinian DG algebra.
Recall that by Proposition \ref{AequivU(A)} we have an equivalence
$$f_{A*}:D((U(A)\otimes
\mR)^{op})\rarr D_{\infty}(\bar{A}^{op}_{\mR^{op}}).$$ Moreover,
the following diagram of functors commutes up to an isomorphism:

$$
\begin{CD}
D((U(A)\otimes \mR)^{op}) @>f_{A*}>>  D_{\infty}(\bar{A}^{op}_{\mR^{op}})\\
@VV\bfR i^! V                       @VV\bfR i^! V\\
D(U(A)^{op}) @>f_{A*}>> D_{\infty}(\bar{A}^{op}).
\end{CD}
$$

Hence, the functor $f_{A*}$ induces an equivalence of groupoids
$\coDEF_{\mR}(U(A))\rarr \coDEF_{\mR}(A)$ and we obtain the
required equivalence of pseudo-functors.\end{proof}

\begin{cor} Let $A$ be an augmented $A_{\infty}$-algebra and let
$\cB$ be a DG algebra quasi-isomorphic to $A$. Then the
pseudo-functor $\coDEF(A)$ and $\coDEF(\cB)$ are equivalent.
\end{cor}
\begin{proof} Indeed, by Proposition \ref{codef-invar} the
pseudo-functors $\coDEF(A)$ and $\DEF(U(A))$ are equivalent, and
by Corollary \ref{codef(B)=codef(C)} the pseudo-functors
$\coDEF(U(A))$ and $\coDEF(\cB)$ are equivalent.
\end{proof}

\begin{cor}\label{codef^h} Let $A$ be an admissible $A_{\infty}$-algebra, and $\mR$ be an artinian negative DG algebra. Then for any
$(S,\sigma)\in \coDEF_{\mR}(A)$ there exists a morphism of DG
algebras $\hat{S}\rarr \mR$ such that the pair $(T,id)$, where
$T=\Hom_{\hat{S}^{op}}(\mR,B\bar{A}\otimes_{\tau_A} A)$, defines
an object of $\coDEF_{\mR}(A)$ which is isomorphic to
$(S,\sigma)$.\end{cor}
\begin{proof} This follows easily from Proposition \ref{codef-invar}, the proof of Lemma
\ref{codef=codef^h} in the case $\B=U(A)$, and Lemma
\ref{remark}.\end{proof}

\section{The bicategory $2^\prime\text{-}\adgalg$ and deformation
pseudo-functor $\DEF$}

It turns out that the deformation pseudo-functor $\Def$ lifts
naturally to a different version of a bicategory of augmented DG
algebras. We denote this bicategory $2^\prime\text{-}\adgalg$. It
differs from $2\text{-}\adgalg$ in two respects: the 1-morphisms
are objects in $D(\cB \otimes \cC ^{op})$ (instead of $D(\cB ^{op}
\otimes \cC)$) and 2-morphisms go in the opposite direction. We
will relate the bicategories $2\text{-}\adgalg$ and $2^\prime
\text{-}\adgalg$ (and the pseudo-functors $\coDEF$ and $\DEF$) in
section \ref{comparison} below.

\begin{defi}\label{2'-adgalg} We define the bicategory $2^\prime \text{-}\adgalg$ of augmented  DG algebras as follows.
The objects are augmented DG algebras. For DG algebras $\cB, \cC$
the collection of 1-morphisms $1\text{-}\Hom(\cB,\cC)$ consists of
pairs $(M,\theta)$, where
\begin{itemize}
\item $M\in D(\cB  \otimes \cC ^{op})$ and there exists an
isomorphism (in $D(\cC ^{op})$)  $\cC\to \nu _*M$ (where $\nu
_*:D(\cB \otimes \cC ^{op})\to D(\cC ^{op})$ is the functor of
restriction of scalars corresponding to the natural homomorphism
$\nu :\cC ^{op}  \to \cB \otimes \cC ^{op}$); \item and $\theta
:M{\stackrel{\bL}{\otimes }}_{\cC}k\to k$ is an isomorphism in
$D(\cB)$.
\end{itemize}
The composition of 1-morphisms
$$1\text{-}\Hom(\cB,\cC)\times 1\text{-}\Hom(\cC,\cD)\to
1\text{-}\Hom(\cB,\cD)$$
 is defined by the tensor product $\cdot
 \stackrel{\bL}{\otimes}_{\cC}\cdot $.
 Given
1-morphisms $(M_1,\theta _1), (M_2,\theta _2)\in
1\text{-}\Hom(\cB,\cC)$ a 2-morphism $f: (M_1,\theta _1)\to
(M_2,\theta _2)$ is an isomorphism (in $D(\cB  \otimes \cC
^{op})$) $f:M_1\to M_2$ such that $\theta _1=\theta _2\cdot
((f){\stackrel{\bL}{\otimes}}_{\cC}k)$. So in particular the
category $1\text{-}\Hom(\cB ,\cC)$ is a groupoid. Denote by
$2^\prime\text{-}\dgart$ the full subbicategory of
$2^\prime\text{-}\adgalg$ consisting of artinian DG algebras.
Similarly we define the full subbicategories
$2^\prime\text{-}\dgart_+$, $2^\prime\text{-}\dgart_-$,
$2^\prime\text{-}\art$, $2^\prime\text{-}\cart$ (I, Definition
2.3).
\end{defi}

\begin{remark}\label{remark_on_'adgalg} The exact analogue of Remark \ref{remark_on_adgalg} holds for the
bicategory $2^\prime\text{-}\adgalg$.
\end{remark}

For any augmented DG algebra $\cB$ we obtain a pseudo-functor
$h^\prime_{\cB}$ between the bicategories $2^\prime\text{-}\adgalg$
and {\bf Gpd} defined by $h^\prime_{\cB}(\cC)=1\text{-}\Hom (\cB
,\cC)$.

Note that a usual homomorphism of DG algebras $\gamma :\cB \to \cC$
defines the structure of a $\cB$-module on $\cC$ with the canonical
isomorphism of DG $\cB$-modules
$\cC{\stackrel{\bL}{\otimes}}_{\cC}k$. Thus it defines a 1-morphism
$(\cC ,\id)\in 1\text{-}\Hom (\cB ,\cC)$. This way we get a
pseudo-functor $\cF^\prime:\adgalg \to 2^\prime\text{-}\adgalg$,
which is the identity on objects.

\begin{remark}\label{2'-adgalg_vs_alg}
The precise analogue of Lemma \ref{2-adgalg_vs_alg} holds for the
bicategory $2^\prime \text{-}\adgalg$ and the pseudo-functor $\cF
^\prime$.
\end{remark}

\begin{prop}\label{DEF} There exist a pseudo-functor  $\DEF(E)$ from
$2^\prime\text{-}\dgart$ to {\bf Gpd} and which is an extension to
$2^\prime\text{-}\dgart$ of the pseudo-functor $\Def (E)$, i.e.
there is an equivalence of pseudo-functors $\Def(E)\simeq
\DEF(E)\cdot \cF ^\prime$.
\end{prop}

\begin{proof} Let $\cR$,
$\cQ$ be artinian DG algebras. Given $(M,\theta )\in 1\text{-}\Hom
(\cR ,\cQ)$ we define the corresponding functor
$$M^*: \Def _{\cR}(E)\to \Def _{\cQ}(E)$$ as follows
$$M^*(S):=S\stackrel{\bL}{\otimes}_{\cR}M$$
for $(S,\sigma)\in \Def _{\cR}(E)$. Then we have the canonical
isomorphism
$$M^*(S)\stackrel{\bL}{\otimes }_{\cQ}k=
S\stackrel{\bL}{\otimes }_{\cR}(M\stackrel{\bL}{\otimes }_{\cQ}k)
\stackrel{\theta}{\to} S\stackrel{\bL}{\otimes
}_{\cR}k\stackrel{\sigma}{\to}E.$$
 So that $M^*(S)\in \Def
_{\cQ}(E)$ indeed.

Given another artinian DG algebra $\cQ ^\prime$ and $M^\prime \in
1\text{-}\Hom (Q,Q^\prime)$ there is a natural isomorphism of
functors
$$M^{\prime *}\cdot M^*=(M\stackrel{\bL}{\otimes
}_{\cQ}M^\prime)^*.$$ Also a 2-morphism $f\in 2\text{-}\Hom (M,M_1)$
between $M,M_1\in 1\text{-}\Hom (\cR ,\cQ)$ induces an isomorphism
of corresponding functors $M^*\stackrel{\sim}{\to}M^*_1$.

Thus we obtain a pseudo-functor $\DEF(E):2^\prime \text{-}\dgart \to
{\bf Gpd}$, such that $\DEF(E)\cdot \cF ^\prime =\Def(E)$.
\end{proof}

We denote by $\DEF_+(E)$, $\DEF_-(E)$, $\DEF_0(E)$, $\DEF_{\cl}(E)$
the restriction of the pseudo-functor $\DEF (E)$ to subbicategories
$2^\prime\text{-}\dgart _+$, $2^\prime\text{-}\dgart _-$,
$2^\prime\text{-}\art$ and $2^\prime\text{-}\cart$ respectively.

\begin{prop}\label{DEF_invar} A quasi-isomorphism $\delta :E_1\to E_2$ of DG $\cA
^{op}$-modules induces an equivalence of pseudo-functors
$$\delta _*:\DEF (E_1)\to \DEF (E_2)$$
defined by $\delta _*(S,\sigma)=(S,\delta \cdot \sigma)$.
\end{prop}

\begin{proof} This is clear.
\end{proof}

\begin{prop}\label{DEF_invar1} Let $F:\cA \to \cA ^\prime$ be a DG functor which
induces a quasi-equivalence $F^{\text{pre-tr}}:\cA
^{\text{pre-tr}}\to \cA ^{\prime \text{pre-tr}}$ (this happens for
example if $F$ is a quasi-equivalence). Then for any $E\in D(\cA
^{op})$ the pseudo-functors $\DEF _-(E)$ and $\DEF _-(\bL F^*
(E))$ are equivalent (hence also $\DEF _-(F_*(E^\prime))$ and
$\DEF _-(E^\prime )$ are equivalent for any $E^\prime \in D(\cA
^{\prime 0})$).
\end{prop}

\begin{proof}
The proof is similar to the proof of I, Proposition 10.4. Let $\cR
,\cQ \in \dgart _-$ and $M\in 1\text{-}\Hom (\cR ,\cQ)$. The DG
functor $F$ induces a commutative functorial diagram
$$\begin{array}{ccc} D(\cA _{\cR} ^{op}) & \stackrel{\bL (F\otimes
\id)^*}{\longrightarrow} & D(\cA _{\cR}^{\prime 0})\\
\bL i^* \downarrow & & \downarrow \bL i^*\\
D(\cA ^{op}) & \stackrel {\bL F^*}{\longrightarrow} & D(\cA
^{\prime 0})
\end{array}$$
(and a similar diagram for $\cQ$ instead of $\cR$) which is
compatible with the functors
$$M^* : D(\cA _{\cR}^{op})\to D(\cA ^{op}_{\cQ})\quad \text{and}\quad M^*
:D(\cA _{\cR}^{\prime 0}\to D(\cA _{\cQ}^{\prime 0}).$$ Thus we
obtain a morphism of pseudo-functors
$$F^*: \DEF _-(E)\to \DEF _-(\bL F^*(E)).$$
By I, Corollary 3.15 the functors $\bL F^*$ and $\bL (F\otimes
\id)^*$ are equivalences, hence this morphism $F^*$ is an
equivalence.
\end{proof}

\begin{cor}\label{def(B)=def(C)} Assume that DG algebras $\cB$ and $\cC$ are
quasi-isomorphic. Then the pseudo-functors $\DEF _-(\cB)$ and
$\DEF _-(\cC)$ are equivalent.
\end{cor}

\begin{proof} We may assume that there exists a morphism of DG algebras $\cB \to \cC$ which
is a quasi-isomorphism. Then put $\cA =\cB$ and $\cA ^\prime =\cC$
in the last proposition.
\end{proof}

The following Theorem is stronger then I, Corollary 11.15 for the
pseudo-functors $\Def_-$ and $\Def_-^h$.

\begin{theo}\label{def=def^h} Let $E\in \A^{op}\text{-mod}$ be a
DG module. Suppose that the following conditions hold:

a) $\Ext^{-1}(E,E)=0$;

b) the graded algebra $\Ext(E,E)$ is bounded above.

Let $F\rarr E$ be a quasi-isomorphism with h-projective $F$. Then
the pseudo-functors $\Def_-(E)$ and $\Def_-^h(F)$ are equivalent.
\end{theo}
\begin{proof} Replace the
pseudo-functor $\Def(E)$ by the equivalent pseudo-functor
$\Def(F)$. Fix some negative artinian DG algebra $\mR\in
\dgart_-$.

Due to I, Corollary 11.4 a) it suffices to prove that for each
$(S,id)\in \Def_{\mR}^h(F)$ one has $i^*(S)=\bfL i^*(S)$. Consider
the DG algebra $\B=\End(F)$. First we will prove the following
special case:
\begin{lemma} The pseudo-functors $\Def_-(\B)$ and $\Def_-^h(\B)$ are equivalent.\end{lemma}
\begin{proof} Take some $(S,\sigma)\in \Def^h(\B)$. Let $A$ be a strictly unital
minimal model of $\B$, and let $f:A\rarr \B$ be a strictly unital
$A_{\infty}$ quasi-isomorphism. By our assumption on
$Ext(E,E)\cong H(\B)$, $A$ is bounded above.

By Theorem \ref{invar_MC} there exists an object $\alpha\in
\cM\cC_{\mR}(A)$ such that $S\cong
\B\otimes_{f_{\mR}^*(\alpha)}\mR$. The DG $\mR^{op}$-modules
$\B\otimes_{f_{\mR}^*(\alpha)}\mR$ and
$f_*(\B\otimes_{f_{\mR}^*(\alpha)}\mR)$ are naturally identified.
Further, by Proposition \ref{invar_h_def1} we have natural
homotopy equivalence (in $A^{op}_{\mR^{op}}\text{-mod}_{\infty}$)
$$\gamma:\mR\otimes_{\alpha} A\rarr f_*(\B\otimes_{f^*(\alpha)}\mR).$$
Thus, it remains to prove that $$i^*(A\otimes_{\alpha} \mR)=\bfL
i^*(A\otimes_{\alpha} \mR).$$ We claim that $\mR\otimes_{\alpha}
A$ is h-projective. Indeed, since $A$ is bounded above and $\mR\in
\dgart_-$, this DG $\mR^{op}$-module has an increasing filtration
by DG $\mR^{op}$-submodules $A^{\geq i}\otimes \mR$ with
subquotients being free DG $\mR^{op}$-modules $A^i\otimes\mR$.
Thus $A\otimes_{\alpha} \mR$ satisfies property (P) as DG
$\mR^{op}$-module and hence is h-projective. Lemma is proved.
\end{proof}
Now take some $(S,id)\in \Def^h(F)$. We claim that $S$ is
h-projective. Recall the DG functor
$$\Sigma_{\mR}:(\B\otimes\mR)^{op}\text{-mod}\rarr \A_{\mR}^{op}\text{-mod},\quad \Sigma(M)=M\otimes_{\B} F.$$
From I, Proposition 9.2 e) we know that $S\cong \Sigma_{\mR}(S')$
for some $(S',id)\in \Def^h(\B)$. By the above Lemma and I,
Proposition 11.2, DG $(\B\otimes\mR)^{op}$-module $S'$ is
h-projective. Since the DG functor $\Sigma_{\mR}$ preserves
h-projectives, it follows that $S$ is also h-projective. Theorem
is proved.
\end{proof}

The next proposition is the analogue of Proposition
\ref{codef(B)=codef(E)} for the pseudo-functor $\DEF _-$. Note
that here we do not need boundedness assumptions on the
h-projective DG module.

\begin{prop}\label{def(B)=def(E)} Let $E\in \A^{op}\text{-mod}$ be a DG module. Suppose that the following
conditions hold:

a) $\Ext^{-1}(E,E)=0$;

b) the graded algebra $\Ext(E,E)$ is bounded above.

Put $\B=\bfR \Hom(E,E)$. Then pseudo-functors $\DEF_-(\B)$ and
$\DEF_-(E)$ are equivalent.
\end{prop}

\begin{proof} Take some h-projective $F$ quasi-isomorphic to
$E$ and replace $\DEF_-(E)$ by the equivalent pseudo-functor
$\DEF_-(F)$. We may assume that $\B=\End(F)$.

By I, Proposition 9.2 e) the DG functor
$\Sigma=\Sigma^{F}:\B^{op}\text{-mod}\rarr \A^{op}\text{-mod}$,
$\Sigma(N)=N\otimes_{\B} F$ induces an equivalence of
pseudo-functors $$\Def^h(\Sigma):\Def^h(\B)\rarr \Def^h(F).$$ By
Lemma \ref{def=def^h} we have that the pseudo-functors $\Def_-(F)$
and $\Def_-^h(F)$ (resp. $\Def_-(\B)$ and $\Def_-^h(\B)$) are
equivalent. We conclude that $\Sigma$ also induces an equivalence
of pseudo-functors
$$\Def_-(\Sigma):\Def_-(\B)\rarr \Def_-(F).$$
Let us prove that it extends to an equivalence
$$\DEF_-(\Sigma):\DEF_-(\B)\rarr \DEF_-(F).$$
Let $\mR,\mQ\in \dgart_-$, $M\in 1\text{-}\Hom(\mR,\mQ)$. We need
to show that the functorial diagram

$$
\begin{CD}
\DEF_{\mR}(\B) @>\DEF_{\mR}(\Sigma)>> \DEF_{\mR}(F)\\
@VM^*VV  @VM^*VV\\
\DEF_{\mQ}(\B) @>\DEF_{\mQ}(\Sigma)>> \DEF_{\mQ}(F).
\end{CD}
$$
commutes. This follows from the natural isomorphism
$$N\otimes_{\B}F\otimes_{\mR}M\cong N\otimes_{\mR} M\otimes_{\B} F.$$
\end{proof}

\section{Deformation pseudo-functor $\DEF$ for an augmented
$A_{\infty}$-algebra}

Let $A$ be an augmented $A_{\infty}$-algebra. We are going to
define the pseudo-functor $\DEF(A):2'\text{-}\dgart\rarr {\bf
Gpd}$.

Let $\mR$ be an artinian DG algebra. An object of the groupoid
$\DEF_{\mR}(A)$ is a pair $(S,\sigma)$, where $S\in
D_{\infty}(\bar{A}^{op}_{\mR^{op}})$, and $\sigma$ is an
isomorphism (in $D_{\infty}(\bar{A}^{op})$)
$$\sigma:\bfL i^*(S)\rarr A.$$

A morphism $f:(S,\sigma)\rarr (T,\tau)$ in $\DEF_{\mR}(A)$ is an
isomorphism (in $D(\bar{A}^{op}_{\cR^{op}})$) $f:S\rarr T$ such
that
$$\tau\circ \bfL i^*(f)=\sigma.$$

This defines the pseudo-functor $\DEF(A)$ on objects. Further, let
$(M,\theta)\in 1\text{-}\Hom(\mR,\mQ)$. Define the corresponding
functor
$$M^*:\DEF_{\mR}(A)\rarr \DEF_{\mQ}(A)$$ as follows. For an
object $(S,\sigma)\in \DEF_{\mR}(A)$ put
$$M^*(S)=S\otimes_{\mR}^{\bfL} M\in D_{\infty}(\bar{A}_{\mQ^{op}}^{op}).$$

Then we have natural isomorphisms in $D_{\infty}(\bar{A})$:
$$M^*(S)\otimes_{\mQ}^{\bfL} k=S\otimes_{\mR}^{\bfL} (M\otimes_{\mQ}^{\bfL}k)\cong S\otimes_{\mR}^{\bfL} k\cong A$$
(the second isomorphism is induced by $\theta$). Thus, $M^*$ is a
functor form $\DEF_{\mR}(A)$ to $\DEF_{\mQ}(A)$.

If $\mQ'$ is another artinian DG algebra and $(M',\theta')\in
1\text{-}\Hom(\mQ,\mQ')$ then there is a natural isomorphism of
functors
$$(M'\otimes_{\mQ}^{\bfL} M)^*\cong M'^*\cdot M^*.$$

Further, if $f\in 2\text{-}\Hom((M,\theta),(M,\theta_1))$ is a
$2$-morphism between objects $(M,\theta),(M,\theta_1)\in
1\text{-}\Hom(\mR,\mQ)$ then it induces an isomorphism between
corresponding functors $M^*\rarr M_1^*$.

Thus we obtain a pseudo-functor $\DEF(A):2\text{-}\dgart\rarr {\bf
GPd}$. We denote by $\DEF_-(A)$ its restriction to the
sub-$2$-category $2\text{-}\dgart_-$.

\begin{prop}\label{def-invar} Let $A$ be an augmented $A_{\infty}$-algebra and $U(A)$ its universal DG algebra. Then there is a natural
equivalence of pseudo-functors $\DEF(U(A))\cong  \DEF(A)$.
\end{prop}
\begin{proof}The proof is the same as of Proposition \ref{codef-invar} and we omit it.\end{proof}

\begin{cor}\label{def(A)=def(B)} Let $A$ be an augmented $A_{\infty}$-algebra and let
$\cB$ be a DG algebra quasi-isomorphic to $A$. Then the
pseudo-functor $\DEF(A)$ and $\DEF(\cB)$ are equivalent.
\end{cor}
\begin{proof} Indeed, by Proposition \ref{codef-invar} the
pseudo-functors $\DEF(A)$ and $\DEF(U(A))$ are equivalent, and by
Corollary \ref{def(B)=def(C)} the pseudo-functors $\DEF(U(A))$ and
$\DEF(\cB)$ are equivalent.
\end{proof}

\begin{cor}\label{def^h} Let $A$ be an admissible $A_{\infty}$-algebra, and $\mR$ be an artinian negative DG algebra. Then for any
$(S,\sigma)\in \DEF_{\mR}(A)$ there exists an $\alpha\in
\cM\cC_{\mR}(A)$ such that the pair $(T,id)$, where
$T=A\otimes_{\alpha} \cR$, defines an object of $\DEF_{\mR}(A)$
which is isomorphic to $(S,\sigma)$.\end{cor}
\begin{proof} This follows easily from Proposition \ref{def-invar} and the proof of Lemma
\ref{def=def^h} in the case $\B=U(A)$.\end{proof}

\section{Comparison of pseudo-functors $\coDEF _-(E)$ and $\DEF
_-(E)$} \label{comparison}

We have proved in I, Corollary 11.9 that under some conditions on
$E$ the pseudo-functors $\coDef _-(E)$ and $\Def _-(E)$ from
$\dgart _-$ to $\bf{Gpd}$ are equivalent. Note that we cannot
speak about an equivalence of pseudo-functors $\coDEF _-(E)$ and
$\DEF _-(E)$ since they are defined on different bicategories. So
our first goal is to establish an equivalence of the bicategories
$2\text{-}\adgalg$ and $2^\prime \text{-}\adgalg$ in the following
sense: we will construct pseudo-functors
$$\cD :2\text{-}\adgalg \to 2^\prime \text{-}\adgalg,$$
$$\cD^\prime :2^\prime \text{-}\adgalg \to 2\text{-}\adgalg,$$
which have the following properties

1) $\cD$ (resp. $\cD^\prime $) is the identity on objects;

2) for each $\cB ,\cC \in Ob (2\text{-}\adgalg)$ they define
mutually inverse equivalences of groupoids
$$\cD :\Hom _{2\text{-}\adgalg}(\cB ,\cC)\to \Hom _{2^\prime
\text{-}\adgalg}(\cB ,\cC),$$
$$\cD^\prime : \Hom _{2^\prime
\text{-}\adgalg}(\cB ,\cC)\to \Hom _{2\text{-}\adgalg}(\cB
,\cC).$$

Fix augmented DG algebras $\cB,$ $\cC$ and let $M$ be a DG $\cC
\otimes \cB ^{op}$-module. Define the DG $\cB\otimes \cC
^{op}$-module $\cD (M)$ as
$$\cD (M):=\bR \Hom _{\cC}(M,\cC).$$

Further, let $N$ be a DG $\cB \otimes \cC ^{op}$-module. Define
the DG $\cB ^{op}\otimes \cC$-module $\cD^\prime (N)$ as
$$\cD ^\prime (N)=\bR \Hom _{\cC ^{op}}(N,\cC).$$

\begin{prop}\label{equivalence_D} The operations $\cD$, $\cD^\prime$  as above induces the
pseudo-functors
$$\cD :2\text{-}\adgalg \to 2^\prime \text{-}\adgalg,$$
$$\cD^\prime :2^\prime \text{-}\adgalg \to 2\text{-}\adgalg,$$
so that the properties 1) and 2) hold.
\end{prop}

\begin{proof} To simplify the notation denote by $\Hom (-,-)$ and
$\Hom ^\prime (-,-)$ the morphisms in the bicategories
$2\text{-}\adgalg$ and $2^\prime \text{-}\adgalg$ respectively.

We will prove that for augmented DG algebras $\cB$ and $\cC$  we
have a (covariant) functor
$$\cD :\Hom (\cB ,\cC)\to \Hom ^\prime (\cB ,\cC),$$
and the functorial diagram
$$\begin{array} {ccccc}
\Hom (\cB _1,\cB _2) & \times & \Hom (\cB _2,\cB _3) & \to & \Hom
(\cB _1 ,\cB _3)\\
\cD \downarrow & & \cD \downarrow & & \cD \downarrow \\
\Hom ^\prime(\cB _1,\cB _2) & \times & \Hom ^\prime(\cB _2,\cB _3) &
\to & \Hom ^\prime(\cB _1 ,\cB _3)
\end{array}$$
commutes for every triple of augmented algebras $\cB _1,$ $\cB _2,$
$\cB _3$.

Let $(M, \theta )\in 1\text{-}\Hom (\cB ,\cC)$. Choose a
quasi-isomorphism $f:\cC \to \nu _*M$ of DG $\cC$-modules. It
induces the quasi-isomorphism
$$\cD (f):\nu _*\cD (M)\to \bR \Hom _{\cC}(\cC, \cC)=\cC$$
of DG $\cC ^{op}$-modules. Moreover, we claim that the
quasi-isomorphism $\theta :k\stackrel{\bL}{\otimes }_{\cC}M\to k$
induces a quasi-isomorphism
$$\cD (\theta):\cD (M)\stackrel{\bL}{\otimes }_{\cC}k\to k^*=k.$$
Indeed, we may and will assume that the DG $\cC \otimes \cB
^{op}$-module $M$ is h-projective. Then by I, Lemma 3.23 it is
also h-projective as a DG $\cC$-module. Therefore by Lemma
\ref{h-projectives} a) below
$$\cD (M)\stackrel{\bL}{\otimes }_{\cC}k =\Hom _{\cC}(M,\cC)\otimes
_{\cC}k.$$ Note that the obvious morphism of DG $\cB$-modules
$$\delta :\Hom _{\cC}(M,\cC)\otimes
_{\cC}k\to \Hom _{\cC}(M,k)$$ is a quasi-isomorphism. Indeed, the DG
$\cC$-module $M$ is homotopy equivalent to $\cC$. Hence it suffices
to check that $\delta$ is an isomorphism when $M=\cC$, which is
obvious, since both sides are equal to $k$. Now notice the obvious
canonical isomorphisms
$$\Hom _{\cC}(M,k)=\Hom _k(k\otimes _{\cC}M,k)= (k\otimes
_{\cC}M)^*\stackrel{\theta ^*}{\longleftarrow}k^*=k.$$

Thus indeed, $(\cD (M), \cD (\theta))$ is an object in $\Hom ^\prime
(\cB ,\cC)$ and therefore we have a (covariant) functor
$$\cD :\Hom (\cB ,\cC)\to \Hom ^\prime (\cB ,\cC).$$

Let now $\cB _1, \cB _2, \cB _3 \in Ob (2\text{-}\adgalg)$ and
$M_1\in 1\text{-}\Hom (\cB _1,\cB _2)$, $M_2\in 1\text{-}\Hom (\cB
_2,\cB _3)$. Then
$$M_2\stackrel{\bL}{\otimes }_{\cB _2}M_1\in 1\text{-}\Hom (\cB
_1,\cB _3), \quad \text{and} \quad \cD (M_1)\stackrel{\bL}{\otimes
}_{\cB _2}\cD(M_2)\in 1\text{-}\Hom ^\prime(\cB _1,\cB _3).$$ We
claim that the DG $\cB _1\otimes \cB ^{op}_3$-modules
$$\cD (M_2\stackrel{\bL}{\otimes }_{\cB _2}M_1) \quad \text{and}
\quad \cD (M_1)\stackrel{\bL}{\otimes }_{\cB _2}\cD(M_2)$$ are
canonically quasi-isomorphic.

Indeed, we may and will assume that $M_1$ and $M_2$ are
h-projective as DG $\cB _2\otimes \cB _1^{op}$- and $\cB _3\otimes
\cB _2^{op}$-modules respectively. Then by Lemma
\ref{h-projectives} below it suffices to prove that the morphism
of DG $\cB _1\otimes \cB _3^{op}$-modules
$$\epsilon :\Hom _{\cB _2}(M_1,\cB _2)\otimes _{\cB _2}\Hom _{\cB
_3}(M_2,\cB _3)\to \Hom _{\cB _3}(M_2\otimes _{\cB _2}M_1,\cB _3)$$
defined by
$$\epsilon(f\otimes g)(m_2\otimes
m_1):=(-1)^{\bar{f}(\bar{g}+\bar{m}_2)}g(m_2f(m_1))$$ is a
quasi-isomorphism. To prove that $\epsilon$ is a quasi-isomorphism
we may replace the DG $\cB _2$-module $M_1$ by $\cB _2$. Then
$\epsilon$ is an isomorphism.

Thus, the operation $\cD$ induces a pseudo-functor

$$\cD :2\text{-}\adgalg \to 2^\prime \text{-}\adgalg.$$

Analogously, the operation $\cD^\prime$ induces a pseudo-functor
$$\cD^\prime :2^\prime \text{-}\adgalg \to 2\text{-}\adgalg.$$ It
is clear that for $M\in 1\text{-}\Hom(\cB,\cC)$ (resp. $N\in
1\text{-}\Hom^\prime (\cB,\cC)$) the canonical morphism $M\to
\cD^\prime \cD(M)$ (resp. $N\to \cD\cD^\prime (M)$) is an
isomorphism. Thus, the compositions $\cD^\prime \cD$ and
$\cD\cD^\prime $ are equivalent to the identity.

Proposition is proved.

\end{proof}

\begin{cor}\label{h_B=h_B'D} For any augmented DG algebra $\cB$ the pseudo-functor
$\cD :2\text{-}\adgalg \to 2^\prime \text{-}\adgalg$ induces a
morphism of pseudo-functors
$$h_{\cB}\to h^\prime _{\cB}\cdot \cD,$$
which is an equivalence.

Similarly, the pseudo-functor $\cD ^\prime :2^\prime
\text{-}\adgalg\to 2\text{-}\adgalg$ induces an equivalence of
pseudo-functors
$$h^\prime_{\cB}\to h _{\cB}\cdot \cD^\prime.$$
\end{cor}

\begin{proof} This is clear.
\end{proof}

\begin{lemma}\label{h-projectives} Let $\cB_1, \cB_2, \cB _3 \in Ob(2\text{-}\adgalg)$, $M_1\in
1\text{-}\Hom (\cB _1 ,\cB _2)$, $M_2\in 1\text{-}\Hom (\cB _2
,\cB _3)$. Assume that $M_1$ and $M_2$ are h-projective as DG $\cB
_2\otimes \cB _1^{op}$- and $\cB _3\otimes \cB _2^{op}$-modules
respectively. Then

a) The DG $\cB _2^{op}$-module $\Hom _{\cB _2}(M_1,\cB_2)$ is
h-projective.

b) The DG $\cB _3$-module $M_2\otimes _{\cB _2}M_1$ is h-projective.
\end{lemma}

\begin{proof} a). Since $M_1$ is h-projective as a DG $\cB_2 \otimes
\cB _1 ^{op}$-module, it is also such as a DG $\cB_2$-module (I,
Lemma 3.23). We denote this DG $\cB_2$-module again by $M_1$.

Choose a quasi-isomorphism of DG $\cB_2$-modules $f:\cB_2 \to
M_1.$ This is a homotopy equivalence since both $\cB_2$ and $M_1$
are h-projective. Thus it induces a homotopy equivalence of DG
$\cB_2^{op}$-modules
$$f^*:\Hom _{\cB_2}(\cB_2 ,\cB_2)\to \Hom _{\cB_2}(M_1,\cB_2).$$
But the DG $\cB_2^{op}$-module $\Hom _{\cB_2}(\cB_2,\cB_2)=\cB_2$
is h-projective. Hence so is $\Hom _{\cB_2}(M _1,\cB_2)$.

b). The proof is similar. Namely, the DG $\cB _3$-module $M_2\otimes
_{\cB _2}M_1$ is homotopy equivalent to $M_2\otimes _{\cB _2}\cB
_2=M_2$, which is homotopy equivalent to $\cB _3$.
\end{proof}

\begin{theo}\label{def=codef1} Assume that the DG $\cA ^{op}$-module $E$ has the
following properties.

i) $\Ext ^{-1}(E,E)=0.$

ii) There exists a bounded above h-projective or h-injective DG
$\cA ^{op}$-module $P$ quasi-isomorphic to $E$.

iii) There exists a bounded below h-projective or h-injective  DG
$\cA ^{op}$-module $I$ which is quasi-isomorphic to $E$.

Then the pseudo-functors $\coDEF _-(E)$ and $\DEF _(E)-\cdot \cD$
from $2\text{-}\dgart _-$ to $\bf{Gpd}$ are equivalent.

Hence also the pseudo-functors $\DEF _-(E)$ and $\coDEF _(E)-\cdot
\cD ^\prime$ from $2^\prime\text{-}\dgart _-$ to $\bf{Gpd}$ are
equivalent.
\end{theo}

\begin{proof} Let $\cR \in \dgart _-$. Recall (I, Theorem 11.13)
the DG functor
$$ \epsilon _{\cR}:\cA ^{op}_{\cR}\text{-mod}\to \cA
^{op}_{\cR}\text{-mod}$$ defined by
$$\epsilon _{\cR}(M)=M\otimes _{\cR}\cR ^*;$$ and the corresponding derived
functor
$$\bL \epsilon _{\cR}:D(\cA ^{op}_{\cR})\to D(\cA
^{op}_{\cR}).$$

We know (I, Theorem 11.13) that under the assumptions i), ii), iii)
this functor induces an equivalence of groupoids
$$\bL \epsilon _{\cR}:\Def _{\cR}(E)\to \coDef _{\cR}(E).$$

Let now $\cQ \in \dgart _-$ and $M\in 1\text{-}\Hom (\cR ,\cQ)$. It
suffices to prove that the functorial diagram
$$\begin{array}{ccc}
\Def _{\cR}(E) & \stackrel{\bL \epsilon _{\cR}}{\to} & \coDef
_{\cR}(E)\\
\cD (M)^*\downarrow & & \downarrow M^!\\
\Def _{\cQ}(E) & \stackrel{\bL \epsilon _{\cQ}}{\to } & \coDef
_{\cQ}(E)
\end{array}
$$
naturally commutes.

Choose a bounded above h-projective or h-injective $P$
quasi-isomorphic to $E$. By I, Theorem 11.6 a) the groupoids $\Def
_{\cR}(E)$ and $\Def _{\cR}^{\h}(P)$ are equivalent. Hence given
$(S,\id)\in \Def _{\cR}^{\h}(P)$ it suffices to prove that there
exists a natural isomorpism of objects in $D(\cA _{\cQ}^{op})$
$$M^!\cdot \bL \epsilon _{\cR}(S)\simeq \bL \epsilon _{\cQ}\cdot \cD
(M)^*(S),$$ i.e.
$$\bR \Hom _{\cR ^{op}}(M,S\stackrel{\bL}{\otimes }_{\cR}\cR ^*)\simeq
S\stackrel{\bL }{\otimes }_{\cR}\bR \Hom
_{\cQ}(M,\cQ)\stackrel{\bL}{\otimes }_{\cQ}\cQ^*.$$ We may and
will assume that the DG $\cQ \otimes \cR ^{op}$-module $M$ is
h-projective. In the proof of I, Lemma 11.7 we showed that  the DG
$\cA ^{op}_{\cR}$-module $S$ is h-projective as a DG $\cR
^{op}$-module. Therefore by Lemma \ref{h-projectives} a) it
suffices to prove that the morphism of DG $\cA
^{op}_{\cQ}$-modules
$$\eta :S\otimes _{\cR}\Hom _{\cQ}(M,\cQ)\otimes _{\cQ}\cQ ^*\to
\Hom _{\cR ^{op}}(M,S\otimes _{\cR}\cR ^*)$$ defined by
$$\eta (s\otimes f\otimes g)(m)(r)=sg(f(mr))$$
is a quasi-isomorphism.

It suffices to prove that $\eta$ is a quasi-isomorphism of DG $\cQ
^{op}$-modules. Notice that just the $\cR^{op}$-module structure
on $S$ is important for us. Furthermore we may assume that $S$
satisfies property (P) as DG $\cR^{op}$-module (I, Definition
3.2). Thus it suffices to prove that $\eta$ is a quasi-isomorphism
if $S=\cR$. Then
$$\eta :\Hom _{\cQ}(M,\cQ)\otimes _{\cQ}\cQ ^*\to
\Hom _{\cR ^{op}}(M,\cR ^*).$$ We have the canonical isomorphisms
$$\Hom _{\cR ^{op}}(M,\Hom _k(\cR ,k))=\Hom _{k}(M\otimes
_{\cR}\cR,k)=M^*.$$ Also, since the DG $\cQ ^{op}$-module $M$ is
homotopy equivalent to $\cQ$, we have the homotopy equivalences
$$\Hom _{\cQ}(M,\cQ)\otimes _{\cQ}\cQ ^*\simeq \Hom _{\cQ}(\cQ,\cQ)\otimes _{\cQ}\cQ
^* \simeq \cQ ^*\simeq M^*.$$
\end{proof}

The next theorem is closely related to the previous one. It
asserts the stronger statement in the case when $E$ is a DG
algebra considered as a DG module over itself.

\begin{theo}\label{def=codef} Let $\B$ be a DG algebra.
Suppose that the following conditions hold:

a) $H^{-1}(\B)=0$;

b) the cohomology algebra $H(\B)$ is bounded above and bounded
below. Then the pseudo-functors $\coDEF_-(\B)$ and
$\DEF_-(\B)\cdot \D$ from $2\text{-}\dgart_-$ to ${\bf Gpd}$ are
equivalent.
\end{theo}

\begin{proof} Let $\mR$ be a negative artinian DG algebra. Recall the DG functors
$$\epsilon_{\mR}:(\B\otimes \mR)^{op}\text{-mod}\rarr (\B\otimes \mR)^{op}\text{-mod},\quad \epsilon_{\mR}(M)=
M\otimes_{\mR}\mR^*,$$
$$\eta_{\mR}:(\B\otimes \mR)^{op}\text{-mod}\rarr (\B\otimes \mR)^{op}\text{-mod},\quad \eta_{\mR}(M)=
\Hom_{\mR}(\mR^*,M).$$ By I, Proposition 4.7 they induce
quasi-inverse equivalences
$$\epsilon_{\mR}:\Def_{\mR}^h(\B)\rarr \coDef_{\mR}^h(\B),$$
$$\eta_{\mR}:\coDef_{\mR}^h(\B)\rarr \Def_{\mR}^h(\B).$$

By Theorem \ref{def=def^h} the pseudo-functors $\Def_-(\B)$ and
$\Def_-^h(\B)$ are equivalent. By Lemma \ref{codef=codef^h} the
pseudo-functors $\coDef_-(\B)$ and $\coDef_-^h(\B)$. It follows
that the derived functors $\bfL \epsilon_{\mR}$, $\bfR \eta_{\mR}$
induce mutually inverse equivalences
$$\bfL \epsilon_{\mR}:\Def_{\mR}(\B)\rarr \coDef_{\mR}(\B),$$
$$\bfR \eta_{\mR}:\coDef_{\mR}(\B)\rarr \Def_{\mR}(\B).$$

Let now $\mQ\in \dgart_-$ and $M\in 1\text{-}\Hom(\mR,\mQ)$. It
suffices to prove that the functorial diagram
$$
\begin{CD}
\Def_{\mR}(\B) @>\bfL \epsilon_{\mR}>> \coDef_{\mR}(\B)\\
@V\D(M)^*VV  @VM^!VV\\
\coDef_{\mQ}(\B) @>\bfL \epsilon_{\mQ}>> \coDef_{\mQ}(\B)
\end{CD}
$$
naturally commutes. This fact is absoulutely analogous to the
analogous fact from the proof of the previous theorem.\end{proof}

\part{Pro-representability theorems}

\section{Pro-representability of the pseudo-functor $\coDEF_-$}

The next theorem claims that under some conditions on the DG
algebra $\cC$ that the functor $\coDEF_-(\cC)$ is
pro-representable.

\begin{theo}\label{pro-repr1} Let $\cC$ be a DG algebra such that the cohomology algebra
$H(\cC)$ is admissible finite-dimensional. Let $A$ be a strictly
unital minimal model of $\cC$. Then the pseudo-functor
$\coDEF_-(\cC)$ is pro-representable by the DG algebra
$\hat{S}=(B\bar{A})^*$. That is, there exists an equivalence of
pseudo-functors $\coDEF_-(\cC)\simeq h_{\hat{S}}$ from
$2\text{-}\dgart_-$ to {\bf Gpd}.
\end{theo}

As a corollary, we obtain the following

\begin{theo}\label{pro-repr2} Let $E\in \A^{op}\text{-mod}$. Assume that the following
conditions hold:

a) the graded algebra $\Ext(E,E)$ is admissible
finite-dimensional;

b) $E$ is quasi-isomorphic to a bounded below $F$ which is
h-projective or h-injective.

Then the pseudo-functor $\coDEF_-(E)$ is pro-representable by the
DG algebra $\hat{S}=(B\bar{A})^*$, where $A$ is a strictly unital
minimal model of $\bR \Hom(E,E)$.
\end{theo}
\begin{proof} By Proposition \ref{codef(B)=codef(E)} the pseudo-functors $\coDEF_-(E)$ and
$\coDEF_-(\bR \Hom(E,E))$ are equivalent. So it remains to apply
Theorem \ref{pro-repr1}.\end{proof}

\begin{proof} Note that we have natural quasi-isomorphism of DG algebras $U(A)\rarr \cC$, hence
$\coDEF_-(\cC)\simeq \coDEF_-(U(A))$. Further, by Proposition
\ref{codef-invar} we have $\coDEF_-(U(A))\simeq \coDEF_-(A)$. We
will construct an equivalence of pseudo-functors
$\Theta:h_{\hat{S}}\rarr \coDEF_-(A)$.

Consider the $A_{\infty}$ $\bar{A}^{op}_{\hat{S}^{op}}$-module
$B\bar{A}\otimes A$. Choose a quasi-isomorphism $B\bar{A}\otimes
A\rarr J$, where $J$ is an h-injective $A_{\infty}$
$\bar{A}^{op}_{\hat{S}^{op}}$-module. Note that $J$ is also
h-injective as a DG $\hat{S}^{op}$-module.

Given an artinian DG algebra $\mR$ and a $1$-morphism
$(M,\theta)\in 1\text{-}\Hom(\hat{S},\mR)$ we define

$$\Theta(M):=\Hom_{\hat{S}^{op}}(M,J).$$

We have $\bfR \Hom_{\mR^{op}}(k,\Hom_{\hat{S}^{op}}(M,J))=\bfR
\Hom_{\hat{S}^{op}}(k\otimes_{\mR}^{\bfL} M,J)$. Hence the
quasi-isomorphism $\theta:k\otimes_{\mR}^{\bfL} M\rarr k$ induces
a quasi-isomorphism

$$\bfR \Hom_{\mR^{op}}(k,\Theta(M))\simeq \bfR \Hom_{\hat{S}^{op}}(k,J)=\Hom_{\hat{S}^{op}}(k,J),$$
and by Proposition \ref{admissible_prelim} the last term is
canonically quasi-isomorphic to $A$ as an $A_{\infty}$
$A^{op}$-module.

If we are given with another artinian DG algebra $\mQ$ and a
$1$-morphism $(N,\delta)\in 1\text{-}\Hom(\mR,\mQ)$, then the
object $\Theta(N\otimes_{\mR}^{\bfL} M)$ is canonically
quasi-isomorphic to the object $\bfR \Hom(N,\Theta(M))$. Thus,
$\Theta$ is a morphism of pseudo-functors.

It remains to prove that for each $\mR\in 2\text{-}\dgart_-$ the
induced functor $\Theta_{\mR}:1\text{-}\Hom(\hat{S},\mR)\rarr
\coDEF_{\mR}(A)$ is an equivalence of groupoids. So fix a DG
algebra $\mR\in$ $2\text{-}\dgart_-$.

{\noindent}{\bf Surjective on isomorphism classes.} Let
$(S,\sigma)$ be an object of $\coDEF_{\mR}(A)$. By Corollary
\ref{codef^h}, there exists a morphism of DG algebras
$\phi:\hat{S}\rarr \mR$ such that the pair $(T,id)$, where
$T=\Hom_{\hat{S}^{op}}(\mR,B\bar{A}\otimes A)$, defines an object
of $\coDEF_{\mR}(A)$ which is isomorphic to $(S,\sigma)$. Further,
by Proposition \ref{hom=rhom1} the morphism
$\Hom_{\hat{S}^{op}}(\mR,B\bar{A}\otimes A)\rarr
\Hom_{\hat{S}^{op}}(\mR,J)$ is quasi-isomorphism. Therefore, the
object $(T,id)$ is isomorphic to $\Theta(M)$, where $M=\mR$ is DG
$\hat{S}^{op}\otimes \mR$-module via the homomorphism $\phi$.

{\noindent}{\bf Full and Faithful.} Consider the above $\Theta$ as
a contravariant DG functor from $\mR\otimes
\hat{S}^{op}\text{-mod}$ to
$\bar{A}^{op}_{\mR^{op}}\text{-mod}_{\infty}$. Define the
contravariant DG functor
$\Phi:\bar{A}^{op}_{\mR^{op}}\text{-mod}_{\infty}\rarr \mR\otimes
\hat{S}^{op}\text{-mod}$ defined by the similar formula:

$$\Phi(N)=\Hom_{\bar{A}^{op}}(N,J).$$

These DG functors induce the corresponding DG functors between
derived categories $$\Theta:D(\mR\otimes \hat{S}^{op})\rarr
D_{\infty}(\bar{A}^{op}_{\hat{S}^{op}}),\quad
\Phi:D_{\infty}(\bar{A}_{\hat{S}^{op}}^{op})\rarr D(\mR\otimes
\hat{S}^{op}).$$ Denote by $\langle k\rangle\subset D(\mR\otimes
\hat{S}^{op})$ and $\langle A\rangle\subset
D_{\infty}(\bar{A}_{\hat{S}^{op}}^{op})$ the triangulated
envelopes of the DG $\mR\otimes \hat{S}^{op}$-module $k$ and
$A_{\infty}$ $\bar{A}^{op}_{\mR^{op}}$-module $A$ respectively.

\begin{lemma} The functors $\Theta$ and $\Phi$ induce mutually inverse anti-equivalences of the
triangulated categories $\langle k\rangle$ and $\langle
A\rangle$.\end{lemma}
\begin{proof}For $M\in \mR\otimes \hat{S}^{op}\text{-mod}$, and
$N\in \bar{A}^{op}_{\mR^{op}}\text{-mod}_{\infty}$ we have the
functorial closed morphisms $$\beta_M:M\rarr \Phi(\Theta(M)),\quad
\beta_M(x)_1(f)=(-1)^{|f||x|}f(x),\quad \beta_M(x)_n=0\text{ for
}n\geq 2;$$
$$\gamma_N:N\rarr \Theta(\Phi(N)),\quad (\gamma_N)_n(a_1,\dots,a_{n-1},y)(f)=(-1)^{n(|a_1|+\dots+|a_{n-1}|+|y|)} f_n(a_1,\dots,a_{n-1},y).$$

By Proposition \ref{admissible_prelim} the $A_{\infty}$
$\bar{A}^{op}_{\hat{S}^{op}}$-module $\Theta(k)$ is
quasi-isomorphic to $A$. Further, $\Phi(A)$ is quasi-isomorphic to
$J$ and hence to $k$. Therefore, $\beta_k$ and $\gamma_A$ are
quasi-isomorphisms, and Lemma is proved.
\end{proof}

Note that for $(M,\theta)\in 1\text{-}\Hom(\hat{S},\mR)$ (resp.
for $(S,\sigma)\in \coDEF_{\mR}(A)$) $M\in \langle k\rangle$
(resp. $S\in \langle A\rangle$). Hence the functor
$\Theta_{\mR}:1\text{-}\Hom(\hat{S},\mR)\rarr \coDEF_{\mR}(A)$ is
fully faithful. This proves the theorem.\end{proof}

\section{Pro-representability of the pseudo-functor $\DEF _-$}

Pro-representability Theorems \ref{pro-repr1} and \ref{pro-repr2}
imply analogous results for the pseudo-functor $\DEF_-$. Namely,
we have the following Theorems.

\begin{theo} \label{pro-repr3} Let $\cC$ be a DG algebra such that the cohomology algebra $H(\cC)$ is admissible finite-dimensional.
Let $A$ be a strictly unital minimal model of $\cC$. Then the
pseudo-functor $\DEF_-(\cC)$ is pro-representable by the DG
algebra $\hat{S}=(B\bar{A})^*$. That is, there exists an
equivalence of pseudo-functors $\DEF_-(\cC)\simeq h_{\hat{S}}'$
from $2'\text{-}\dgart_-$ to {\bf Gpd}.
\end{theo}
\begin{proof} By Theorem \ref{pro-repr1} we have the equivalence
$$\coDEF_-(\cC)\simeq h_{\hat{S}}$$
of pseudo-functors from $2\text{-}\dgart$ to {\bf Gpd}.

By Theorem \ref{def=codef} we have the equivalence
$$\coDEF_-(\cC)\simeq \DEF_-(\cC)\cdot \D$$
of pseudo-functors from $2\text{-}\dgart$ to {\bf Gpd}.

Further, by Corollary \ref{h_B=h_B'D}
$$h_{\hat{S}}\simeq h_{\hat{S}}'\cdot \D.$$
Hence $\DEF_-(\cC)\cdot \D\simeq h_{\hat{S}}'\cdot \D$ and
therefore $$\DEF_-(\cC)\simeq h_{\hat{S}}'.$$\end{proof}

We get the following corollary.

\begin{theo}\label{pro-repr4} Let $E\in \A^{op}\text{-mod}$. Assume that the graded algebra
$\Ext(E,E)$ is admissible finite-dimensional. Then the
pseudo-functor $\DEF_-(E)$ is pro-representable by the DG algebra
$\hat{S}=(B\bar{A})^*$, where $A$ is a strictly unital minimal
model of the DG algebra $\bfR \Hom(E,E)$.
\end{theo}
\begin{proof} Indeed, by Proposition \ref{def(B)=def(E)} the
pseudo-functors $\DEF_-(E)$ and $\DEF_-(\bfR \Hom(E,E))$ are
equivalent. And by Theorem \ref{pro-repr3} the pseudo-functors
$\DEF_-(\bfR \Hom(E,E))$ and $h_{\hat{S}}'$ are equivalent.
\end{proof}

We would like to mention here several examples.

\begin{example}\label{k-points} Let $X$ be a commutative scheme over $k$ of finite type, and let $x\in X(k)$ be a
regular $k\text{-}$point. Take the skyscraper sheaf $\cO_x\in
D^b_{coh}(X).$ Then one can show that
$\Ext^{\cdot}(\cO_x,\cO_x)\cong\Lambda(T_{x}X),$ and the DG
algebra $\bR\Hom(\cO_{x},\cO_{x})$ is formal. It follows that
$H^i(\hat{S})=0$ for $i\ne 0,$ and $H^0(\hat{S})\cong
k[[t_1,\dots,t_n]],$ where $n=\dim_x X.$\end{example}

\begin{example}\label{NJac}Let $X$ be a proper curve of genus $g$ over $k$ and $\cL\in D^b_{coh}(X)$ a line bundle over $X.$
Then $\Ext^0(\cL,\cL)=k,$ $\Ext^1(\cL,\cL)=k^g,$ and
$\Ext^i(\cL,\cL)=0$ for $i\ne 0,1.$ It follows that the DG algebra
$\hat{S}$ is concentrated in degree zero and is isomorphic to the
algebra of non-commutative power series in $g$
variables.\end{example}

\begin{example}\label{NGr}Let $V$ be a vector space of dimension $n,$ and let $W\subset V$ be a subspace of dimension $m,$
$1\leq m\leq n-1.$ Put $E=\cO_{\PP(W)}\in D^b_{coh}(\PP(V)).$ One
can show that the graded algebra $A=\Ext^{\cdot}(E,E)$ is
isomorphic to $\sum\limits_{0\leq i\leq
n-m}\Sym^i(W^{\vee})\otimes \Lambda^i(V/W).$ The later algebra can
be shown to be quadratic Koszul. Again, one can show that the DG
algebra $\bR\Hom(E,E)$ is formal. It follows that $H^i(\hat{S})=0$
for $i\ne 0,$ and $H^0(\hat{S})$ is a (completion of) Koszul dual
to $A.$ For $m\ne 1,$ we have that the algebra $H^0(\hat{S})$ is
non-commutative.
\end{example}

In the proof of Theorem \ref{pro-repr1} we showed that the bar
complex $B\bar{A}\otimes_{\tau_A} A$ is the "universal
co-deformation" of the $A_{\infty}$ $\bar{A}^{op}$-module $A$.
However, Theorem \ref{pro-repr3} is deduced from Theorem
\ref{pro-repr1} without finding the analogous "universal
deformation" of the $A_{\infty}$ $\bar{A}^{op}$-module $A$. We do
not know if this "universal deformation" exists in general (under
the assumptions of Theorem \ref{pro-repr3}). But we can find it
and hence give a direct proof of Theorem \ref{pro-repr3} if the
minimal model $A$ of $\cC$ satisfies an extra assumption (*)
below.

For the rest of this section we assume that the DG algebra $\cC$
has an augmented minimal model $A$.

\begin{defi}\label{cond_*} Let $A$ be an augmented $A_{\infty}$-algebra. Consider $k$ as a left $A_{\infty}$
$A$-module. We say that $A$ satisfies the condition (*) if the
canonical morphism
$$k\to \Hom _{\bar{A} ^{op}}(\Hom _{\bar{A}}(k,A),A)$$ of left $A_{\infty}$
$\bar{A}$-modules is a quasi-isomorphism.
\end{defi}

\begin{example} Let $A$ be an augmented $A_{\infty}$-algebra. If $k$ lies in $\Perf(A)$
then $A$ satisfies the condition (*).

In particular, suppose that $A$ is homologically smooth and
compact. That is, the diagonal $A_{\infty}$ $A-A$-bimodule $A$
lies in $\Perf(A-A)$ (smoothness), and $\dim H(A)<\infty$
(compactness). Then the $A_{\infty}$ $A$-module is perfect iff it
has finite-dimensional total cohomology. Thus, $k\in \Perf(A)$ and
$A$ satisfies the condition (*).
\end{example}

\begin{example}\label{Gorenstein_*} Let $A$ be an augmented $A_{\infty}$-algebra which is left and right
 Gorenstein of dimension $d$.This means that
$$\Ext ^p_{\bar{A}}(k,A)=\left\{ \begin{array}{ll}
k, & \text{if p=d}\\
0, & \text{otherwise,}\\
\end{array} \right. $$ and
$$\Ext ^p_{\bar{A}^{op}}(k,A)=\left\{ \begin{array}{ll}
k, & \text{if p=d}\\
0, & \text{otherwise.}\\
\end{array} \right. $$ Then $A$ satisfies the
condition (*).
\end{example}

For the rest of this section assume that $A$ is admissible,
finite-dimensional and satisfies the condition (*).

Denote by $\cE$ the $A_{\infty}$
$\bar{A}^{op}_{\hat{S}^{op}}$-module
$$\cE :=\Hom _{A}(k,A).$$
This $A_{\infty}$-module is isomorphic to $A\otimes \hat{S}$ as a
graded $(\hat{S}^{op})^{gr}$-module and can be given explicitly by
the formula

\begin{equation}\label{univ_def} m_n^{\cE}(m,a_1,\dots,a_{n-1})=
m_n^{\cM\cC_{\infty}^{\hat{S}}(A)(0,\dots,0,\tau_A)}(m,a_1\otimes
\one_{\hat{S}},\dots,a_{n-1}\otimes \one_{\hat{S}}).\end{equation}

\begin{remark}\label{CotimesS} The definition of the $A_{\infty}$-category $\cM\cC_{\infty}^{\hat{S}}(A)$ is the same as
if $\hat{S}$ would be artinian. It is correct because $\hat{S}$ is
complete in $\m$-adic topology and $A$ is finite-dimensional. In
the above formula $\tau_A$ is considered as an element of
$A\otimes \hat{S}=\Hom_k(B\bar{A},A)$. We denote the $A_{\infty}$
$\bar{A}^{op}_{\hat{S}^{op}}$-module $\cE$ by
$A\otimes_{\tau_A}\hat{S}$.
\end{remark}

We claim that $\cE$ is the "universal deformation" of $A$. This is
justified by Theorem \ref{pro-repr_special} below. Let us start
with a few lemmas.

%The following lemma is the analogue of Lemma \ref{CotimesR}.

%\begin{lemma}\label{BotimesG^*} Let $\cB$ be a finite dimensional DG algebra and $\cG$
%be a DG coalgebra. Then the DG algebras $\cB \otimes \cG ^*$ and
%$\Hom _k(\cG ,\cB)$ are isomorphic.
%\end{lemma}

%\begin{proof} The isomorphism $\nu :\cB \otimes \cG ^*\to \Hom (\cG
%,\cB)$ is given by the formula
%$$\nu (b\otimes f)(g)=bf(g),$$
%for $b\in \cB, f\in \cG ^*, g\in \cG$.
%\end{proof}

%\begin{remark}\label{CotimesS} The inclusion map $k\hookrightarrow \cC$ induces an
%isomorphism of graded $(\cC \otimes \hat{S})^{op}$-modules
%$$\cE =\Hom _{\cC}(\cC \otimes _{\tau _{\cC}}B\cC,\cC)\to \Hom
%_k(B\cC ,\cC).$$ If we further use the isomorphism of the previous
%lemma we find that the DG $(\cC \otimes \hat{S})^{op}$-module
%$\cE$ is isomorphic to $\cC \otimes \hat{S}$ with the differential
%$d_{\cC}\otimes 1+1\otimes d_{\hat{S}}+\tau _{\cC}$; here $\tau
%_{\cC}$ stands for the left multiplication by the universal
%twisting cochain $\tau _{\cC}\in \Hom _k(B\cC, \cC)=\cC \otimes
%\hat {S}$. We denote this DG $(\cC \otimes \hat{S})^{op}$-module
%$\cC \otimes _{\tau _{\cC}}\hat{S}$.
%\end{remark}

\begin{lemma}\label{E_h-proj} The object $\cE$ considered as a DG $\hat{S}^{op}$-module is
h-projective.
\end{lemma}

\begin{proof} Notice that the stupid filtration $A^{\geq i}$ of the complex $A$ is finite.
Since $A$ is admissible it follows that the differential
$m_1^{\cE}$ preserves the $(\hat{S}^{op})^{gr}$-submodule $A^{\geq
i}\otimes \hat{S}$. Hence the DG $\hat{S}^{op}$-module
$\cE=A\otimes_{\tau_A}\hat{S}$ has a finite filtration by DG
$\hat{S}^{op}$-submodules $A^{\geq i}\otimes \hat{S}$ with
subquotients being free $\hat{S}^{op}$-modules $A^i\otimes
\hat{S}$. Thus the DG $\hat{S}^{op}$-module $\cE$ is h-projective.
\end{proof}

\begin{lemma}\label{Eotimesk=A} The $A_{\infty}$ $\bar{A}^{op}$-module $\cE
\stackrel{\bL}{\otimes}_{\hat{S}}k$ is canonically
quasi-isomorphic to $A$.
\end{lemma}

\begin{proof} By Lemma \ref{E_h-proj} and Remark \ref{CotimesS} we have
$$\cE \stackrel{\bL}{\otimes}_{\hat{S}}k=\cE \otimes_{\hat{S}}k=A \otimes _{\tau _A}\hat{S}\otimes _{\hat{S}}k,$$
and the last $A_{\infty}$ $A ^{op}$-module is isomorphic to $A$
since $\tau _A\in A \otimes \m$, where $\m\subset \hat{S}$ is the
augmentation ideal.
\end{proof}

Now we are ready to define a morphism of pseudo-functors
$$\Psi : h^\prime_{\hat{S}}\to \DEF _-(A).$$
Let $\cR \in \dgart _-$ and $M=(M,\theta )\in
1\text{-}\Hom(\hat{S},\cR)$. We put $$\Psi(M):=\cE
\stackrel{\bL}{\otimes }_{\hat{S}}M\in
D_{\infty}(\bar{A}^{op}_{\cR^{op}}).$$ Notice that the structure
isomorphism $\theta :M\stackrel{\bL}{\otimes }_{\cR}k\to k$
defines an isomorphism
$$\bL i^*\Psi(M)=\Psi(M)\stackrel{\bL}{\otimes}_{\cR}k=\cE
\stackrel{\bL}{\otimes }_{\hat{S}}M\stackrel{\bL}{\otimes
}_{\cR}k\stackrel{\sim}{\to} \cE \stackrel{\bL}{\otimes
}_{\hat{S}}k,$$ and the last term is canonically quasi-isomorphic
to $A$ as an $A_{\infty}$ $\bar{A} ^{op}$-module by Lemma
\ref{Eotimesk=A}. Hence $\Psi (M)$ is indeed an object in the
groupoid $\DEF _{\cR}(A)$.

If $\delta :M\to N$ is a 2-morphism, where $M,N \in
1\text{-}\Hom(\hat{S}, \cR)$, then $\Psi (\delta):\Psi (M)\to \Psi
(N)$ is a morphism of objects in the groupoid $\DEF _{\cR}(A)$.
Thus $\Psi$ is indeed a morphism of pseudo-functors.

\begin{theo}\label{pro-repr_special} The morphism $\Psi : h^\prime_{\hat{S}}\to \DEF _-(A)$
is an equivalence.
\end{theo}

\begin{proof}
It remains to show that for each $\cR \in \dgart _-$ the induced
functor
$$\Psi :1\text{-}\Hom (\hat{S},\cR)\to \DEF _{\cR}(A)$$
is an equivalence of groupoids.

We fix $\cR$.

\noindent{\bf Surjective on isomorphism classes.}

Let $(S,\sigma )\in \DEF _{\cR}(A)$. By Corollary \ref{def^h}
there exists an element $\alpha\in \cM\cC_{\cR}(A)$ such that the
pair $(T,\id)$, where $T=A\otimes_{\alpha} \cR$, defines an object
of $\DEF_{\cR}(A)$ and $(T,\id)$ is isomorphic to $(S,\sigma )$ in
$\DEF _{\cR}(A)$. The element $\alpha $ corresponds to a (unique)
admissible twisting cochain $\tau:\cR ^*\to A$, which in turn
corresponds to a homomorphism of DG coalgebras $g_{\tau}:\cR ^*
\to B\bar{A}$ (Proposition \ref{tw_coch}). By dualizing we obtain
a homomorphism of DG algebras $g^*_{\tau}:\hat{S}\to \cR$ and
hence the corresponding object $M_{\alpha}=({}_{\hat{S}}\cR
_{\cR},\id )\in 1\text{-}\Hom (\hat{S},\cR).$

\begin{lemma}\label{surjective_special} The object $\Psi (M_{\alpha})\in \Def _{\cR}(A)$ is isomorphic to
$(T,\id)$.
\end{lemma}

\begin{proof}
By Remark \ref{CotimesS}
$$\cE
=A \otimes _{\tau _A}\hat{S},$$ and hence by Lemma \ref{E_h-proj}
$$\Psi (M_{\alpha})=(A \otimes _{\tau _A}\hat{S})\otimes
_{g^*_{\tau}}\cR.$$ Notice that the image of $\tau_A$ under the
map $$\one_A\otimes g_{\tau}^*:A \otimes \hat{S}\to A \otimes
\cR$$ coincides with $\tau$. Thus $\Psi (M_{\alpha})=T$.
 \end{proof}

\noindent{\bf Full and faithful.}

Let us define a functor $\Pi :\DEF _{\cR}(A)\to 1\text{-}\Hom
(\hat{S},\cR)$ as follows: for $S=(S, \sigma )\in \DEF _{\cR}(A)$
we put
$$\Pi (S):=\Hom _{\bar{A} ^{op}}(\cE ,S)\in D(\hat{S}^{op}\otimes \cR).$$
We claim that $\Pi (S)$ is  an object in
$1\text{-}\Hom(\hat{S},\cR)$, i. e. it is quasi-isomorphic to
$\cR$ as a DG $\cR ^{op}$-module and the isomorphism $\sigma$
defines an isomorphism $\Pi (S)\stackrel{\bL}{\otimes
}_{\cR}k\stackrel {\sim}{\to}k$.

Indeed, again by Corollary \ref{def^h} we may and will assume that
$(S,\sigma)=(T,\id)$, where $T=A \otimes_{\alpha} \cR$, $\alpha\in
\cM\cC_{\mR}(A)$. We have
$$\Pi (T)=\Hom _{A ^{op}}(\cE ,A \otimes _{\alpha} \cR)=\Hom _{A ^{op}}( \Hom _{A}(k,A ),A)\otimes _{\alpha}\cR.$$ Since
the $A_{\infty}$-algebra $A$ satisfies the condition (*) the last
term as a DG $\cR ^{op}$-module is canonically quasi-isomorphic to
$k\otimes \cR=\cR.$ Thus we have a canonical isomorphism $\Pi
(S)\stackrel{\bL}{\otimes }_{\cR}k\stackrel {\sim}{\to}k$.

Note that the functors $\Psi$ and $\Pi$ are adjoint:
$$\Hom _{\DEF _{\cR}(A)}(\Psi (M),S)=\Hom _{1\text{-}\Hom
(\hat{S},\cR)}(M,\Pi(S)).$$

Now let us consider $\Psi$ and $\Pi$ as functors simply between
the derived categories $D(\hat{S}\otimes \cR ^{op})$ and
$D_{\infty}(\bar{A} ^{op}_{\cR ^{op}}$. (They remain adjoint).
Denote by $\langle k\rangle\subset D(\hat{S}\otimes \cR ^{op})$
and $\langle A \rangle \subset D_{\infty}(\bar{A} ^{op}_{\cR
^{op}}$ the triangulated envelopes of the DG module $k$ and the
$A_{\infty}$ $\bar{A}^{op}_{\cR^{op}}$-module $A$ respectively.
Let $(S,\sigma)\in \DEF _{\cR}(A)$. By Corollary \ref{def^h} we
may and will assume that $(S,\sigma)=(T,\id)$, where $T=A
\otimes_{\alpha} \cR$, $\alpha\in \cM\cC_{\mR}(A)$. Hence $S\in
\langle A\rangle$. Choose $(M, \theta )\in 1\text{-}\Hom
(\hat{S},\cR)$. Since the DG algebra $\hat{S}\otimes \cR ^{op}$ is
local and complete by Lemma \ref{D_f=<k>} we have $M\in \langle
k\rangle$. Therefore it suffices to prove the following lemma.

\begin{lemma}\label{psi=pi^-1} The functors $\Psi$ and $\Pi$ induce mutually inverse
equivalences of triangulated categories $\langle k\rangle$ and
$\langle A \rangle$.
\end{lemma}

\begin{proof} It suffices to prove that the adjunction maps
$k\to \Pi \Psi(k)$ and $\Psi \Pi(A)\to A$ are isomorphisms.

We have $\Pi \Psi (k)=\Hom _{\bar{A}^{op}}(\cE ,\cE
\stackrel{\bL}{\otimes}_{\hat{S}}k)=\Hom _{A ^{op}}(\cE ,A)$
(Lemma \ref{Eotimesk=A}). Hence $k\to \Pi \Psi(k)$ is a
quasi-isomorphism because $A$ satisfies property (*).

Vice versa, $\Psi \Pi (A)=\cE
\stackrel{\bL}{\otimes}_{\hat{S}}(\Hom _{\bar{A} ^{op}}(\cE
,A))=\cE \stackrel{\bL}{\otimes}_{\hat{S}}k$, since $A$ satisfies
property (*). But $\cE \stackrel{\bL}{\otimes} _{\hat{S}}k=A$ by
Lemma \ref{Eotimesk=A}.

This proves the lemma.
\end{proof}
Theorem is proved.
\end{proof}

\subsection{Explicit equivalence $\DEF_-(E)\cong h_{\hat{S}}'$}

Let $E\in \cA^{op}\text{-mod}$. Suppose that the graded algebra
$\Ext(E,E)$ is admissible and finite-dimensional. Let $A$ be a
strictly unital minimal model of the DG algebra $\bR \Hom (E,E)$.
Suppose that $A$ satisfies the condition (*) above. Further, let
$F\to E$ be a quasi-isomorphism with h-projective $F$,
$\cC=\End(F)$ and let $f:A\to \cC$ be a strictly unital
$A_{\infty}$-quasi-isomorphism. By Theorem \ref{pro-repr_special},
the $A_{\infty}$ $\bar{A}^{op}_{\hat{S}^{op}}$-module
$\Hom_{\bar{A}}(k,A)$ is the "universal deformation" of the
$A_{\infty}$ $\bar{A}^{op}$-module $A$. It follows from the
equivalence $\DEF_-(A)\cong \DEF_-(\cC)$ (Corollary
\ref{def(A)=def(B)}) that that the $(\cC\otimes
\hat{S})^{op}$-module
$$\Hom_{\bar{A}}(k,\cC)=\cC\otimes_{f^*(\tau_A)}\hat{S}$$ is the
"universal deformation" of the DG $\cC^{op}$-module $\cC$.

Put
$$\cF=\Hom_{\bar{A}}(k,\cC)\otimes_{\cC} F=(\cC\otimes_{f^*(\tau_A)} \hat{S})\otimes_{\cC} F.$$
Then $\cF$ is a DG $\cA_{\hat{S}}^{op}$-module. We claim that it
is a "universal deformation" of the DG module $E$. More precisely,
we get the following

\begin{cor}\label{explicit} Let $E$ and $\cF$ be as above. Then the functors
$\Phi_{\cR}:D(\hat{S}\otimes \cR^{op})\to D(\cA_{\cR}^{op})$,
$$\Phi_{\cR}(M)=\cF\stackrel{\bL}{\otimes}_{\hat{S}}M,$$ induce the
equivalence of pseudo-functors $$\Phi:h_{\hat{S}}'\to \DEF_-(E)$$
from $\dgart_-$ to {\bf Gpd}.
\end{cor}

\begin{proof} Indeed, the morphism $\Phi:h_{\hat{S}}'\to
\DEF_-(E)$ is isomorphic to the composition of the equivalence
$$\Psi:h_{\hat{S}}'\to \DEF_-(A)$$ from Theorem
\ref{pro-repr_special}, the equivalence $$\DEF_-(A)\cong
\DEF_-(\cC)$$ from Corollary \ref{def(A)=def(B)}, and the
equivalence
$$\DEF_-(\Sigma):\DEF_-(\cC)\to \DEF_-(E)$$ from the proof of
Proposition \ref{def(B)=def(E)}.
\end{proof}

\section{Classical pro-representability}

Recall that for a small groupoid $\cM$ one denotes by $\pi _0(\cM)$
the set of isomorphism classes of objects in $\cM$.

All our deformation functors have values in the 2-category of
groupoids ${\bf Gpd}$. We may compose those pseudo-functors with
$\pi _0$ to obtain functors with values in the category $\Set$ of
sets. Classically pro-representability theorems are statements
about these compositions. Out pro-representability Theorems
\ref{pro-repr1}, \ref{pro-repr2}, \ref{pro-repr3}, \ref{pro-repr4}
have some "classical" implications which we discuss next.

\begin{defi}\label{cat_alg} Denote by $\alg$ and $\calg$ the full subcategories of
the category $\adgalg$ (I, Section 2) consisting of local (!)
augmented algebras (resp. local commutative augmented algebras)
concentrated in degree zero. That is we consider the categories of
usual local augmented (resp. commutative local augmented)
algebras. Then we have the full subcategories $\art \subset \alg$
and $\cart \subset \calg$ of (local augmented) artinian (resp.
commutative artinian) algebras (I, Definitions 2.1-2.3). Note that
for $\cB ,\cC \in \alg$ the group of units of $\cC$ acts by
conjugation on the set $\Hom (\cB ,\cC)$. We call this the adjoint
action. The orbits of this action define an equivalence relation
on $\Hom (\cB ,\cC)$ and we denote by $\alg /\ad$ the
corresponding quotient category, where $\Hom _{\alg /\ad}(\cB
,\cC)$ is the set of equivalence classes. Let
$$q:\alg \to \alg /\ad$$
be the quotient functor. We obtain the corresponding full
subcategory $\art /\ad \subset \alg /\ad$.
\end{defi}

\begin{remark}\label{remark_on_cat_alg} Note that if $\cB ,\cC \in \alg$ and $\cC$ is
commutative then the adjoint action on $\Hom (\cB ,\cC)$ is trivial.
\end{remark}

Recall the pseudo-functor $\cF :\adgalg \to 2\text{-}\adgalg$ from
Section \ref{2-adgalg_codef}. We denote also by $\cF$ its
restriction to the full subcategory $\alg$. Since the functor $q$
and the pseudo-functor $\cF$ are the identity on objects we will
write $\cB$ instead of $q(\cB)$ or $\cF (\cB)$ for $\cB \in \alg$.

Fix $\cB \in \alg$. We consider two functors from $\alg$ to $\Set$
which are defined by $\cB$: $h_{\cB}\cdot q$ and $\pi _0 \cdot
h_{\cB} \cdot \cF$. Namely, for $\cC \in \alg$:
$$h_{\cB}\cdot q(\cC)=\Hom _{\alg /\ad}(\cB ,\cC),$$
$$\pi _0 \cdot h_{\cB} \cdot \cF (\cC)=\pi _0(1\text{-}\Hom
_{2\text{-}\adgalg}(\cB ,\cC ).$$

\begin{lemma}\label{pro-repr_vs_pro-repr_cl} For any $\cB \in \alg$ the above functors $h_{\cB}\cdot
q$ and  $\pi _0 \cdot h_{\cB} \cdot \cF$ from $\alg $ to $\Set$
are isomorphic.
\end{lemma}

\begin{proof} This is proved in Lemma \ref{2-adgalg_vs_alg} a), b).
\end{proof}

\begin{cor}\label{pro-repr_cl_comm} For any $\cB \in \alg $ the functors $h_{\cB}$ and $\pi
_0 \cdot h _{\cB}\cdot \cF $ from $\calg$ to $\Set$ are
isomorphic.
\end{cor}

\begin{proof} This follows from Lemma \ref{pro-repr_vs_pro-repr_cl} and Remark \ref{remark_on_cat_alg}.
\end{proof}

\begin{defi}\label{Koszul} Let $A$ be an augmented $A_{\infty}\text{-}$algebra. We call $A$
Koszul if the DG algebra $\hat{S}:=(B\bar{A})^*$ is
quasi-isomorphic to $H^0(\hat{S})$.
\end{defi}

Note that the augmented $A_{\infty}\text{-}$algebras coming from
Examples \ref{k-points}, \ref{NJac}, \ref{NGr} are formal and
quadratic Koszul, hence Koszul in our sence.

\begin{lemma}\label{invar_h_B} Let $\phi :\cB \to \cC$ be a quasi-isomorphism of
augmented DG algebras. Then it induces a morphism $\phi
_*:h_{\cC}\to h_{\cB}$ of pseudo-functors from $2\text{-}\adgalg$
to ${\bf Gpd}$. This morphism is an equivalence.
\end{lemma}

\begin{proof} Indeed, for $\cE \in 2\text{-}\adgalg$ and $M\in
1\text {-}\Hom (\cC ,\cE)$ denote by $\phi _*M\in 1\text {-}\Hom
(\cB ,\cE)$ the DG $\cB \otimes \cE ^{op}$-module obtained from
$M$ by restriction of scalars. This functor $\phi _*$ defines an
equivalence of derived categories
$$\phi _*:D(\cC \otimes \cE ^{op})\to D(\cB \otimes \cE ^{op})$$
since $\phi$ is a quasi-isomorphism. Hence it defines an equivalence
of groupoids
$$\phi _*:1\text {-}\Hom (\cC ,\cE)\to 1\text {-}\Hom (\cB ,\cE).$$
\end{proof}

\begin{theo}\label{pro-repr_cl1} Let $\cC$ be a DG algebra such
that the strictly unital minimal model $A$ of $\cC$ (Definition
\ref{admissible}) is a Koszul $A_{\infty}$-algebra. Put
$\hat{S}=(B\bar{A})^*$. Then

a) there exists an isomorphism of functors from $\art$ to $\Set$
$$h_{H^0(\hat{S})}\cdot q\simeq \pi _0\cdot \coDef _0(\cC);$$

b) there exists an isomorphism of functors from $\cart $ to $\Set$
$$h_{H^0(\hat{S})}\simeq \pi _0 \cdot \coDef _{\cl}(\cC).$$
\end{theo}

\begin{proof} a) Note that the DG algebra $\hat{S}$ is concentrated
in nonpositive degrees. hence we have a natural homomorphism of
augmented DG algebras $\hat{S}\to H^0(\hat{S})$ which is a
quasi-isomorphism. Hence by Lemma \ref{invar_h_B} the
pseudo-functors
$$h_{\hat{S}}, h_{H^0(\hat{S})}:2\text{-}\adgalg \to {\bf Gpd}$$
are equivalent. Notice that $\hat{S}$ is a local algebra and the
homomorphism $\hat{S}\to H^0(\hat{S})$ is surjective. Hence the
algebra $H^0(\hat{S})$ is also local.

By Theorem \ref{pro-repr1} we have an equivalence of
pseudo-functors
$$\coDEF _0(\cC)\simeq h_{\hat{S}}:2\text{-}\art \to {\bf Gpd}.$$
Thus $\coDEF _0(\cC)\simeq h_{H^0(\hat{S})}$. By Proposition
\ref{coDEF}
$$\coDEF _-(\cC)\cdot \cF \simeq \coDef _-(\cC).$$
Therefore
$$\coDef _0(\cC)\simeq h_{H^0(\hat{S})}\cdot \cF :\art \to {\bf
Gpd}.$$ Finally, by Lemma \ref{pro-repr_vs_pro-repr_cl}
$$\pi _0\cdot \coDef _0(\cC)\simeq h_{H^0(\hat{S})}\cdot q:\art \to
\Set.$$ This proves a).

b) This follows from a) and Remark \ref{remark_on_cat_alg}
\end{proof}

\begin{remark}\label{pro-repr_all} Under the assumptions of
Theorem \ref{pro-repr_cl1} the same conclusion holds for
pseudo-functors $\coDef ^{\h}(\cC),$ $\Def (\cC),$ $\Def
^{\h}(\cC)$ instead of $\coDef (\cC)$. Indeed, by Lemmas
\ref{def=def^h}, \ref{codef=codef^h} and Theorem \ref{def=codef}
there are equivalences of pseudo-functors
$$\coDef _-(\cC)\simeq \coDef
^{\h}_-(\cC)\simeq \DEF _-(\cC)\simeq \Def _-^{\h}(\cC).$$
\end{remark}

\begin{theo}\label{pro-repr_cl2} Let $E\in \cA ^{op}\text{-mod}$. Assume that $E$ is
quasi-isomorphic to a bounded below $F\in \cA ^{op}\text{-mod}$
which is h-projective or h-injective. Also assume that the graded
algebra $\Ext(E)$ is admissible and finite-dimensional, and the
strictly unital minimal model $A$ of the DG algebra $\End(F)$ is
Koszul. Put $\hat{S}=(B\bar{A})^*$. Then

a) there exists an isomorphism of functors from $\art $ to $\Set$
$$h_{H^0(\hat{S})}\cdot q\simeq \pi _0\cdot \coDef _0(E);$$

b) there exists an isomorphism of functors from $\cart $ to $\Set$
$$h_{H^0(\hat{S})}\simeq \pi _0 \cdot \coDef _{\cl}(E).$$
\end{theo}

\begin{proof} By I, Proposition 11.16 the
pseudo-functors $\coDef _-(E)$ and $\coDef _-(\End(F))$ are
equivalent. So the theorem follows from Theorem
\ref{pro-repr_cl1}.
\end{proof}

\begin{theo}\label{pro-repr_cl3} Let $E\in \cA ^{op}\text{-mod}$. Assume that
the graded algebra $\Ext(E)$ is admissible and finite-dimensional,
and the strictly unital minimal model $A$ of the DG algebra $\bR
\Hom(E,E)$ is Koszul. Put $\hat{S}=(B\bar{A} )^*$. Then

a) there exists an isomorphism of functors from $\art $ to $\Set$
$$h_{H^0(\hat{S})}\cdot q\simeq \pi _0\cdot \Def _0(E);$$

b) there exists an isomorphism of functors from $\cart $ to $\Set$
$$h_{H^0(\hat{S})}\simeq \pi _0 \cdot \Def _{\cl}(E).$$
\end{theo}

\begin{proof} By I, Proposition 11.16 the
pseudo-functors $\Def _-(E)$ and $\Def _-(\bR \Hom(E,E))$ are
equivalent. So the theorem follows from Theorem \ref{pro-repr_cl1}
and Remark \ref{pro-repr_all}.
\end{proof}

If, in addition, the $A_{\infty}$-algebra $A$ in the above Theorem
satisfies condition (*) (Definition \ref{cond_*}), then the
equivalences a), b) can be made explicit. Namely, we get the
following

\begin{cor}\label{explicit1} Let $E$, $\hat{S}$ be as in Theorem \ref{pro-repr_cl3} and $\cF$
be as in Corollary \ref{explicit}. Then the equivalence
$h_{H^0(\hat{S})}\cdot q\to \Def_0(E)$ of functors from $\art$ to
$\Set$, and the equivalence $h_{H^0(\hat{S})}\to \Def_0(E)$ of
functors from $\cart$ to $\Set$ are induced by the functors
$\Phi_{\cR}:D(H^0(\hat{S})\otimes \cR^{op})\to D(\cA_{\cR}^{op})$,
$$\Phi_{\cR}(M)=\cF\stackrel{\bL}{\otimes}_{\hat{S}}M.$$
\end{cor}

\begin{proof} This follows from Corollary \ref{explicit}
\end{proof}

\end{document}